\documentclass[a4paper,12pt]{article}
\usepackage{amssymb,latexsym,amsmath,amsthm,amscd,graphpap}
\usepackage{amsmath}
\usepackage{amsfonts}
\usepackage{amssymb}

\newtheorem{theorem}{Theorem}[section]
\newtheorem{prop}[theorem]{Proposition}
\newtheorem{defi}[theorem]{Definition}

\newtheorem{lemma}[theorem]{Lemma}

\newtheorem{coro}[theorem]{Corollary}

\newenvironment{demo}{ \noindent \emph{\textbf{Proof:}}}{\hfill$\square$\\
\vspace{0.4cm}}
\newenvironment{rem}{ \noindent \emph{\textbf{Remark :}}}{  \bigskip }
\newenvironment{rems}{ \noindent \emph{\textbf{Remarks :}}} {
\bigskip }

\newcommand{\Rm}{\mathbb{R}}

\newcommand{\Nm}{\mathbb{N}}

\newcommand{\Com}{\mathbb{C}}

\newcommand{\Gm}{\mathfrak{G}}

\newcommand{\Cm}{\mathcal{C}}
\newcommand{\Nc}{\mathcal{N}}

\newcommand{\Om}{\mathcal{O}}
\newcommand{\Em}{\mathcal{E}}

\newcommand{\Um}{\mathcal{U}}
\newcommand{\Vm}{\mathcal{V}}

\newcommand{\Zm}{\mathcal{Z}}
\newcommand{\Zr}{\mathbb{Z}}

\newcommand{\no}{n$^{\text{o}}$}

\renewcommand{\Re}{\mathrm{Re}}
\renewcommand{\Im}{\mathrm{Im}}
\renewcommand{\dim}{\mathrm{dim}}
\newcommand{\codim}{\mathrm{codim}}
\newcommand{\ind}{\mathrm{ind}}
\newcommand{\lcro}{[\![}
\newcommand{\rcro}{]\!]}

\numberwithin{equation}{section}
\def\HB {\hfill\break}

\begin{document}

\title{\bf Generic Morse-Smale property for the parabolic equation on
the circle}

\author{Romain JOLY and Genevi\`eve RAUGEL }

\date{}

\maketitle
\vspace{1cm}

\begin{abstract}
In this paper, we show that, for scalar reaction-diffusion
equations $u_t=u_{xx}+f(x,u,u_x)$ on the circle $S^1$,
the Morse-Smale property is generic with respect
to the non-linearity $f$. In \cite{CR}, Czaja and Rocha have proved that
any connecting orbit, which connects two hyperbolic periodic orbits,
  is transverse and that there does not exist any homoclinic orbit,
  connecting a hyperbolic periodic orbit to itself. In \cite{JR},
we have shown that,
generically with respect to
the non-linearity $f$, all the equilibria and periodic
orbits are hyperbolic.
Here we complete these results by
showing that any connecting orbit between two hyperbolic equilibria
with distinct Morse indices or between a hyperbolic equilibrium and
a hyperbolic periodic orbit is automatically transverse. We also
show that, generically with respect to $f$, there does not exist
any connection between equilibria with the same Morse index. The above
properties, together with the existence of a compact global attractor and
the Poincar\'e-Bendixson property,
allow us to deduce that, generically with respect to $f$,
the non-wandering set consists in a finite number of
hyperbolic equilibria and periodic orbits .
The main tools in the proofs include the lap number property, exponential
dichotomies and the Sard-Smale theorem. The proofs also require a
careful analysis of the asymptotic behavior of solutions of the
linearized equations along the connecting orbits.

\

\end{abstract}

\vspace{1cm}
\section{Introduction}

In the study of the dynamics of flows or semi-flows generated by
systems of ordinary differential or partial differential
equations arising in physics or biology, global stability (also
called structural stability) is a very important property.  Indeed,
often, one only knows approximate values of the various coefficients in the
equations; or else, in order to numerically determine the solutions of the
equations, one introduces a space and (or) time discretized system.
Therefore, one often studies a system, which is an approximation of
the original dynamical system. If the dynamics of the original system
are globally stable, then the qualitative
global behaviour of the solutions remains unchanged under small perturbations
of the system and the knowledge of the dynamics of this approximate
system is sufficient in practice. Unfortunately, in general dynamical
systems bifurcation phenomena can take place and thus
drastic changes in the dynamics can arise. However, one may hope that
such phenomena almost never happen in the considered class of
dynamical systems or that the systems, which are robust, are dense or
generic in the considered class (see Section 1.2 below for the
definition of genericity).
\vskip 3mm

In the 1960's and 1970's, such structural stability problems have been
extensively studied in the frame of vector fields (and also iterates of
diffeomorphisms) on compact smooth manifolds $\mathcal{M}^n$ of finite
dimension $n \geq 1$. In this context (see \cite{Smale60}), Smale introduced
the notion of Morse-Smale dynamical systems, that is, systems
for which the non-wandering set consists only in a
finite number of hyperbolic equilibria and hyperbolic periodic orbits
and the intersections of the stable and unstable manifolds
of equilibria and periodic orbits are all transversal.
Palis and Smale have shown that Morse-Smale vector fields on compact
manifolds $\mathcal{M}^n$ are structurally stable (\cite{Pal69},
\cite{PaSm}). Moreover, the class of Morse-Smale vector fields on a
compact manifold $\mathcal{M}^2$ of dimension $2$ is generic in the
class of all $C^1$-vector fields (\cite{Pei}).
In the case of the sphere $S^2$, the proof is a consequence of the genericity
of the Kupka-Smale property and of the Poincar\'{e}-Bendixson theorem.
In the same way, one also shows that the Morse-Smale vector fields are
generic in the class of ``dissipative" vector fields on $\Rm^2$.  In
the case of a general two-dimensional manifold (especially in the
non-orientable case), the proof is more delicate.
The Morse-Smale property is also generic in the class of all gradients vector
fields on a Riemannian manifold $\mathcal{M}^n$, $n \geq 1$.  In
the simple case of gradient systems (where the non-wandering set is
reduced to equilibrium points), this genericity property is an
immediate consequence of the genericity of the Kupka-Smale vector
fields (\cite{Kupka,Smale2,Peixoto2}).
On compact manifolds $\mathcal{M}^n$ of finite
dimension $n \geq 3$, the Morse-Smale vector fields are
still plentiful. There had been some hope that
Morse-Smale systems (and hence stable dynamics) could be generic in the class
of general vector fields. But unfortunately, in dimension higher than two, the
Morse-Smale vector fields are no longer dense in the
set of all vector fields.  In particular, transverse homoclinic
orbits connecting a hyperbolic periodic orbit to itself may exist,
giving rise to chaotic behaviour (\cite{Smale65}). And this cannot be
removed by small perturbations.

\vskip 3mm

The above mentioned results give us a hint on what can be expected in
the case of dynamical systems generated by partial differential
equations (PDE's in short). The results available in infinite
dimensions are still rather partial. Like in the case of vector
fields on compact manifolds or on $\Rm^n$, one can define Morse-Smale
dynamical systems.  W. Oliva \cite{Oliva}
(see also \cite{HMO}) has proved that the Morse-Smale dynamical
systems $S(t)$, generated by ``dissipative" parabolic equations or
more generally by ``dissipative in finite time smoothing equations",
are structurally stable, in the sense that the restrictions of the
flows to the compact global attractors are topologically equivalent
under small perturbations.  Already, in 1982, he had proved the
structural stability of Morse-Smale maps, which in turn implies the
stability of Morse-Smale gradient semi-flows generated by PDE's.  We
recall that a dynamical system $S(t)$ on a Banach space, generated by an
evolutionary PDE, is gradient if it admits a strict Lyapunov
functional, which implies that the non-wandering set reduces to the
set of equilibria.

As in the finite-dimensional case, one expected the density or
genericity of the Morse-Smale gradient semi-flows within the set of
gradient semi-flows generated by a given class of PDE's.  Already, in
1985, D. Henry \cite {He85} proved the noteworthy property that the
stable and unstable manifolds of two equilibria of the
reaction-diffusion equation with separated boundary conditions on the
interval $(0,1)$ intersect transversally (see also \cite{Ang86} for
another proof in the case of hyperbolic equilibria).  One of the main
tools in his proof was the decay property of the zero number (also
called Sturm number or lap number; see Section 2.1), which will also
be often used in this paper.  Since, as shown by Zelenyak in \cite{Ze},
the reaction-diffusion equation on the interval $(0,1)$ with separated
boundary conditions is gradient, the transversality property of Henry
implies the first known result of genericity of Morse-Smale systems in
a class of PDE's.
The scalar reaction-diffusion equation, defined on a bounded domain
$\Omega$ of $\Rm^n$, $n \geq 2$, is no longer gradient in general.
However, it is gradient if one considers non-linearities $f(x,u)$, depending
only on $x$ and the values of the function $u$ (and not of the values of its
derivatives). In this class of gradient parabolic equations,
Brunovsk\'{y} and Pol\'{a}\v{c}ik \cite{BrunPola} have shown in 1997
that the Morse-Smale property is generic with respect to the
non-linearity $f(x,u)$.  Later, the genericity of
the gradient Morse-Smale flows in the class of gradient flows
generated by the damped wave equations (with fixed damping) defined on
any bounded domain has been proved by Brunovsk\'{y} and Raugel in
\cite{BrunRaug} (for the case of variable damping, we refer to
\cite{RJo}).

It must be emphasized that, whereas the proof of the structural stability
follows the lines (with some adjustments) of the proof given on compact
finite-dimensional manifolds, the proof of the genericity of Morse-Smale
property requires other approaches. Indeed, perturbing a semi-flow generated
by a PDE in order to make it Morse-Smale has an interest only if one is able to
perform it within the same class of equations. In the case of general vector
fields on finite-dimensional manifolds, one can perturb the vector field in a
local manner with all the freedom one needs. In the case of PDE's,
the perturbed equations must remain in the considered class. Therefore,
the perturbations are constrained.  An analogous problem involving
constrained perturbations has been studied in the finite-dimensional
case by Robbin \cite{Robbin}. An additional problem arises in the
case of PDE's, namely the perturbations could a priori really non-local.
Typically, perturbing the nonlinearity changes the semi-flow in a
large part of the phase space in a way which is hard to understand.

\vskip 3mm

At first glance, the non-density of the Morse-Smale vector fields on
compact manifolds of dimension $n \geq 3$ gave only little hope that
infinite-dimensional Morse-Smale semi-flows are dense in some classes
of non-gradient PDE's.  However, Fiedler and Mallet-Paret (\cite{FMP})
showed in 1989 that the scalar reaction-diffusion equation \eqref{eq}
on $S^1$ satisfies the Poincar\'{e}-Bendixson property (which, as we
recalled earlier, played an important role in the proof of density of
Morse-Smale vector fields in the set of all dissipative vector fields
in $\Rm^2$).  More recently, in 2008, Czaja and Rocha (\cite{CR})
proved that, for the scalar reaction-diffusion equation on $S^1$, the
stable and unstable manifolds of hyperbolic periodic orbits always
intersect transversally.  These results gave us some hope that
Morse-Smale dynamics could be generic for scalar reaction-diffusion
equations on $S^1$ since they are generic for two-dimensional vector
fields.  In 2008, we proved that the equilibria and periodic orbits
are hyperbolic, generically with respect to the non-linearity
(\cite{JR}).  The results of Fiedler, Rocha and Wolfrum  (\cite{FRW})
together with the generic hyperbolicity property of \cite{JR} imply that
the Morse-Smale property is generic in the special class of
reaction-diffusion equations with spatially homogeneous
non-linearities $f(u,u_x)$.

Here, we complete the global qualitative picture of the
scalar reaction-diffusion equations \eqref{eq} on $S^1$ in the case of a
general non-linearity $f(x,u,u_x)$ and conclude the
proof of the genericity of the Morse-Smale systems in this class.
These results indicate a similarity between scalar
reaction-diffusion equations on $S^1$ and two-dimensional vector fields
and take place in a more general correspondence
between parabolic equations and finite-dimensional vector fields in
any space dimension, as noticed in \cite{JR2}.
For scalar parabolic equations on bounded domains $\Omega$ in $\Rm^d$,
$d\geq 2$, the properties of zero number do no longer
hold, the Poincar\'e-Bendixson property fails and the Morse-Smale property is
no longer generic.  But, the genericity of the
Kupka-Smale property still holds
like in the case of vector fields in dimension $n\geq 3$, see 
\cite{BrunJolyRaug}.

\subsection{The parabolic equation on the circle: earlier results}

In this paper, we consider the following scalar reaction-diffusion
equation on $S^1$,
\begin{equation}
\label{eq}
\left\{
\begin{array}{rl}
u_t(x,t)&=u_{xx}(x,t)+f(x,u(x,t),u_x(x,t))~, \quad (x,t)\in S^1\times (0,
+\infty)~,\\
u(x,0)&=u_0(x)~, \quad x\in S^1~,
\end{array}
\right.
\end{equation}
where $f$ belongs to the space $C^2(S^1\times \Rm\times\Rm,\Rm)$ and
$u_0$ is given in the Sobolev space $H^s(S^1)$, with $s\in(3/2,2)$
(so that $H^s(S^1)$ is continuously embedded into $C^{1 +
\alpha}(S^1)$ for $\alpha=s -3/2$).

Eq. \eqref{eq} defines a local dynamical system $S_f(t)$ on
$H^s(S^1)$ (see \cite{He81} and \cite{Polacik-hand}) by setting
$S_f(t)u_0=u(t)$, where $u(t)$ is the solution of \eqref{eq}
(if the dependence in $f$ of the dynamical system is
not important, we simply denote $S(t)$ instead of $S_f(t)$). \HB
In order to obtain a global dynamical system, we impose some additional
conditions on $f$, namely
we assume that there exist a function $k(.) \in\Cm^0(\Rm^+,\Rm^+)$
and constants
$ \varepsilon>0$ and $\kappa>0$ such that
\begin{equation}
\label{condf1}
\begin{split}
\forall R>0,~\forall \xi\in\Rm~, \quad &\sup_{(x,u)\in S^1\times[-R,R]}
|f(x,u,\xi)|\leq k(R)(1+|\xi|^{2-\varepsilon})~,\cr
    \forall |u|\geq \kappa,~\forall x \in S^1~, \quad
& \quad  \quad uf(x,u,0)\leq 0~.
\end{split}
\end{equation}
Then, Eq. \eqref{eq} defines a global dynamical system
$S_f(t)$ in $H^s(S^1)$ (see \cite{Polacik-hand}).  Moreover, $S_f(t)$
  admits a compact global attractor ${\cal A}_f$, that is,
there exists a compact set ${\cal A}_f$ in $H^s(S^1)$ which is invariant
(i.e. $S_f(t) {\cal A}_f= {\cal A}_f$, for any $t\geq 0$) and attracts every
bounded set of $H^s(S^1)$.

The most interesting part of the dynamics of \eqref{eq} is contained in the
attractor ${\cal A}_f$. Our purpose is to describe these dynamics, at least 
for a
dense set of nonlinearities $f$. We introduce the set
$\Gm=\Cm^2(S^1\times \Rm\times\Rm,\Rm)$
endowed with the {\sl Whitney topology}, that is, the topology generated by the
neighbourhoods
\begin{equation}\label{Whitney}
\{g\in\Gm~/~|D^i f(x,u,v)-D^i g(x,u,v)|\leq\delta(u,v),~\forall
i\in\{0,1,2\},~
\forall (x,u,v)\in S^1\times\Rm^2\}~,
\end{equation}
where $f$ is any function in $\Gm$ and $\delta$ is any positive
continuous function (see \cite{GG}). It is well known that $\Gm$ is a Baire
space, which means that any countable intersection of open and dense
sets is dense in $\Gm$ (see \cite{GG} for instance). We say that a set is
{\sl generic} if it countains a countable intersection of open and dense
sets and we say that the parabolic equations on the circle \eqref{eq} 
satisfy a property
generically (with respect to the non-linearity) if this property
holds for any $f$ in a generic subset of $\Gm$. The notion of 
genericity is a common notion
for defining ``large" subsets of Baire spaces, replacing the notion of
``almost everywhere" of $\Rm^d$.

\vskip 3mm

Before describing the properties of \eqref{eq}, we recall a few
basic notions in dynamical systems, for the reader convenience.
For any  $u_0\in H^1(S^1)$, the $\omega$ and
$\alpha$-{\sl limit sets} of $u_0$ are defined
respectively by
\begin{equation*}
\begin{split}
\omega(u_0)=&\{v\in H^1(S^1) \,|\,  \hbox{ there exists a sequence }
    t_n \in \Rm^+, \hbox{ such that }~ \cr
& \hskip 27mm  t_n \mathrel {\mathop \rightarrow_{n \rightarrow 
+\infty}} +\infty~
   \hbox{and}~  S(t_n)u_0 \mathrel {\mathop\rightarrow_{n \rightarrow
   +\infty}}  v\} \cr
\alpha(u_0)=&\{v\in H^1(S^1) \,|\,  \hbox{ there exist a negative
orbit } u(t), \,t \leq 0, \hbox{ with } u(0)=u_0, \cr
& \hbox{ and a sequence }
t_n \in \Rm^+, \hbox{ such that }
  t_n \mathrel {\mathop \rightarrow_{n \rightarrow +\infty}} +\infty~
   \hbox{and}~  u(-t_n) \mathrel {\mathop\rightarrow_{n \rightarrow
   +\infty}}
v\}
\end{split}
\end{equation*}
The {\sl non-wandering set} is the set of points $u_0\in H^1(S^1)$ such that
for any neighbourhood $\Nc$ of $u_0$, $S(t)\Nc\cap\Nc\neq\emptyset$ for
arbitrary large times $t$. In particular, equilibria and periodic 
orbits belongs
to the non-wandering set. \HB
A {\sl critical element} means either
an equilibrium point or a periodic solution of \eqref{eq}. \HB
Let $e \in H^s(S^1)$ be an equilibrium point of \eqref{eq}.
As usual, one introduces the linearized operator $L_e$ around the
equilibrium $e$ (see Section 2, for the precise definition and the
spectral properties of $L_e$). \HB
We say that $e$ is a {\sl hyperbolic} equilibrium point if the
intersection of the spectrum $\sigma(L_e)$ of $L_e$ with the imaginary
axis is empty.  The {\sl Morse index} $i(e)$ is the (finite) number of
eigenvalues of $L_e$ with positive real part (counted with their
multiplicities). \HB
If $e$ is a hyperbolic equilibrium point of \eqref{eq}, there exists a
neighbourhood $U_e$ of $e$ such that the local stable and unstable
sets
\begin{equation*}
\begin{split}
W^s_{loc}(e) \equiv W^s(e,U_e) = &\{u_0 \in H^s(S^1) \, | \, S_f(t)u_0
\in U_e\, , \, \forall t \geq 0\} \cr
W^u_{loc}(e) \equiv W^u(e,U_e) = &\{u_0 \in H^s(S^1) \, | \,
S_f(t)u_0 \hbox{ is well-defined for }t\leq 0 \hbox{ and } \cr
& \hskip 27mm S_f(t)u_0
\in U_e\, , \, \forall t \leq 0\} ~.
\end{split}
\end{equation*}
are (embedded) $C^1$-submanifolds of $H^s(S^1)$ of codimension $i(e)$ and
dimension $i(e)$ respectively.\HB
We also define the global stable and unstable sets
\begin{equation*}
\begin{split}
& W^s(e) = \{u_0 \in H^s(S^1) \, | \, S_f(t)u_0
\mathrel {\mathop\rightarrow_{t \rightarrow
   +\infty}} e\}~, \cr
& W^u(e) = \{u_0 \in H^s(S^1) \, | \, S_f(t)u_0\hbox{ is well-defined 
for }t\leq 0 \hbox{ and }
S_f(t)u_0 \mathrel {\mathop\rightarrow_{t \rightarrow
   - \infty}} e\}~.
   \end{split}
\end{equation*}
Since the parabolic equation \eqref{eq}, as well as the corresponding adjoint
equation, satisfy the backward uniqueness property (see \cite{Bardos-Tartar}),
Theorem 6.1.9 of \cite{He81} implies that $W^s(e)$ and
$W^u(e)$ are injectively immersed $C^1$-manifolds in $H^s(S^1)$ of
codimension $i(e)$ and dimension $i(e)$ respectively (see
also  \cite{ChChHa} and \cite{HJR}). We remark that
$W^u(e) = \cup_{t \geq 0} S_{f}(t)W^u(e,U_e)$ is the union of $C^1$ 
embedded submanifolds
of $H^s(S^1)$ of dimension $i(e)$ (see \cite{CR}).

Let next $\Gamma =\{\gamma (x,t) | t \in [0,p]\}$ be a periodic orbit
of \eqref{eq} of minimal period $p$.  The linearized equation around
$\Gamma$ defines an evolution operator $\Pi(t,0): \varphi_0 \in
H^s(S^1) \rightarrow \Pi(t,0)\varphi_0 =\varphi(t) \in H^s(S^1)$,
where $\varphi(t)$ is the solution of the linearized equation. The
operator $\Pi(p,0)$ is called the {\sl period map} (see Section 2 for
the precise definition of $\Pi(p,0)$ and its spectral properties) . \HB
The periodic orbit $\Gamma$ or the periodic solution
$\gamma(t)$ is {\sl hyperbolic} if the intersection of the spectrum
$\sigma(\Pi(p,0))$ with the unit circle in $\Com$ is reduced to $1$
and $1$ is a simple (isolated) eigenvalue. The {\sl Morse index} $i(\Gamma)$
is the (finite) number of eigenvalues of $\Pi(p,0)$ of modulus strictly
larger than $1$ (counted with their multiplicities). \HB
   By \cite[Theorem 14.2
and Remark 14.3]{Ruelle}) or \cite{HaRa09} (see also
\cite{HJR}), if $\gamma(t)$ is a hyperbolic periodic orbit, there
exists a small
neighbourhood $U_{\Gamma}$ of $\Gamma$ in $H^s(S^1)$ such that
\begin{equation}
\label{stunsper}
\begin{split}
&W^s_{loc}(\Gamma) \equiv W^s(\Gamma,U_{\Gamma}) = \{u_0 \in H^s(S^1)
\, | \, S_f(t)u_0
\in U_{\Gamma}\, , \, \forall t \geq 0\} \cr
&W^u_{loc}(\Gamma) \equiv W^u(\Gamma,U_{\Gamma}) = \{u_0 \in H^s(S^1)
\, | \, S_f(t)u_0
\in U_{\Gamma}\, , \, \forall t \leq 0\} \cr
\end{split}
\end{equation}
are (embedded) $C^1$-submanifolds of $H^s(S^1)$ of codimension $i(\Gamma)$ and
dimension $i(\Gamma) +1$ respectively.  We also define the global
stable and unstable sets
\begin{equation*}
\begin{split}
&W^s(\Gamma) = \{u_0 \in H^s(S^1) \, | \, S_f(t)u_0
\mathrel {\mathop\rightarrow_{t \rightarrow
   +\infty}}  \Gamma\} ~, \cr
& W^u(\Gamma) = \{u_0 \in H^s(S^1) \, | \, S_f(t)u_0
\hbox{ is well-defined for }t\leq 0 \hbox{ and }
S_f(t)u_0 \mathrel {\mathop\rightarrow_{t \rightarrow
   -\infty}} \Gamma\} ~.
\end{split}
\end{equation*}
Again, Theorem 6.1.9 of \cite{He81} implies that $W^s(\Gamma)$ and
$W^u(\Gamma)$ are injectively immersed $C^1$-manifolds in $H^s(S^1)$ of
codimension $i(\Gamma)$ and dimension $i(\Gamma)+1$ respectively (see
also \cite{ChChHa} and \cite{HJR}). We also notice
(see Lemma 6.1 of \cite{CR} and also \cite{HaRa09} or
\cite{HJR}) that
$W^u(\Gamma) = \cup_{t \geq 0} S_f(t) W^u_{loc}(\Gamma)$
is a union of embedded submanifolds of $H^s(S^1)$ of dimension
$i(\Gamma)+1$.

Let $q^{\pm}$ be two hyperbolic critical elements. We say that
$W^u(q^-)$ and $W^s_{loc}(q^+)$ {\sl intersect transversally} (or are
transverse) and we denote it by
$$
W^u(q^-) \pitchfork W^s_{loc}(q^+)~,
$$
if, at each intersection point $u_0 \in W^u(q^-) \cap W^s_{loc}(q^+)$,
$T_{u_0}W^u(q^-)$ contains
a closed complement of $T_{u_0}W^s_{loc}(q^+)$ in $H^s(S^1)$. By
convention, two manifolds which do not intersect are always transverse.

\vskip 3mm

We now describe all the known properties of the dynamics of \eqref{eq}.
First, we mention that in the particular case where 
$f(x,u,u_x)=f(x,u)$ does not
depend on the values of the derivative $u_x$, the dynamical system 
$S_{f}(t)$ generated by
\eqref{eq} is gradient, that is, admits a strict Lyapunov functional.
In particular, it has no periodic orbits and the non-wandering set is
reduced to equilibria. As a direct consequence of \cite{BrunPola}, the
Morse-Smale property is generic with respect to the non-linearity
$f(x,u)$.

In the general case where $f(x,u,u_x)$ depends on the three
variables, all the two-dimensional dynamics
  can be realized on  locally (non-stable) invariant manifolds of the flow
  $S_{f}(t)$ of \eqref{eq} (see \cite{SandFied}) and thus  periodic 
orbits can exist
  (\cite{AngenFied}). Hence, the dynamics can be more complicated. However, 
like for
    vector-fields in $\Rm^2$, the remarkable Poincar\'{e}-Bendixson
    property holds (\cite{FMP}).

\begin{theorem}\label{Poincare-B} {\bf
Fiedler-Mallet-Paret (1989)} {\bf (Poincar\'{e}-Bendixson property)}\\
For any $u_0 \in H^s(S^1)$, the $\omega-$limit set $\omega(u_0)$ of $u_0$,
   satisfies exactly one of the following possibilities.\\
i) Either $\omega(u_0)$ consists of a single periodic orbit,\\
ii) or the $\alpha-$ and $\omega-$limit sets of any $v \in \omega(u_0)$
consist only of equilibrium points.\\
\end{theorem}

Theorem \ref{Poincare-B} is the first step towards showing that the
non-wandering set generically reduces to a finite number of equilibria
and periodic orbits.  One of the main ingredients of the proof of
Theorem \ref{Poincare-B} is the {\sl Sturm property} (also called zero
number or lap number property).  More precisely, for any $\varphi \in
C^1(S^1)$, we define the zero number $z(\varphi)$ as the (even) number
of strict sign changes of $\varphi$.  If $v(x,t)$ is the solution of a
scalar linear parabolic equation on a time interval $I= (0, \tau)$,
then, $z(v(\cdot,t))$ is finite, for any $t \in I$, and nonincreasing
in $t$.  Moreover, $z(v(\cdot,t))$ drops at $t=t_0$, if and only if
$v(\cdot, t_0)$ has a multiple zero (the detailed statement of these
properties is recalled in Section 2.1).  Sturm property is very
specific to the parabolic equation in space dimension one.

The second major step on the way to the proof of the genericity of
the Morse-Smale property has been the paper of Czaja and Rocha
\cite{CR}. Inspired by the transversality results of Fusco and Oliva
(\cite{FuOliva}) for special classes of vector fields on $\Rm^n$,
Czaja and Rocha have proved the following
fundamental and nice transversality properties.
\begin{theorem}\label{CzaRocha} {\bf Czaja-Rocha (2008)} \HB
1) There does not exist  any solution $u(t)$ of
\eqref{eq} converging to a same hyperbolic periodic orbit $\Gamma$, as
$t$ goes to $\pm \infty$.\HB
2) Let $\Gamma^\pm$ be two hyperbolic periodic orbits. Then,
$$
W^u(\Gamma^-) \pitchfork W^s_{loc}(\Gamma^+)~.
$$
Moreover, if the intersection $W^u(\Gamma^-)\cap
W^s_{loc}(\Gamma^+)$ is not empty, then $i(\Gamma^-) > i(\Gamma^+)$.
\end{theorem}

Among other arguments, the proof of Theorem \ref{CzaRocha} (\cite{CR}) uses the
decay properties of the zero number as well as
the filtrations of the
phase space with respect to the asymptotic behaviour of the solutions
of the linearized equation around an orbit connecting
two hyperbolic periodic orbits like in \cite{ChChHa}.

The above results hold under the assumption of hyperbolicity of the
periodic orbits. To complete Theorem \ref{CzaRocha}, we have proved
in \cite{JR} that, generically with respect to the non-linearity $f$, 
the periodic
orbits are all hyperbolic, which means that the above hyperbolicity
assumption is not so restrictive.

\begin{theorem} \label{th-JR} {\bf Joly-Raugel (2008)}
There exists a generic subset $\mathcal{O}_h$ of $\Gm$ such that, for
any $f \in \mathcal{O}_h$, all the equilibria and the periodic
solutions of \eqref{eq} are hyperbolic.
\end{theorem}

Besides the Sard-Smale theorem (recalled in Appendix A), one
of the main ingredients of Theorem \ref{th-JR} is again the zero number
property of the difference of two solutions of \eqref{eq} or of the 
solutions of
the linearized equations around equilibria or periodic solutions.

\subsection{Main new results}

In this paper, we prove that, generically with respect to $f$, the
semi-flow $S_f(t)$ generated by Eq. \eqref{eq} on $S^1$ is Morse-Smale. To this
end, we first complete the automatic transversality results of Czaja
and Rocha as follows.

\begin{theorem}\label{th1}
{\bf Automatic transversality results}\\
1) If $e_-$ and $e_+$ are two hyperbolic equilibrium points of
\eqref{eq} with different Morse indices,
then the unstable manifold $W^u(e_-)$ transversally intersects the
local stable manifold $W^s_{loc}(e_+)$.\\
2) If $\Gamma$ (respectively $e$) is a hyperbolic periodic orbit
(respectively a hyperbolic equilibrium point) of
\eqref{eq}, then
$$
W^u(e) \pitchfork W^s_{loc}(\Gamma)~,
$$
and
$$
W^u(\Gamma) \pitchfork W^s_{loc}(e)~.
$$
\end{theorem}

The proof of Theorem \ref{th1} is similar to the one of Theorem
\ref{CzaRocha}. However, the proof of the automatic transversality of the
connecting orbits between two hyperbolic equilibria of different
Morse indices requires a tricky argument, in addition to
those of \cite{CR}. We also emphasize that, even if most of the
ideas in our proof are basically similar to the ones of
\cite{CR}, they are used in a different way.
In fact, as in \cite{CR}, the basic tools are the same as the ones of 
\cite{He85},
\cite{Angenent}, \cite{FuOliva}, namely a careful analysis of
the asymptotic behavior of solutions converging to an equilibrium or a periodic
orbit (Appendix \ref{AppendixD}) combined with a systematic
application of  the Sturm properties (recalled in Theorem
\ref{th-zeros}).
\vskip 1mm

Theorems \ref{CzaRocha} and \ref{th1} show that
  a non-transverse connecting orbit arises only as an orbit connecting two
equilibrium points with same Morse index. We call such an orbit a
{\sl homoindexed orbit}. Since every two-dimensional flow can be
realized in a (locally) invariant manifold of the semi-flow
of a parabolic equation on $S^1$ (\cite{SandFied}), we know that
homoindexed orbits, in particular homoclinic orbits, may occur in the flow of
\eqref{eq}. However, as we prove here, such connecting orbits can be broken
generically in $f$.

\begin{theorem}\label{th2}
{\bf Generic non-existence of homoindexed connecting orbits}\\
   There exists a generic subset $ \mathcal{O}_{M} \subset
\mathcal{O}_h$ of $\Gm$ such that, for
any $f \in \mathcal{O}_M$, there does not exist any
solution $u(t)$ of \eqref{eq} such
that $u(t)$ converges, when $t$ goes to $\pm \infty$, to two
equilibrium points with the same Morse index.
In particular, homoclinic orbits are  precluded.
\end{theorem}

To prove Theorem \ref{th2}, we actually prove the genericity of the
transversality of the homoindexed orbits, which at once implies the
genericity of non-existence of such orbits.  In order to prove the
genericity of transversality and to get a meaningful result, we need
to perturb \eqref{eq} by arbitrary small perturbations, in such a way
that the perturbed semi-flow is still generated by a scalar parabolic
equation on $S^1$.  As already mentionned, perturbing the
non-linearity acts on the phase-space in a non-local way.
In the context of proving generic tranversality in the class of gradient
scalar-reaction diffusion equations, these problems were first
circumvented by Brunovsk\'y and Pol\'a\v{c}ik in \cite{BrunPola}.
They employed an equivalent formulation of transversality which
appeared earlier in \cite{Robbin}, \cite{HaLi86}, \cite{Sala90}, but
remained almost unnoticed for some time (such formulation has however
been used, already in the 1980's, by Chow, Hale and Mallet-Paret
\cite{CHM-P} in global bifurcation problems of heteroclinic and
homoclinic orbits).  This equivalent formulation says that $0$ is a
regular value of a certain mapping $\Phi$ (depending on the
perturbation parameter), defined on a space of functions of time with
values into the state space of the equation.  It is noteworthy that,
in this formalism, the elements, the image of which are $0$, are
precisely the trivial and non-trivial connecting orbits.  Using this
equivalent formulation of the transversality together with the
Sard-Smale theorem, given in Appendix A, Brunovsk\'y and Pol\'a\v{c}ik
achieved their proof of genericity of the transversality of stable
and unstable manifolds of equilibria with respect to the
non-linearity. \HB
In the proof of Theorem \ref{th2}, we follow the lines of the one of
\cite{BrunPola} (see also \cite{BrunRaug} and \cite{RJo}), by
introducing this equivalent formulation of transversality.  But, as in
\cite{BrunRaug}, the equivalent regular value formulation of
transversality takes place in a space of sequences (obtained by a time
discretization of the semi-flow associated with \eqref{eq}).  As
there, in the application of the Sard-Smale theorem, one returns to
the continuous time only in the last step, when one verifies the
non-degeneracy condition of the Sard-Smale, which gives rise to a
functional condition (see Theorem \ref{BDPhisurj}).
To find a perturbation satisfying this functional condition, we use 
as a central
argument the one-to-one property of homoindexed connecting orbits
stated in Proposition
\ref{injectif} (which is again a consequence of the Sturm property).
Notice that parabolic strong unique
continuation properties  lead to a weaker version of Proposition
\ref{injectif}, which is actually sufficient to show the genericity of
transversality and holds in any space dimension (see
\cite{BrunJolyRaug}).
\vskip 2mm
Using Theorems \ref{Poincare-B}, \ref{CzaRocha}, \ref{th-JR}, \ref{th1}
and \ref{th2}, we finally prove the following genericity of Morse-Smale systems
of type \eqref{eq}.

\begin{theorem}\label{corointro} {\bf Genericity of Morse-Smale property} \\
For any $f$ in the generic set $\mathcal{O}_M$,
$S_f(t)$ is a Morse-Smale dynamical system, that is,\HB
1) the non-wandering set consists only in a finite number of equilibria
and periodic orbits, which are all hyperbolic.
\HB
2) the unstable manifolds of all equilibria and periodic orbits
transversally intersect the local stable manifolds of all equilibria
and periodic orbits.
\end{theorem}

Since the second part of Theorem \ref{corointro} is a direct
consequence of Theorems \ref{CzaRocha}, \ref{th-JR}, \ref{th1}
and \ref{th2}, it only remains to show that the
non-wandering set is trivial. This is done by using the Poincar\'e-Bendixson
Property (Theorem \ref{Poincare-B}) together with arguments similar to the ones
used for vector fields in $\Rm^2$ 

\vskip 2mm

\begin{rem} All the above theorems are stated under the
``dissipative" condition \eqref{condf1} on the non-linearity $f$.
Actually, Theorems \ref{Poincare-B}, \ref{CzaRocha},
\ref{th-JR} and \ref{th1} still hold without assuming
\eqref{condf1} and the proofs are not more involved. Theorem \ref{th2} is also
still true, even if the proof
is a little more technical (but not more difficult). However, Theorem
\ref{corointro} is not
true in general without a dissipative condition on $f$ or more
generally without knowing that $S_f(t)$ admits a compact global
attractor. Indeed, if $S_f(t)$ has no compact global attractor,
already the number of equilibria (and periodic orbits) can be infinite.
For the reader convenience and for avoiding unnecessary technicalities, we
have chosen to impose the dissipative condition \eqref{condf1} in the
whole paper.
\end{rem}

\vskip 3mm

The paper is organized as follows. In Section 2, we recall the
fundamental properties of the zero number on $S^1$ as well as the
useful spectral properties of the linearized equations around
equilibrium points or periodic orbits. Section 3 is devoted to
relations between the lap number and the Morse indices and to the
one-to-one property of homoindexed orbits. In all these results, the
zero number plays an important role. While some of these results are
primordial in the proof of our main theorems,
others are stated for sake of completeness in the description of the
properties of Eq.  \eqref{eq} on the circle.
Section 4 contains the proof of the automatic transversality results stated in
Theorem \ref{th1}. In Section 5, we prove Theorem \ref{th2}, that is, 
the generic
non-existence of orbits connecting two hyperbolic equilibria
with same Morse index. Section 6 is focused on the study of the
non-wandering set and on the proof of Theorem \ref{corointro}.
The appendices
contain the necessary background for reading this paper.
In Appendix \ref{AppendixA}, we recall the Sard-Smale theorem in the form used
in Section 5. Appendix \ref{AppendixB} contains the basic definitions and
properties of exponential dichotomies and their applications to the functional
characterization of the transversality of the trajectories of the parabolic
equation. This Appendix plays an important role in the core of the 
proof of Theorem
\ref{th2}. Finally, in the Appendix \ref{AppendixD}, we describe the 
asymptotics
of the solutions of the linearized equations around connecting 
orbits, which are
one of the main ingredients of the proof of Theorem \ref{th1}.

\vskip 3mm

{\noindent \bf Acknowledgements:} The authors wish to thank P.
Brunovsk\'y, R. Czaja and C. Rocha for fruitful discussions.


\section{Preliminaries and auxiliary results}

In the introduction, we have already seen that (even without Assumption
\eqref{condf1}) for any $u_0 \in H^s(S^1)$, $s \in (3/2,2)$,
Equation  \eqref{eq} admits a local mild solution $u(t) \in C^0 ([0,
\tau_{u_0}), H^s(S^1))$ (see \cite{He81}). Moreover, this solution
$u(t)$ is classical and belongs to $C^0((0, \tau_{u_0}), H^2(S^1)) \cap
C^1((0, \tau_{u_0}), L^2(S^1)) \cap C^{\theta} ((0,\tau_{u_0}),
H^s(S^1))$, where
$\theta =1 -s/2$. In addition, the function $u_t(t) : t \in (0,
\tau_{u_0}) \mapsto u_t(t) \in H^{\ell}(S^1)$, $0 \leq \ell <2$, is locally
H\"{o}lder-continuous (see \cite[Theorem 3.5.2]{He81}). Since
$u(t)$ is in $C^0 ((0, \tau_{u_0}), H^2(S^1))$, the term $f(x,u,u_x)$
belongs to $C^0((0, \tau_{u_0}), H^1(S^1))$ and thus $u_{xx}=u_t
-f(x,u,u_x)$ is in $C^0((0, \tau_{u_0}), H^1(S^1))$.
    In particular, $u(t)$ belongs to $C^0((0,\tau_{u_0}), H^3(S^1))$,
    which is continuously embedded into
$C^0((0,\tau_{u_0}), C^2(S^1))$. If moreover the condition
\eqref{condf1} holds, then $\tau_{u_0}= \infty$ for
every $u_0$.

In the course of this paper, we often need to consider the linearized
equation along a bounded trajectory $u(t) \equiv S_f(t)u_0$, $t \in
\Rm$, of Eq. \eqref{eq}, that is, the equation
    \begin{equation}
\label{2linu}
v_t=v_{xx}+ D_u f(x,u,u_x)v + D_{u_x} f(x,u,u_x)v_x~,~
t\geq \sigma ~,
\quad v(\sigma,x)= v_0 ~,
\end{equation}
where $v_0$ belongs to $L^2(S^1)$. Since $u(t)$ is in
$C^{\theta}(\Rm,H^s(S^1))$, the coefficients $D_uf(x,u(t),u_x(t))$
and $D_{u_x}f(x,u(t),u_x(t))$ are locally H\"{o}lder-continuous from
$\Rm$ into $C^0(S^1)$. Thus,
we deduce from \cite[Theorem 7.1.3]{He81}
that, for any $v_0 \in L^2(S^1)$, for any $\sigma \in \Rm$, there
exists a unique (classical)
solution $v(t) \in C^0 ([\sigma,+ \infty),L^2(S^1)) \cap
C^0 ((\sigma,+ \infty), H^{\ell}(S^1))$ of
\eqref{2linu} such that $v(\sigma)=v_0$, where $\ell$ is any real
number with $0  \leq \ell <2$.
Setting $T_{u}(t,\sigma)v_0=v(t)$, we
define a family of continuous linear evolution operators on
$L^2(S^1)$. We remark that the evolution operator
$T_{u}(t,\sigma)$, $t \geq \sigma$, associated with the trajectory
$u$ is injective and that its range is dense in $H^s(S^1)$.

\vskip 1mm

We complete these generalities by remarking that,
in what follows, we sometimes consider the difference $w=
u_1-u_2$ between two solutions $u_1\in C^0([0,+\infty), H^s(S^1))$ and
$u_2 \in C^0([0,+\infty), H^s(S^1))$ of Eq.
\eqref{eq}. The difference $w$ is a solution of the following
linear equation,
    \begin{equation}
\label{2u12}
w_t= w_{xx}+a(x,t)w_x + b(x,t)w~,~ t\geq \sigma ~,
\quad w(\sigma,x)= w_0 ~,
\end{equation}
where
\begin{equation}
\label{2abu12}
\left\{\begin{array}{l}
a(x,t)=\int_0^1 f'_{\partial_x u}(x, \theta u_2
+(1-\theta) u_1,\partial_x (\theta u_2
+(1-\theta) u_1)) d\theta   \\
b(x,t)=\int_0^1 f'_u(x,\theta u_2
+(1-\theta) u_1,\partial_x (\theta u_2
+(1-\theta) u_1))d\theta \end{array}\right.
\end{equation}
We emphasize that, since $u_t$ is locally H\"{o}lder-continuous from
$(0,+\infty)$ into $H^{\ell}(S^1)$, $3/2 < \ell <2$ and $u_x$ is
continous from $(0,+\infty)$ into $H^2(S^1)$, the coefficients
$a(x,t)$ and $b(x,t)$ belong to $C^1(S^1 \times (0, +\infty), \Rm)$.

\subsection{The lap number property}

The lap number property, or zero number property, is the fundamental
property of one-dimensional scalar parabolic equations. We recall that, for
any $\varphi \in C^1(S^1)$,
the zero number $z(\varphi)$ is defined as the (even) number of strict
sign changes of $\varphi$.

\begin{theorem}\label{th-zeros}
1) Let $T>0$, $a\in\Cm^1(S^1\times [0,T] ,\Rm)$ and $b\in\Cm^0(S^1\times
[0,T] ,\Rm)$. Let $v:S^1\times (0,T)\rightarrow \Rm$ be a classical bounded
non-trivial solution of
$$
  v_t= v_{xx} +a(x,t) v_x +b(x,t)v~,
$$
Then, the number $z(v(t))$ of zeros of $v(t)$ is finite and
non-increasing in time $t \in [0,T]$ and strictly decreases at $t=t_0$ if and
only if $x\mapsto v(x,t_0)$ has a multiple zero.

2) If $u$ and $v$ are two solutions in $C^0([0,T],H^s(S^1))$
of \eqref{eq}, then $u_t$ and $u-v$ satisfy the lap number property
stated in Statement 1) on the time interval $(0,T]$.
\end{theorem}

Such kind of results goes back to Sturm \cite{Sturm} in the case where
$a$ and $b$
are time-independent. The non-increase of the number of zeros in the
time-dependent problems has been obtained in \cite{Nickel} and
\cite{Matano}. The property of strict decay first appeared in
\cite{AngenFied} in the case of analytic coefficients. It has
been generalized in \cite{Matano2} and \cite{Angenent}. \HB
We notice that the statement 2) is a direct consequence of the first
statement and of the remarks made at the beginning of this
section.


\subsection{The spectrum of the linearized operators}\label{subsecLin}

We recall here the Sturm-Liouville properties of the linearized
operators associated to
Eq. \eqref{eq}. These results mainly come from \cite{Angenent}
and \cite{AngenFied},
together with a property obtained in \cite{JR}.

Let $e \in H^s(S^1)$ be an equilibrium point of \eqref{eq}.
We introduce the linearized operator$L_e$ on $L^2(S^1)$, with domain
$H^2(S^1)$, defined by
\begin{equation}
\label{MapLe}
L_ev= v_{xx}+ D_uf(x,e,e_x)v +
D_{u_x}f(x,e,e_x)v_x~,
\end{equation}
and consider the linearized equation around $e$, given by
\begin{equation}
\label{eqeline}
\begin{split}
v_t(x,t) &=L_e v(x,t)
~, \quad (x,t)\in S^1\times (0,
+\infty)~,
\cr
v(x,0) & = v_0(x)~, \quad x\in S^1~.
\end{split}
\end{equation}
The operator $L_e$ is a sectorial operator and a Fredholm operator with compact
resolvent.
Therefore, its spectrum consists of a sequence of isolated
eigenvalues of finite multiplicity.
Let $(\lambda_i)_{i\in\Nm}$ be the spectrum of $L_e$,  the
eigenvalues being repeated according to
their multiplicity and being ordered so that $\Re (\lambda_{i+1})\leq
\Re (\lambda_i)$.

\begin{prop}\label{spectreLe}
The first eigenvalue $\lambda_0$ is real and simple and the
corresponding eigenfunction
$\varphi_0\in H^2(S^1)$ does not vanish on $S^1$. The other
eigenvalues go by pairs
$(\lambda_{2j-1},\lambda_{2j})$ and $\Re (\lambda_{2j+1})
< \Re (\lambda_{2j})$ for
all $j\geq 0$. The pair $(\lambda_{2j-1},\lambda_{2j})$ consists of
either two simple complex conjugated
eigenvalues, or two simple real eigenvalues with
$\lambda_{2j}<\lambda_{2j-1}$, or a real eigenvalue with
multiplicity equal to two. Finally, if $\varphi$ is a real function
belonging to the two-dimensional
generalized eigenspace corresponding to
$(\lambda_{2j-1},\lambda_{2j})$, then $\varphi$ has
exactly $2j$ zeros which are all simple.
\end{prop}

Let $\Gamma = \{ \gamma(x,t) \, | \, t \in [0,p] \}$ be a periodic orbit
of \eqref{eq} of minimal period $p$. We
consider the linearized equation
\begin{equation}
\label{linearized}
\varphi_t =\varphi_{xx}+ D_{u}f(x,\gamma,\gamma_{x})\varphi
+ D_{u_x}f(x,\gamma,\gamma_{x})\varphi_x~, \quad t\geq \sigma~, \quad
\varphi (x,\sigma) = \varphi_0(x)~.
\end{equation}
Let $s\in (3/2,2)$, we introduce the operator $\Pi(t,\sigma): H^s(S^1)
\longrightarrow H^s(S^1)$,
       defined by $\Pi(t,\sigma)\varphi_0= \varphi(t)$ where
$\varphi(t)$ is the solution of the linearized equation
\eqref{linearized}. The operator $\Pi(p,0)$ is called the {\sl period
map}. Due to the regularization properties of the parabolic equation,
$\Pi(p,0)$ is compact. Its spectrum consists of zero and a sequence
of eigenvalues $(\mu_i)_{i\in\Nm}$ converging to zero, where we
repeat the eigenvalues according to their multiplicity and order them
such that $|\mu_{i+1}\leq |\mu_i|$. Notice that $0$ is not
an eigenvalue of $\Pi(p,0)$ due to the backward uniqueness property
of the parabolic equation. Moreover, we also remark that, since
\eqref{eq} is an autonomous equation, $1$ is always
an eigenvalue of $\Pi(p,0)$ with eigenfunction $\gamma_t(0)$.

\begin{prop}\label{spectre-Pi}
The first eigenvalue $\mu_0$ is real and simple and the corresponding
eigenfunction $\varphi_0\in H^s(S^1)$ does not vanish on $S^1$. The
other eigenvalues go by pairs $(\mu_{2j-1},\mu_{2j})$ and
$|\mu_{2j+1}|< |\mu_{2j}|$ for all $j\geq 0$. The pair
$(\mu_{2j-1},\mu_{2j})$ consists of either two simple complex conjugated
eigenvalues, or two simple real eigenvalues of the same sign, or a real
eigenvalue with multiplicity equal to two. In particular, $-1$ is
never an eigenvalue. Finally, if
$\varphi$ is a real function belonging to the two-dimensional
generalized eigenspace corresponding to $(\mu_{2j-1},\mu_{2j})$, then
$\varphi$ has exactly $2j$ zeros which are all simple.
\end{prop}


\section{Main consequences of the lap number properties on connecting 
orbits}\label{LapMorseInj}

In this section, we describe several important consequences of the
lap number theorem on the properties of orbits connecting hyperbolic
critical elements. Most of them will be used in the core of the
proofs of Theorems \ref{th1} and \ref{th2}.

\subsection{Relations between Morse indices and lap numbers}\label{Sequper1}

In this paragraph, we consider orbits connecting hyperbolic equilibrium
points and hyperbolic periodic orbits and give several inequalities 
involving Morse
indices and numbers of zeros. The ideas of these properties were already
contained in \cite{CR}, where they have been proved in the case of connections
between two periodic orbits. We complete the results of \cite{CR} by
considering the cases where equilibrium points are also involved.
Actually, some of the proofs are slightly simpler in these cases since, if
$e$ is an equilibrium point and $\gamma(t)$ a periodic orbit, $e-\gamma(t)$ is
periodic, whereas the difference between two periodic orbits may be only
quasiperiodic.

The properties stated in this section have their own interest. In the other
parts of this paper, we will use them in the case of a connection between an
equilibrium point and a periodic orbit. Therefore, we mainly restrict the
proofs to this case. The omitted proofs are similar.

\vskip 2mm

The following theorem corresponds to \cite[Theorems 5.2 and
6.2]{CR}. We recall that $z(v)$ denotes the strict number of sign
changes of the function $x\mapsto v(x)$.

\begin{theorem}\label{4CRsu}
Let $\gamma(x,t)$ be a hyperbolic periodic orbit of \eqref{eq} of
minimal period $p$ and let $\Gamma=\{\gamma(t),~t\in [0,p)\}$. \HB
1) For any $u_0 \in W^s_{loc}(\Gamma) \setminus \Gamma$, there exist
$a \in \Gamma$ and $\kappa >0$ such that $\lim_{t \to \infty}
e^{\kappa t} \| S_f(t) u_0 -S_f(t) \gamma (a)\|_{H^s}=0$.
Moreover,
    \begin{equation}
\label{4zstable}
z(u_0 -\gamma (a)) \geq \begin{cases} i(\Gamma) +1 =2q  & \hbox{ if }
i(\Gamma) = 2q -1 , \cr
i(\Gamma) + 2 = 2q +2 & \hbox{ if } i(\Gamma) =2q.
\end{cases}
\end{equation}
2) For any $u_0 \in W^u(\Gamma) \setminus \Gamma$, there exist
$\tilde{a} \in \Gamma$ and $\tilde{\kappa} >0$ such that $\lim_{t \to \infty}
e^{ \tilde{\kappa} t} \| S_f(-t) u_0 -S_f(-t) \gamma
(\tilde{a})\|_{H^s}=0$. Moreover,
    \begin{equation}
\label{4zunstab}
z(u_0 -\gamma (a)) \leq \begin{cases} i(\Gamma) -1 =2q - 2 & \hbox{ if }
i(\Gamma) = 2q -1 , \cr
i(\Gamma) = 2q  & \hbox{ if } i(\Gamma) =2q.
\end{cases}
\end{equation}
\end{theorem}

\begin{demo} For sake of completeness, we give
a short proof of statement 1) which is
slightly different from the one of \cite{CR}. The proof of statement 2) is
similar.

The first part of the statement 1) is just a
reminder of the fact that the local stable manifold of $\Gamma$ is
the union of the local strongly stable manifolds of all the points
$\gamma (b)$, $b \in [0, p]$ as explained in Appendix \ref{App-per}. The
second part of 1) directly follows from Corollary \ref{CappD02}.
Indeed, let $v(t)= S_f(t)u_0 -\gamma(t+a) \equiv u(t) - \gamma(t+a)$.
Then,
there exists a complex eigenvalue $\mu_i$ of the period map $\Pi(p+a ,a)$ with
$| \mu_i| <1$ such that $v(np)$ satisfies one of the asymptotic behaviors
(i)-(iv) described in Corollary \ref{CappD02}.
If the index $i(\Gamma)$ is equal to $2q-1$, then $\mu_{2q-1}= 1$
and thus $i \geq
2q$, which implies that the number of zeros $z(v(np))$ is at least
equal to $2q$ for
  $n$ large enough. If the index $i(\Gamma)$ is equal to $2q$, then
$\mu_{2q}=1$ and thus
$i\geq 2q+1$, which implies that the number of zeros $z(v(n))$ is at
least equal to
$2q+2$ for $n$ large enough. Since $v(t)$ is the
difference of two solutions of \eqref{eq}, Theorem \ref{th-zeros} 
shows that these lower
bounds on $z(v(np))$ for large $n\in\Nm$ hold in fact for all $t\in\Rm$.
\end{demo}

Of course, the corresponding properties are true for hyperbolic 
equilibrium points.
Since their proof is similar to the one of Theorem
\ref{4CRsu} (and even simpler), it is omitted.

\begin{theorem}\label{4CRsu2}
Let $e(x)$ be a hyperbolic equilibrium point of \eqref{eq}.\\
1) For any $u_0 \in W^s_{loc}(e) \setminus \{e\}$,
$$
z(u_0 -e) \geq \begin{cases} i(e)+1=2q  & \hbox{ if }
i(e) = 2q -1 , \cr
i(e) = 2q & \hbox{ if } i(e) =2q.
\end{cases}
$$
2) For any $u_0 \in W^u(e) \setminus \{e\}$,
$$
z(u_0 -e) \leq \begin{cases} i(e)-1 =2q - 2 & \hbox{ if }
i(e) = 2q -1 , \cr
i(e) = 2q  & \hbox{ if } i(e) =2q.
\end{cases}
$$
\end{theorem}

The following two lemmas, which are rather simple, are useful in the following
sections.
\begin{lemma} \label{zconst}
If $e$ is an equilibrium point of \eqref{eq} and $\gamma (t)$ is a periodic
solution of \eqref{eq} of minimal period $p > 0$, then the zero
number $z(e- \gamma (t))$ is constant and thus, for any time $t$, the function
$x\mapsto e(x)- \gamma (x,t)$ has no multiple zero. The same properties
hold if one considers the difference between two distinct equilibrium points or
two distinct periodic solutions.
\end{lemma}
\begin{demo}
By Theorem \ref{th-zeros}, the number of zeros of $v(t)= e- \gamma
(t)$ is non-increasing and  strictly decreases only at the times $t$
where $v(t)$ has a multiple zero. If $v(t)$ has a multiple zero at $t=t_0$,
then Theorem \ref{th-zeros} and the periodicity of $v$ imply that, for any
$\varepsilon\in(0,p)$,
$$
z(v(t_0-\varepsilon)) > z(v(t_0 +\varepsilon))\geq
z(v(t_0+p-\varepsilon))=z(v(t_0-\varepsilon))~,
$$
which leads to a contradiction. Thus $v(t)$ has no multiple zero.
\end{demo}
\vskip 1mm

\begin{rem}
Since, for each $t$, $v(t)$ has no multiple zero and
$H^s(S^1)$ is embedded in $\Cm^1(S^1)$, there exists a neighbourhood
$B_{H^s}(v(t), 2
\varepsilon_{v(t)})$ on which the zero number is constant.  Since the curve
$\{e- \gamma (t) \, | \, t \in [0, p]\}$ is compact,
there exists a finite covering $\cup_{i=1}^{n}
B_{H^s(S^1)}(e- \gamma (t_i),\varepsilon_i)$ of the set $\{e- \gamma
(t) \, | \, t \in \Rm\}$ in $H^s(S^1)$, on which the zero number in
constant.
\end{rem}

\begin{lemma}\label{zconst2}
Let $e^{-}$ be a hyperbolic equilibrium of \eqref{eq} and let
$\Gamma^{+}$ be a hyperbolic
  periodic orbit of minimal period $p>0$. Let $u_0\in W^u(e^-) \cap
W^s_{loc}(\Gamma^{+})$ and let $u(t)$, $t\in\Rm$, be the solution of \eqref{eq}
with $u(0)=u_0$. Let $a \in [0,p)$ be such that
$\lim_{t\rightarrow +\infty}\| u(t)-\gamma^+(a+t)\|_{H^s}=0$. Then, for any
time $t\in\Rm$,
$$
z(u(t)-\gamma^+(a+t))\leq z(\gamma^+(a)-e^-)\leq z(u(t)-e^-)~.
$$
The same property holds for any orbit connecting hyperbolic equilibrium
points or hyperbolic periodic orbits.
\end{lemma}

\begin{demo}
We set $v_+(t)=u(t)-\gamma^+(a+t)$ and $v_-(t)=u(t)-e^-$. We notice that the
lap number property stated in Theorem \ref{th-zeros} holds for $v_\pm$.
Lemma \ref{zconst} shows that
$z(\gamma^+(a)-e^-)=z(\gamma^+(a+t)-e^-)$. Moreover, $\gamma^+(a)-e^-$ has no
multiple zeros due to Lemma \ref{zconst} and thus its number of zeros is
stable with respect to small enough perturbations in  $H^s(S^1)$. For
any large $t_0$, $v_+(t_0)$ is small enough so that
$$z(\gamma^+(a)-e^-)=z(\gamma^+(a+t_0)-e^-)=z(v_+(t_0)+\gamma^+(a+t_0)-e^-)=z(v_
-(t_0))~.$$
Applying Theorem \ref{th-zeros}, we get that for all $t\leq t_0$,
$z(\gamma^+(a)-e^-)\leq z(v_-(t))$.\\
The inequality $z(v_{+}(t)) \leq z(\gamma^+(a)- e^{^-})$, for all $t
\in \Rm$, is proved in a similar way. \HB
The proof is the same in the case of orbits connecting hyperbolic
equilibria or hyperbolic periodic orbits. These inequalities have
been previously proved in the case of orbits connecting two
hyperbolic periodic orbits in \cite[Theorem 7.3]{CR}
\end{demo}

As a direct consequence of Lemma \ref{zconst2}, Theorem \ref{4CRsu} and
Proposition \ref{4CRsu2}, we obtain the following result.

\begin{coro} \label{i+i-}  Let $e^{\pm}$ and $\Gamma^{\pm}$ be
hyperbolic equilibria and periodic orbits of Eq. \eqref{eq}. \HB
1) If $W^u(e^-) \cap W^s_{loc}(\Gamma^{+})\neq \emptyset$, then
$$
i(\Gamma^{+}) +1\leq i(e^-)~.
$$
Moreover, if $i(\Gamma^{+}) =2q^{+}$, $q^{+}
>0$, then $i(\Gamma^{+})+2\leq i(e^-) $. \HB
2) If $W^u(\Gamma^-) \cap W^s_{loc}(e^{+})\neq \emptyset$, then
$$
i(e^{+}) \leq i(\Gamma^{-})~.
$$
Moreover, if $i(\Gamma^{-}) =2q -1$, $q\geq 1$,
then $i(e^{+}) +1 \leq  i(\Gamma^{-})$.\HB
3) If $W^u(e^-) \cap W^s_{loc}(e^{+})\neq \emptyset$, then
$$
i(e^+)\leq i(e^-)
$$
and the equality is possible if and only if $i(e^+)=i(e^-)$ is even.
\end{coro}

\subsection{One-to-one property of homoindexed orbits}\label{Sinjectif}

The proposition proved in this section plays a central role in the
construction of a suitable perturbation to break homoindexed orbits
(see Section 5).  A similar one-to-one property has been previously used
in \cite{JR} to make the periodic orbits hyperbolic.  As already
indicated, the proof relies on the decay of the lap-number stated in
Theorem \ref{th-zeros}.  In space dimension higher than one, there is
no equivalent of the lap-number property.  However, one
can use unique continuation properties to obtain a weaker but
useful equivalent of Proposition \ref{injectif}, see
\cite{BrunJolyRaug}.

\begin{prop} \label{injectif}
Let $e_{\pm}$ be two hyperbolic equilibria such that $i(e_{-})= 
i(e_{+})=m=2m'$ is even,
and let $u(t)$ be a connecting orbit between $e_{-}$ and $e_{+}$.
In the case where $e_{-}= e_{+}=e$, we assume that $u(t) \ne e$.
Then the following properties hold: \HB
1) For any $t \in \Rm$, for any $x \in S^1$, $(u(x,t), \partial_x u(x,t)) \ne
(e_{\pm}(x), \partial_x e_{\pm}(x))$.\HB
2)  The map $(x,t) \in S^1 \times \Rm \mapsto (x, u(x,t), \partial_x
u(x,t))$ is one to one.
\end{prop}

\begin{demo} As in Section 2.2, we denote by
$L_{e_{\pm}}$ the linearized operator around the equilibrium
$e_{\pm}$ and by $\lambda_i^{\pm}$, $i \geq 0$, its eigenvalues,
counted with their multiplicities.

Since $u(t)$ converges to $e_\pm$ when $t$ goes to $\pm\infty$, according to
Corollary \ref{CappD0}, there exists an eigenvalue $\lambda^\pm_{i^\pm}$ of
$L_{e_\pm}$ such that the asymptotic behaviors of $u$ are given by
\begin{equation}
\label{5u-}
u(t) = e_\pm + e^{\Re (\lambda^\pm_{i^\pm}) t} \psi^\pm_{i^\pm}(t)
+ o(e^{\Re(\lambda^\pm_{i^\pm}) t})~~~\text{ when }~t\rightarrow \pm\infty
\end{equation}
where $\psi^\pm_{i^\pm}$ corresponds to one of the possible asymptotic
behaviors (i)-(iv) described at the beginning of Section
\ref{Sheteind}. We emphasize
that the term $o(e^{\Re(\lambda^\pm_{i^\pm}) t})$ has to be understood in the
sense of the $H^s(S^1)$ topology (and thus this term is also small in
the $C^1$-sense). We
  recall that if $i^\pm=2q-1$ or $2q$ then $\psi^\pm_{i^\pm}$ has exactly
$2q$ zeros which are simple. Finally, notice that since $i(e_{-})= i(e_{+})=m$
is even, $\lambda^\pm_{m}<0<\lambda^\pm_{m-1}$ are simple real
eigenvalues and that $i^-\leq m-1$ and $i^+\geq m$.

\vskip 2mm

We are now ready to prove the first assertion.
We introduce the functions $v_{\pm}(t) = u(t) - e_{\pm}$.
By the second part of Theorem \ref{th-zeros}, the number of zeros $z(v_{\pm})$
is finite and non-increasing and it strictly decreases at some time
$\tau_{\pm}$ if and only if the function $x \mapsto
v_{\pm}(x,\tau_{\pm})$ has a multiple zero.
But Proposition \ref{4CRsu2} and Lemma \ref{zconst2} at once
imply that $z(e_{+}-e_{-}) = m = z(v_{+}(t)) = z(v_{-}(t))$ for any 
$t$ and thus
that the map $x \mapsto u(x,t)-e_{\pm}(x)$ has no multiple
zero. We deduce that
\begin{equation}
   \label{5vectpro}
\lambda^-_{i^-}= \lambda_{m-1}^{-}~,~~\psi^-_{i^-}=
\varphi_{m-1}^{-}~,~~\lambda^+_{i^+}= \lambda_{m}^{+}~\text{ and
}~\psi^+_{i^+}=\varphi_{m}^{+}~,
\end{equation}
where $\varphi_{m-1}^{-}$ and
$\varphi_{m}^{+}$ are eigenfunctions corresponding to the simple real
eigenvalues $\lambda_{m-1}^{-}$ and $\lambda_{m}^{+}$ respectively.

\vskip 2mm

We next prove by contradiction that the second statement
holds. Assume that the map $(x,t) \in S^1 \times \Rm \mapsto (x,
u(x,t), \partial_x
u(x,t))$ is not injective. Then there exist
$x_0$, $t_0 \in \Rm$ and $\tau_0 \in \Rm$, such
that
$$
u(x_0, t_0)=u(x_0, t_0 + \tau_0)~, \quad \partial_x u(x_0, t_0)=
\partial_x u(x_0, t_0 +\tau_0)~.
$$
The function $v(x,t)= u(x, t +  \tau_0)- u(x,t)$  satisfies
$v(x_0,t_0)=0$ and $\partial_x v(x_0, t_0)=0$. It is not identically
zero since it is a non-trivial connecting orbit. Thus, due to the second part
of Theorem \ref{th-zeros}, the zero number $z(v(t))$ is
non-increasing and drops
strictly at $t=\tau_0$. The properties \eqref{5u-} and
\eqref{5vectpro} imply that
    \begin{equation*}
\begin{split}
v(t)=& ~e^{\lambda_{m-1}^{-}t}(e^{\lambda_{m-1}^{-} \tau_0} -1)
\varphi_{m-1}^{-} +o(e^{\lambda_{m-1}^{-}t})~, \quad \hbox{ as } t
\to -\infty \cr
v(t)=& ~e^{\lambda_{m}^{+}t}(e^{\lambda_{m}^{+}\tau_0} -1)
\varphi_{m}^{+} +o(e^{\lambda_{m}^{+}t})~, \quad \hbox{ as } t
\to \infty ~.
\end{split}
\end{equation*}
Thus, $z(v(t)) = m$, for any $t \in \Rm$, which contradicts the fact
that $z(v(t))$ drops at $t=\tau_0$ and the second statement holds.
\end{demo}


\section{Automatic transversality results}

This section is devoted to the proof of Theorem \ref{th1}. In Section
4.1,  we prove the first statement of Theorem \ref{th1}, whereas in
Section 4.2 , we prove the second statement.

\subsection{Automatic transversality of heteroindexed orbits connecting
two equilibria}\label{Sheteind}

Let $e_-$ and $e_+$ be two hyperbolic equilibrium points of
\eqref{eq} with different Morse indices $i(e_-)$ and $i(e_{+})$.
Following the notations of Section \ref{subsecLin}, we denote
$L_{e_\pm}$ the corresponding linearized operators and by
$(\lambda_i^\pm,\varphi_i^\pm)$ their set of eigenvalues and
generalized eigenfunctions.

Assume that $W^u(e_-)\cap W^s_{loc}(e_+)\neq\emptyset$ (otherwise the
intersection is transversal by definition). Notice that Corollary \ref{i+i-}
and the fact that $i(e_-) \neq i(e_+)$ imply that $i(e^+)<i(e^-)$. Let
$u(t)$ be a global
solution of \eqref{eq} with $u(0)\in W^u(e_-)\cap W^s_{loc}(e_+)$. Since
$W^u(e_-)$ is a finite-dimensional manifold, there exists a finite
basis $(v_1^0,\ldots,v_p^0)$ of $T_{u(0)}W^u(e_-)\cap
T_{u(0)}W^s_{loc}(e_+)$. Let $v_k(t)$ be the (global) solutions of
\begin{equation}
\label{3eqvk}
\partial_t v_k=\partial^2_{xx}v_k+ D_u f(x,u,u_x)v_k+ D_{u_x}f(x,u,u_x)
\partial_x v_k~, \quad v_k(0)=v_k^0~.
\end{equation}
Notice that, since $v_k^0$ belongs to the finite-dimensional space
$T_{u(0)}W^u(e_-)$, the solution $v_k(t)$ exists for any $t \in \Rm$.
Corollary \ref{CappD1} gives all the possible precise asymptotic
behaviors of $v_k(t)$ when $t$ goes to $\pm\infty$. For each $k$,
there exist an eigenvalue $\lambda_{i_k^-}^-$ of $L_{e_-}$ with
positive real part such
that, when $t\rightarrow -\infty$, the asymptotic behavior of $v_k(t)$
in $H^s(S^1)$ is described by the following possibilities:
\begin{description}
\item[(i)]   if $\lambda_{i_k^-}^-$ is a simple real eigenvalue with
eigenfunction $\varphi_{i_k^-}^-$,
then there exists $a_k\in\Rm-\{0\}$ such that $v_k(t)=
a_k e^{\lambda_{i_k^-}^- t}\varphi_{i_k^-}^-
+ o(e^{\lambda_{i_k^-}^- t}) \equiv \psi_{i_k^{-}}^{-}(t) +
o(e^{\lambda_{i_k^-}^- t})$.

\item[(ii)]
    If $\lambda_{i_k^-}^-=\lambda_{i_k^-+1}^-$ is a double real
eigenvalue with two
independent eigenfunctions $\varphi_{i_k^-}^-$ and
$\varphi_{i_k^-+1}^-$, then there exist
$(a_k,b_k)\in\Rm^2-\{(0,0)\}$ such that $v_k(t)=a_k e^{\lambda_{i_k^-}^-
t}\varphi_{i_k^-}^- + b_k e^{\lambda_{i_k^-}^- t}\varphi_{i_k^-+1}^- +
o(e^{\lambda_{i_k^-}^- t})\equiv \psi_{i_k^{-}}^{-}(t) +
o(e^{\lambda_{i_k^-}^- t})$.

\item[(iii)]
    If $\lambda_{i_k^-}^-=\lambda_{i_k^-+1}^-$ is an algebraically
double real eigenvalue
with eigenfunction $\varphi_{i_k^-}^-$ and with generalized
eigenfunction $\varphi_{i_k^-+1}^-$, then there exist
$(a_k,b_k)\in\Rm^2-\{(0,0)\}$ such that $v_k(t)=(a_k+b_k t)e^{\lambda_{i_k^-}^-
t}\varphi_{i_k^-}^- + b_k e^{\lambda_{i_k^-}^- t}\varphi_{i_k^-+1}^- +
o(e^{\lambda_{i_k^-}^- t}) \equiv \psi_{i_k^{-}}^{-}(t) +
o(e^{\lambda_{i_k^-}^- t})$.

\item[(iv)]
If $\lambda_{i_k^-}^-=\overline{\lambda_{i_k^-+1}^-}$ is a
(simple) nonreal eigenvalue
with eigenfunction
$\varphi_{i_k^-}^-=\overline{\varphi_{i_k^-+1}^-}$, then there exist
$(a_k,b_k)\in\Rm^2-\{(0,0)\}$ such
that
    \begin{equation*}
    \begin{split}
v_k(t)=e^{\Re(\lambda_{i_k^-}^-)t} &\Big(
(a_k \cos(\Im(\lambda_{i_k^-}^-)t) - b_k \sin( \Im(\lambda_{i_k^-}^-)
t)) \Re(\varphi_{i_k^-}^-) \cr
& - (a_k \sin( \Im(\lambda_{i_k^-}^-) t) +
b_k\cos( \Im(\lambda_{i_k^-}^-) t) ) \Im( \varphi_{i_k^-}^- )\Big) +
o(e^{\Re(\lambda_{i_k^-}^-) t}) \cr
\equiv \psi_{i_k^{-}}^{-}(t) + o(&e^{\Re(\lambda_{i_k^-}^-) t}) ~.
\end{split}
\end{equation*}
\end{description}
The vectors $v_k(t)$ have the same type of behaviors when
$t\rightarrow +\infty$, provided
$\lambda_{i_k^-}^-$ is replaced by an eigenvalue $\lambda_{i_k^+}^+$
of $L_{e_+}$ with
negative real part.

\begin{lemma}\label{3lemme1}
Without loss of generality, we may assume that the behaviors of the
functions $v_k(t)$ when $t$ goes to $-\infty$ are different in the sense that
the corresponding family $(\psi_{i_k^{-}}^{-})$, $1 \leq k \leq p$,
  is free and hence generates a finite-dimensional vector space of
dimension $p$.
\end{lemma}
\begin{demo}
This lemma is a simple consequence of a Gram-Schmidt process.  Without loss of
generality, we can assume that $\lambda_{i_1^-}^-$ has the smallest
real part among the family
$(\lambda_{i_1^-}^-,\ldots,\lambda_{i_p^-}^-)$.  If there exists $k>1$
such that $\Re(\lambda_{i_k^-}^-)=\Re(\lambda_{i_1^-}^-)$ with
asymptotic behavior of type (i) or with asymptotic behavior of type
(ii)-(iv) such that the pairs $(a,b)$ for $v_1(t)$ and $v_k(t)$ are
linearly dependent, then we can replace $v_k$ by $v_k-\alpha v_1$ such
that the real part of $\lambda_{i_k^-}^-$ increases.  Notice that
$(v_1^0,\ldots,v_p^0)$ is still a basis of $T_{u(0)}W^u(e_-)\cap
T_{u(0)}W^s(e_+)$.  Assume now that $\lambda_{i_2^-}^-$ has the
smallest real part among the real parts of the family
$(\lambda_{i_2^-}^-,\ldots,\lambda_{i_p^-}^-)$ (which can be the
same as the real part of $\lambda_{i_1^-}^-$, but not smaller).  If there
exists $k>2$ such that $\Re(\lambda_{i_k^-}^-)=\Re(\lambda_{i_2^-}^-)$
and $v_k(t)$ has an asymptotic behavior linearly dependent of
the one of $v_2(t)$, then we can replace $v_k$ by $v_k-\alpha v_2$
so that the real part of $\lambda_{i_k^-}^-$ increases.  We pursue
the process until the end to obtain Lemma \ref{3lemme1}.
\end{demo}

\noindent{{\textbf{Proof of the first statement of Theorem
\ref{th1}:}}

In what follows, we assume that the basis $(v_1^0,\ldots,v_p^0)$ of
$T_{u(0)}W^u(e_-)\cap T_{u(0)}W^s(e_+)$ has been chosen as in Lemma
\ref{3lemme1}. By definition $W^u(e_-)$ intersects
$W^s(e_+)$ transversally if and only if $T_{u(0)}W^u(e_-) +
T_{u(0)}W^s(e_+)=H^s(S^1)$. We know that
$\dim(T_{u(0)}W^u(e_-))=i(e_-)$ and
$\codim(T_{u(0)}W^s(e_+))=i(e_+)$. Thus, $T_{u(0)}W^u(e_-) +
T_{u(0)}W^s(e_+)=H^s(S^1)$ if and only if
\begin{equation}
\label{dimInter}
\dim\left(T_{u(0)}W^u(e_-)\cap T_{u(0)}W^s(e_+)\right)\leq i(e_-) -
i(e_+) ~.
\end{equation}

The proof of the first statement of Theorem \ref{th1} consists in the
careful study of three different cases.

If $i(e_+)=0$, then $T_{u(0)}W^s(e_+)=H^s(S^1)$ and the
transversality trivially holds.

Assume that $i(e_+)=2q-1$ is odd, which implies that the eigenvalues of
$L_{e_+}$ with negative real part are $(\lambda^+_i)_{i\geq 2q-1}$.
By Proposition \ref{spectreLe}, the corresponding generalized eigenfunctions
$(\varphi^+_i)_{i\geq 2q-1}$ have all at least $2q$ zeros.  Due to
Corollary \ref{CappD1} in Appendix C, each function $v_k(t)$ has at
least $2q$ zeros for large $t$.  Since by Theorem \ref{th-zeros} the
number of zeros of $v_k(t)$ is non-increasing,
$v_k(t)$ has at least $2q$ zeros for every time $t\in\Rm$.  Applying
    Corollary \ref{CappD1} again, we obtain that necessarily
$i^-_k \geq 2q -1$ and that $i(e_{-}) \geq 2q$.
    Since Lemma \ref{3lemme1} states that the asymptotic behaviors of all
$v_k(t)$ are different, there are at most
$i(e_-) - (2q -1)=i(e_-)-i(e_+)$ possible asymptotic behaviors.

Assume now that $i(e_+)=2q\neq 0$ is even.  By Proposition
\ref{spectreLe}, this means that the pair of eigenvalues
$(\lambda_{2q-1}^+,\lambda_{2q}^+)$ is a pair of simple real
eigenvalues satisfying $\lambda_{2q}^+<0<\lambda_{2q-1}^+$.  All
eigenfunctions corresponding to the eigenvalues of $L_{e_+}$ with
negative real part have at least $2q+2$ zeros, except $\varphi^+_{2q}$
which has $2q$ zeros.  Arguing as above, we obtain that $i^-_k \geq
2q-1$ and that $i(e_-) \geq 2q$, that is, that there are at most
$i(e_-) -(2q-1) =i(e_-)-i(e_+)+1$ possible asymptotic behaviors
for the functions $v_k(t)$ when $t\rightarrow -\infty$.  Here is the
point where the fact that $i(e_-) > i(e_+)$ is crucial.  Indeed, this
assumption implies that $\lambda_{2q}^-$ and $\lambda_{2q-1}^-$ have
both positive real parts.  Assume that $p=i(e_-)-i(e_+)+1$.  Then,
according to Lemma \ref{3lemme1}, all the possible asymptotic
behaviors corresponding to the eigenvalues
$\lambda_{2q-1}^-,\ldots,\lambda_{i(e_-)-1}^-$ are taken by the
functions $v_k(t)$.  Without loss of generality, we may assume that
$v_1(t)$ and $v_2(t)$ have an asymptotic behavior corresponding to
$\psi_{2q}^{-}(t)$ and $\psi_{2q -1}^{-}(t)$.
Since, due to  Theorem \ref{th-zeros}, the number of zeros of
$v_k(t)$ is non increasing in time,
the asymptotic behaviors
of $v_1(t)$ and $v_2(t)$, when
$t\rightarrow +\infty$, correspond necessarily to $\lambda^+_{2q}$.
Since $\lambda^+_{2q}$ is a simple eigenvalue, we can find
$\alpha\in \Rm \setminus \{0\}$ such that $v_1(t)+\alpha
v_2(t)=o(e^{\lambda^+_{2q}t})$ when $t\rightarrow +\infty$.
Therefore, by Corollary \ref{CappD1}, $v_1(t)+\alpha v_2(t)$ has
at least $2q+2$ zeros when $t$ tends to infinity.
However, when $t\rightarrow
-\infty$, the asymptotic behavior of $v_1(t) + \alpha v_2(t)$ is
determined by $\psi_{2q}^{-} + \alpha \psi_{2q -1}^{-}$, which does
not identically vanish by Lemma  \ref{3lemme1}. By  Proposition
\ref{spectreLe}, $\psi_{2q}^{-} + \alpha \psi_{2q -1}^{-}$ has
exactly $2q$ zeros. Thus, $v_1(t ) +\alpha
v_2(t)$ has exactly $2q$ zeros for $t$ close to $-\infty$ and we get a
contradiction
with the lap number property stated in Theorem \ref{th-zeros}. As a
consequence, one of the asymptotic behaviors corresponding to
$\psi_{2q}^{-} $ or
$\psi_{2q -1}^{-}$ is not realized by the family $(v_k(t))_{1 \leq k
\leq p}$. Thus,
$p\leq i(e_-)-i(e_+)$ and
\eqref{dimInter} holds.

This concludes the proof of Assertion 1) of Theorem \ref{th1}. \hfill
$\square$\\ \vspace{0.4cm}

\subsection{Connections involving periodic orbits}\label{Sequper2}

This section is devoted to the proof of the second statement of
Theorem \ref{th1}. The proof follows the same lines as the ones of the
proof of the first part of Theorem \ref{th1}. As we shall see, the
proof is even simpler because of the presence of the eigenvalue
$\mu=1$, which implies for example that the dimension of the unstable
manifold of a hyperbolic periodic orbit is larger than its Morse
index. In particular, the difficulties encountered in the second
part of the proof in Section 4.1 (see the case $i(e_{+}=2q$) do not
occur.

\begin{theorem} \label{4e-G+}
1) Let $e^{-}$ (resp. $\Gamma^{+}$) be a hyperbolic equilibrium point (resp.
hyperbolic periodic orbit of period $p^{+} >0$). Then the unstable 
manifold $W^u(e^{-})$ intersects
transversally the local stable manifold $W^s_{loc}(\Gamma^{+})$.\\
2) Let $\Gamma^{-}$ (resp. $e^{+}$) be a hyperbolic
periodic orbit of period $p^{-} >0$ (resp. hyperbolic equilibrium 
point). Then the unstable manifold
$W^u(\Gamma^{-})$ intersects
transversally the local stable manifold $W^s_{loc}(e^{+})$.
\end{theorem}

\begin{demo} We proceed
as in Section \ref{Sheteind}. The main ideas behind this proof are
the same as those of the proof of the transversality orbits connecting
two hyperbolic periodic orbits of \cite{CR}.

Assume that $W^u(e^-)\cap W^s_{loc}(\Gamma^+)\neq \emptyset$ (otherwise the
intersection is transversal by definition). Let $u(t)$ be a global
solution of \eqref{eq} with $u(0)\in W^u(e^-)\cap W^s_{loc}(\Gamma^+)$. By
definition $W^u(e^-)$ intersects
$W^s_{loc}(\Gamma^+)$ transversally if and only if $T_{u(0)}W^u(e^-) +
T_{u(0)}W^s_{loc}(\Gamma^+)=H^s(S^1)$. We know that
$\dim(T_{u(0)}W^u(e^-))=i(e^-)$ and
$\codim(T_{u(0)}W^s_{loc}(\Gamma^+))=i(\Gamma^+)$. Thus, $T_{u(0)}W^u(e^-) +
T_{u(0)}W^s_{loc}(\Gamma^+)=H^s(S^1)$ if and only if
\begin{equation}
\label{dInter2}
\dim\left(T_{u(0)}W^u(e^-)\cap T_{u(0)}W^s_{loc}(\Gamma^+)\right) \leq i(e^-) -
i(\Gamma^+) ~.
\end{equation}
Since $W^u(e_-)$ is a finite-dimensional manifold, there exists a finite
basis $(v_1^0,\ldots,v_p^0)$ of $T_{u(0)}W^u(e_-)\cap
T_{u(0)}W^s_{loc}(e_+)$. Let $v_k(t)$ be the global solutions of Eq.
\ref{3eqvk}, that is,
$$
\partial_t v_k=\partial^2_{xx}v_k+ D_uf(x,u,u_x)v_k+ D_{u_x}f(x,u,u_x)
\partial_x v_k~, \quad v_k(0)=v_k^0~.
$$
The basis $(v_1^0,\ldots,v_p^0)$ of
$T_{u(0)}W^u(e^-)\cap T_{u(0)}W^s_{loc}(\Gamma^+)$ has the asymptotic
behavior  described by (i)-(iv) at the beginning of Section
\ref{Sheteind}, when $t$ goes to $-\infty$. Without loss of
generality, we may assume that this basis has been chosen as in
Lemma \ref{3lemme1}. By Corollary \ref{CappD1}, we know that the
number of zeros of $v_k(t)$ is at most equal to the Morse index
$i(e^{-})$. Let $\Gamma^{+}=\{\gamma^{+}(t) | t \in [0,p^{+})\}$.
Since $u(0)$ belongs to $W^s_{loc}(\Gamma^+)$, there exists $a^{+}
\in [0,p^{+})$ such that $u(0)$ belongs to the local strongly stable
manifold of the point $\gamma^{+}(a^{+})$. Thus,
when $t$ goes to $+\infty$, the asymptotic
behavior of the function $v_k(t)$ is given by Corollary \ref{CappD2} and
corresponds to eigenvalues $\mu^+_{i^+_k}$ of the period map
$\Pi^+(p^{+} +a^{+},a^{+})$ (which coincide with the eigenvalues of the period
map $\Pi^+(p^{+}, 0)$, with
$|\mu^{+}_{i^+_k}|<1$. Hence, $i^+_k>i(\Gamma^+)+1$. Furthermore, if
$i^+_k=2j$ or $2j-1$, then, when $t$ is large enough, $v_k(t)$ has exactly $2j$
zeros which are all simple.

If $i(\Gamma^+)=0$, then $T_{u(0)}W^s_{loc}(\Gamma^+)=H^s(S^1)$ and the
transversality trivially holds. We notice that however this situation
does not arise, since, as poved by \cite{Hirsch}, periodic orbits
are never stable in the case of Eq. \eqref{eq}.

Assume that $i(\Gamma^+)=2q^{+}-1$ is odd, then for $t$
large enough, each function $v_k(t)$ has at least $2q^{+}$ zeros.
Since by Theorem \ref{th-zeros} the
number of zeros of $v_k(t)$ is non-increasing,
$v_k(t)$ has at least $2q^{+}$ zeros for every time $t\in\Rm$.  Applying
    Corollary \ref{CappD1}, we obtain that necessarily
$i^-_k \geq 2q^{+} -1$ and that $i(e_{-}) \geq 2q^{+}$ (we already
know this property by Corollary \ref{i+i-}).
Since Lemma \ref{3lemme1} states that the asymptotic behaviors of all
$v_k(t)$ are different, there are at most
$i(e^-) - (2q^{+} -1)=i(e^-)-i(\Gamma^+)$ possible asymptotic
behaviors when $t$ goes to $-\infty$.

Assume that $i(\Gamma^+)=2q^{+}$ is even, then, for $t$
large enough, each function $v_k(t)$ has at least $2q^{+}+2$ zeros
and, arguing as above, we obtain that
there are at most $i(e^-) - (2q^{+} +1) \leq i(e^-)-i(\Gamma^+)$ possible
asymptotic behaviors when $t$ goes to $-\infty$.

In each case, \eqref{dInter2} is satisfied and the heteroclinic orbit $u(t)$
is transverse.

\vspace{2mm}

The proof of the second statement is similar to the one of statement
1) and also to the proof of the first part of Theorem \ref{th1} given
in Section 4.1. Notice that we do not encounter the difficulty, which
arises in Section 4.1 in the case where $i(e^{+}) =2q$ and
$i(\Gamma^-)=i(e^{+}) =2q$, since we have
an additional dimension at our disposal. Indeed, if $\Gamma^-$ is a 
hyperbolic periodic orbit, then
$\dim(W^u(\Gamma^-)))=i(\Gamma^-) +1$. Therefore, $W^u(\Gamma^-)$ intersects
$W^s_{loc}(e^+)$ transversally at $u(0)$ if and only if
\begin{equation}
\label{dInter3}
\dim\left(T_{u(0)}W^u(\Gamma^-)\cap T_{u(0)}W^s_{loc}(e^+)\right)\leq
i(\Gamma^-) +1  - i(e^+) ~.
\end{equation}
Thus, arguing as in Section \ref{Sheteind}, we prove that
the intersection $W^u(\Gamma^-)\cap W^s(e^+)$ is transversal, even
if $i(\Gamma^-)=i(e^+)$.
\end{demo}


\section{Generic non-existence of homoindexed orbits}\label{Shomoind}

In the previous sections, we have seen that the unstable manifolds of
hyperbolic periodic orbits always intersect transversally the local
stable manifolds of hyperbolic periodic orbits or equilibrium points.
Likewise, the unstable manifolds of hyperbolic equilibrium points
always intersect transversally the local stable manifolds of
hyperbolic periodic orbits.  In Section \ref{Sheteind}, we have proved
that any orbit connecting two hyperbolic equilibrium points $e_{-}$
and $e_{+}$ of different Morse indices $i(e_{-})$ and $i(e_{+})$ is
transverse.  In Section 3.1, we have seen that there does not exist
any connecting orbit between two hyperbolic equilibrium points with same odd
index $i(e_{-})=i( e_{+})= 2m-1$, $m \geq 1$.  Thus, in this section,
it remains to show that generically with respect to the non-linearity
$f$, there does not exist any orbit connecting two equilibrium points
with same even index $i(e_{-})= i(e_{+})= 2m$, $m \geq 1$.  In the
proof of the generic non-existence of homoindexed orbits, we will
actually show that generically with respect to $f$, all the connecting
orbits, connecting equilibria with equal even Morse index, are
transverse, which precludes the existence of homoindexed orbits.  To
show this genericity result, we we shall use a functional
characterization of the transversality of all connecting orbits
$\mathcal{C}(e^{-}(f), e^{+}(f))$ of \eqref{eq} and apply the
Sard-Smale theorem.

\subsection{Preliminaries}

In the introduction, we have assumed that the conditions
\eqref{condf1} hold, which imply that Eq. \eqref{eq} defines
a global dynamical system $S_{f}(t)$ in $H^s(S^1)$ given by
$S_f(t)u_0=u(t)$, where $u(t) \in C^1(\Rm^{+}, H^s(S^1))$ is the
(classical) solution of \eqref{eq}. If the conditions
\eqref{condf1} do not hold, then $S_f(t)$ is only a local dynamical
system. At the end of the introduction, we have remarked that the
automatic transversality as well as the generic transversality
properties are still true, even if Hypothesis \eqref{condf1} does
no longer hold. For this reason, in this section, we do not take into
account this hypothesis.

We recall that $\Gm$ denotes the space $C^2(S^1 \times \Rm \times \Rm, \Rm)$
endowed with
Whitney topology (see \eqref{Whitney}). We fix a
non-linearity $f_0$ in $\Gm$ (satisfying or not
Hypothesis \eqref{condf1}). We assume that
$f_0$ is chosen so that all the equilibria and periodic orbits of the
corresponding equation \eqref{eq} are hyperbolic.
We also consider the set $\Cm_0(e^{-}_0,e^{+}_0) \equiv
\Cm_{f_0}(e^{-}_0,e^{+}_0)$ of all the
orbits $u(t)=S_{f_0}(t) u_0$ of \eqref{eq}, connecting
two (hyperbolic) equilibria $e^{\pm}_0$. Remark that $e^{+}_0$ could be
equal to $e^{-}_0$.

In the next section 5.3, we shall give a functional characterization of
the transversality of the connecting orbits $\Cm(e^{-}(f),e^{+}(f))$ of
\eqref{eq} for $t \in \Rm$, which connect equilibria $e^{-}(f)$ and $e^{+}(f)$,
close to $e^{-}_0$ and $e^{+}_0$, when $f$ belongs to a small enough
neighbourhood of $f_0$ in $\Gm$. We will show that, even if
$\Cm_0(e^{-}_0,e^{+}_0)$ is not a transverse connecting orbit, we can
find $f$ as close to $f_0$ as is wanted so that all the connecting orbits
$\Cm(e^{-}(f),e^{+}(f))$ (with norm less than a given constant) are transverse.

Since $\Gm$ is not metrizable and the classical perturbation
theorems are usually proved in Banach spaces, we will ``replace"
$\Gm$ by a Banach space in the following way. \HB
Since $H^s(S^1)$ is continuously embedded in $C^1(S^1)$, there exists 
a positive integer
$k_0 >1$ such that $\| v\|_{C^1} \leq k_0 \| v\|_{H^s}$, for any $v \in
H^s(S^1)$. Now, for any $M_0$, we
introduce the restriction operator
$R(M_0): g \in \Gm \mapsto Rg \in C^2(S^1 \times [-k_0(M_0+2),k_0 (M_0+2)]
\times [- k_0(M_0+2), k_0(M_0+2)], \Rm)$ defined by
$$
R(M_0) g= g_{|{S^1 \times [-k_0(M_0+2), k_0(M_0+2)] \times [-k_0(M_0+2),
k_0(M_0+2)]}}~.
$$
The map $R(M_0)$ is continuous, open and surjective from $\Gm$ into
$R(M_0)\Gm$.

In what follows, we need the following two auxiliary lemmas. The
first lemma allows to construct appropriate neighborhoods of the
equilibria $e^{\pm}_0$ in $H^s(S^1)$, when $f$ is close to $f_0$ in
a small enough neighborhood of $f_0$ in $C^2$. This lemma is 
classical and is proved
as  \cite[Lemma 4.c.2]{BrunPola} (see also
\cite[Lemma 4.10 ]{BrunRaug}). Its proof mainly uses
the continuous dependence of the equilibria and local unstable or 
stable manifolds
with respect to the non-linearity $f$.

\begin{lemma} \label{Bunifvois} Let $M_0$ be a given positive
constant and $f_0 \in \Gm$ be given such that
all its equilibrium points are hyperbolic. Then, $f_0$ has a finite 
number of equilibria $e_j$,
$1\leq j \leq N_0$ such that $\| e_j\|_{H^s}\leq M_0$.
There exist $r_0 >0$, $R_0 >0$, $R_1 >0$, with $r_0
    < R_0 < R_1$, and a
small neighbourhood $\Vm(f_0) \equiv \Vm(f_0, M_0)$ of $f_0$ in
$R(M_0)\Gm$, depending only on $f_0$ and $M_0$, such that the
following properties hold:
\begin{enumerate}
\item   For any $f \in \Vm(f_0, M_0)$ and any $j$, $1 \leq j \leq
N_0$, there exists an
equilibrium point $e_j(f)$ of $S_f(t)$ in $B_{H^s}(e_j(f_0), r_0)$.  The
equilibrium $e_j(f)$ is unique in the closed ball
$B_{H^s}(e_i(f_0),R_1))$ and
has the same Morse index as $e_j(f_0)$.

\item $R_1$ can be chosen so that $B_{H^s}(e_i(f_0),R_1)) \cap
B_{H^s}(e_j(f_0),R_1)) = \emptyset$, if $i \ne j$.

\item There exist small neighbourhoods $\Nc_j(f)$ of $e_j(f)$,
(with $B_{H^s}(e_j(f_0),R_0) \subset \Nc_j(f) \subset
B_{H^s}(e_j(f_0),R_1))$, which  converge to
$\Nc_j(f_0)$ in $H^s(S^1)$ as $f$ converges to $f_0$ in $R(M_0)\Gm$
and satisfy the following property: \HB
the  local stable set $W^s_{loc,f}(e_j(f),\Nc_j(f))$
and the local unstable set $W^u_{loc,f}(e_j(f),\Nc_j(f))$ are
$C^1$-manifolds of codimension $i(e(f_0))$ and  dimension $i(e(f_0))$
respectively. Moreover, $W^u_{loc,f}(e_j(f),\Nc_j(f)) \cap
W^s_{loc,f}(e_j(f),\Nc_j(f)) = \{e_j(f)\}$.

\item If $u(t)=S_{f}(t) u^*$ is a solution of \eqref{eq} such that
$u(t)$ belongs to $B_{H^s}(e_j(f_0), R_0)$ for all $t \geq t_0$
(respectively, to $B_{H^s}(e_j(f_0), R_0)$ for all $t \leq t_1$), then
$u(t)$, $t \geq t_0$ belongs to the local stable manifold
$W^s_{loc,f}(e_j(f),\Nc_j(f))$ (respectively $u(t)$, $t \leq t_1$
belongs to the local unstable manifold
$W^u_{loc,f}(e_j(f),\Nc_j(f))$).

\item Moreover, there exists a positive constant $c_0$ such that, if
$u_1$, $u_2$ are two solutions of \eqref{eq} with $f=f_1$ and $f=f_2$,
where $f_i$, $i=1,2$ belong to $\Vm(f_0, M_0)$, and if $u_1(t)$, $u_2(t)$
belong to $B_{H^s}(e_j(f_0), R_0)$ for all $t \in J$, where $J =
(-\infty, t_0]$ or $J= [t_0,+\infty)$ for some $t_0 \in \Rm$, then,
\begin{equation}
\label{Bu1u2}
\sup_{t \in J} \| u_1 (t)-u_2(t)\|_{H^s} \leq c_0 \big( \| f_1
-f_2\|_{C^1} + \| u_1 (t_0)-u_2(t_0)\|_{H^s}\big)
\end{equation}
\end{enumerate}
\end{lemma}

We also need the following auxiliary lemma about convergence of
connecting orbits.  Its proof is the same as the
proofs of \cite[Lemma 4.c.3]{BrunPola} and
\cite[Lemma 4.11 ]{BrunRaug}. See the proof of Lemma \ref{lnoerrant2} for
similar arguments.

\begin{lemma} \label{5orbitn}
Let $f_0$, $M_0$ and $r_0   <R_0$ be as in Lemma\ref{Bunifvois}.
Let $e^{-}(f_0)$ and $e^{+}(f_0)$ be two (hyperbolic) equilibria of $f_0$
satisfying the conditions of Lemma \ref{Bunifvois}. Let $\rho_0
>0$  be any positive number such that $r_0 < \rho_0 < R_0$.

Let $f_{\nu} \in {\cal V}(f_0)$ be a sequence converging in $R(M_0)\Gm$ to
some function
$f_{\infty} \in {\cal V}(f_0)$. Assume that, for $\nu= 1,2, \ldots$,
$u_{\nu}$ is a solution of \eqref{eq} for $f=f_{\nu}$ such that,
   \begin{equation}
\label{5condn}
\begin{split}
& u_{\nu}(t) \in B_{H^s}(0, M_0)~, \quad \forall t \in \Rm \cr
& u_{\nu}(t) \in B_{H^s}(e^{-}(f_0)), \rho_0)~, \quad \forall t \in
(-\infty, -t_0]\cr
&  u_{\nu}(t) \in B_{H^s}(e^{+}(f_0)), \rho_0)~, \quad \forall t \in
[t_0, \infty)~,
\end{split}
\end{equation}
where $t_0$ is a positive time. If $e^{-}(f_0) =e^{+}(f_0)=e(f_0)$,
we assume in addition  that there exists a
sequence of times $t_{\nu} \in (-t_0,t_0)$ such that
$u_{\nu}(t_{\nu}) \notin B_{H^s}(e(f_0), R_0)$. \HB
Then, $u_{\nu}$ admits a subsequence $u_{\nu_j}$ that converges in
$C^0_b(\Rm,H^s(S^1))$ to a non trivial connecting orbit $u_{\infty}$
of \eqref{eq} for $f=f_{\infty}$, connecting the equilibria
$e^{-}(f_{\infty})$ and $e^{+}(f_{\infty})$.
\end{lemma}

\begin{rem} If,  in the case where $e^{-}(f_0) =e^{+}(f_0)=e(f_0)$, we
do not require that there exists a
sequence of times $t_{\nu} \in (-t_0,t_0)$ such that
$u_{\nu}(t_{\nu}) \notin B_{H^s}(e(f_0), R_0)$, then the subsequence
$u_{\nu_j}$ could converge in $C^0_b(\Rm,H^s(S^1))$ to an equilibrium
point $e(f_{\infty})$ of $S_{f_{\infty}}(t)$.
\end{rem}


In \cite{BrunPola}, in order to give a functional characterization of
the transversality of the connecting orbits $\Cm (e^{-}, e^{+})$,
Brunovsk\'y and Pol\'a\v{c}ik have introduced a functional defined on
a subspace $\Em$ of the continuous bounded mappings from $\Rm$ into
$L^2$.  In our situation, following the path of \cite{BrunPola}, we
could introduce the spaces,
$$
\Em = C^{1,\delta}(\Rm, L^2(S^1)) \cap C^{0,\delta}(\Rm, H^2(S^1))~,
\quad \delta >0~, \quad \Zm = C^{0,\delta}(\Rm, L^2(S^1))~.
$$
Then, we would fix a non-linearity $f_0$, the equilibrium points of
which are all hyperbolic, and fix two such equilibria $e_0^{-}$ and
$e_0^{+}$. By Lemma \ref{Bunifvois}, if a solution
$u(t)= S_{f_0}(t)
u_0$ belongs to the ball $B_{H^s}(e_0^{+}(f_0),R_0)$ for every time
$t \geq t_0$ and to the ball $B_{H^s}(e_0^{-}(f_0),R_0)$ for every
time $t \leq -t_0$,
where $t_0>0$, then $u(t)$
is a connecting orbit from $e_0^{-}$ to $e_0^{+}$. This leads us to
introduce the open subset $\Um_0 \subset \Em$,
\begin{equation}
\label{BUmcont}
\begin{split}
\Um_0 \equiv{\cal U}_0(e_0^{-}, e_0^{+})
= \{ w(t) \in \Em \, | \, & w(t) \in  B_{H^s}(e_0^{+}(f_0),R_0)
\hbox{ for every }t \geq t_0~, \cr
\quad & w(t) \in  B_{H^s}(e_0^{-}(f_0),R_0)
\hbox{ for every }t \leq -t_0\}~.
\end{split}
\end{equation}
We could finally define the functional $\Phi(w,f): (w,f) \in \Um_0 \times
\Vm(f_0) \mapsto \Phi(w,f) \in \Zm$ by
\begin{equation}
\label{PhiBr}
\Phi(w,f) = w_t(x,t) -w_{xx}(x,t)+f(x,w(x,t),w_x(x,t))~.
\end{equation}
As in  \cite[Lemma 4.b.5 and
Corollary 4.b.6]{BrunPola}, we could
show that, if $(u,f)$ belongs to $\Phi^{-1}(0) \cap (\Um_0 \times
\Vm(f_0))$, then the linearized operator $D_{u}\Phi(u,f)$ is a
Fredholm operator of index $i(e_0^{-}) - i(e_0^{+})$. Moreover, we
could show that, if $0$ is a regular value of the map $u \in \Um_0
\mapsto \Phi(u,f)$, then all the connecting orbits $\tilde u(t)$ such that
$(\tilde u,f) \in \Um_0 \times \Vm(f_0))$ are transverse. This is a consequence
of a functional characterisation of the transversality similar to the
one of Appendix \ref{AppendixB}.  Then, we would apply the Sard-Smale
theorem (Theorem \ref{SardSmale}) to the function $\Phi$ to deduce
that, generically with respect to $f$, $0$ is a regular value of the
map $u \in \Um_0 \mapsto \Phi(u,f)$.

\vskip 2mm

However, we have seen, in \cite{BrunRaug}, that it is more convenient
to use a {\sl discretized} version of the functional $\Phi$, that is,
to work with bounded sequences $(w(n\tau))_{n \in \Nm}$ rather than
with bounded continuous mappings $w(t)$.  In the next section, as in
\cite{BrunRaug}, we shall discretize the time variable and replace the
functional $\Phi(\cdot)$ defined on bounded functions on $\Rm$ by a
discrete analog, defined on bounded sequences.

\subsection{Proof of Theorem \ref{th2}}

As already explained, the proof of Theorem \ref{th2} essentially
consists in using the (discrete) functional characterisation of the
transversality given in Appendix \ref{AppendixB} and in applying the
Sard-Smale theorem to an appropriate discretization of the functional
\eqref{PhiBr}. The application of the Sard-Smale theorem involves
some technical difficulties. The way to overcome them is now well
understood (see \cite{BrunPola} and \cite{BrunRaug}). The
verification of the surjectivity of the functional $\Phi$
is the crucial point in the application of the Sard-Smale theorem.

The proof of Theorem \ref{th2} can be decomposed into several steps.

\vspace{3mm}

\noindent {\bf Step 1: Choice of particular neighborhoods and 
reduction to a simpler problem}
\vskip 2mm
\noindent We introduce the sequence of bounded open sets, $m \in \Nm$, given by
$$
B_{m} = B_{H^s}(0,m)=\{ v \in H^s(S^1) \, | \, \| v\|_{H^s} <m \}~.
$$
Since $H^s(S^1)= \cup_{m} \overline B_{m}$, Theorem \ref{th2} will be proved if
we show that, for each $m$, there exists a generic set in $\Gm$, such
that, for any $f$ in this generic set, any orbit $\tilde u(t)$ of \eqref{eq},
connecting two (hyperbolic) equilibria and satisfying $\tilde u(t) \in
{\overline B}_m$, $t \in \Rm$, is transverse. We recall that, since 
$H^s(S^1)$ is
continuously embedded in $C^1(S^1)$, there exists a positive integer
$k_0$ such that $\| v\|_{C^1} \leq k_0 \| v\|_{H^s}$, for any $v \in
H^s(S^1)$. In  \cite[Proposition 3.2]{JR}, we have shown that the set
$$
\Om^h_{m} \, = \, \{ f \in \Gm \, | \, \hbox{ any equilibrium }e
\hbox{ of
\eqref{eq} with } \| e\|_{C^1(S^1)} \leq k_0 (m +1) \hbox{ is
hyperbolic}\}
$$
is open and dense in  $\Gm$.\\
As in Section 5.1, we want to work in subspaces of $C^2(S^1 \times
[-k_0(m+2), k_0(m+2)] \times [-k_0(m+2), k_0(m+2)], \Rm)$ and hence,
we use the restriction operator $R(m)$ that we simply denote $R$.
We set ${\cal RO}_m= R (\Om^h_{m})$ endowed with the topology of
$C^2(S^1 \times [-k_0(m+2), k_0(m+2)] \times [-k_0(m+2), k_0(m+2)],
\Rm)$, which is a separable Banach space.  The set ${\cal RO}_m$ is an
open subset of $C^2(S^1 \times [-k_0(m+2), k_0(m+2)] \times
[-k_0(m+2), k_0(m+2)], \Rm)$ and the map $R$ is continous, open and
surjective.

As already remarked in \cite[Proof of Theorem 4.c.1, p.
165]{BrunPola} (and also in \cite[Proposition 4.12]{BrunRaug}),
Theorem \ref{th2} will be proved by using the following proposition.

\begin{prop}\label{th2Reduc}
Assume that, for any $m\in\Nm$ and any $f_0 \in {\cal RO}_m$,
there exist a small neighbourhood ${\cal V}_{f_0}$ of $f_0$ in ${\cal 
RO}_m$ (or simply
in $C^2(S^1 \times [-k_0(m+2), k_0(m+2)] \times [-k_0(m+2),
k_0(m+2)], \Rm)$) and a generic set ${\cal G}_{f_0,m}$ in ${\cal V}_{f_0}$ such
that, for any $f \in {\cal G}_{f_0,m}$, any solution $\tilde u(t)$ of 
\eqref{eq},
connecting two (hyperbolic) equilibria and satisfying $\|\tilde 
u(t)\|_{H^s}\leq
m$, for any $t \in \Rm$, is transverse. Then Theorem \ref{th2} holds.
\end{prop}
\begin{demo}
Let $m$ be given. Since ${\cal RO}_m$ is separable, there exists a countable
set of functions $f_i$, ${i\in\Nm}$, such that the family of corresponding
neighbourhoods $({\cal V}_{f_i})$, ${i\in\Nm}$, covers ${\cal RO}_m$. Let
$({\cal G}_{f_i,m})$, ${i\in\Nm}$, be the corresponding generic sets and
et $\tilde{\cal G}_{f_i,m}={\cal
G}_{f_i,m}\cup ({\cal RO}_m\setminus {\cal V}_{f_i})$, which is a 
generic subset
of ${\cal RO}_m$. The set ${\cal G}_m=\cap_{i\in\Nm}\tilde{\cal G}_{f_i,m}$ is
generic in ${\cal RO}_m$. Moreover, for any  $f \in {\cal G}_m$, any
solution $\tilde u(t)$ of \eqref{eq}, connecting two (hyperbolic) 
equilibria and
satisfying $\|\tilde u(t)\|_{H^s}\leq m$, $t \in \Rm$, is transverse.
  Since the map $R$ is continous, open and surjective,
  $R^{-1}({\cal G}_m)$ and $R^{-1}({\cal G}_m) \cap \Om_h$
  (where $\Om_h$ is the generic set introduced in
Theorem \ref{th-JR}) are generic subsets of $\Gm$. Finally, we notice
that ${\cal O}_M=\cap_{m\in\Nm} R^{-1}({\cal G}_m) \cap \Om_h$ is the
generic set given in Theorem \ref{th2}.
\end{demo}

The interest of Proposition \ref{th2Reduc} is that we can now work in
a small neighborhood ${\cal V}_{f_0}$ of $f_0$ in ${\cal RO}_m$
insted of working in ${\cal RO}_m$. This neighborhood can be chosen
as small as is needed.

{}From now on, we fix $m\in\Nm$ and $f_0 \in {\cal RO}_m$. We apply
  Lemma \ref{Bunifvois} with $M_0=m$. Hence, there exists a small
  neigborhood $\Vm(f_0)$ of $f_0$ in ${\cal RO}_m$ such that all the
  properties described in Lemma \ref{Bunifvois} are satisfied. In
  particular, let $e_j(f_0)$, $1 \leq
j \leq N_0$ be the (hyperbolic) equilibrium points of $S_{f_0}(t)$ such that
$\| e_j(f_0) \|_{H^s} \leq m$.
  Since $S_{f_0}(t)$ has also a finite
number of (hyperbolic) equilibrium points in the closed ball
$B_{H^s}(0,m+1)$, there exist two real numbers $k$ and $k_1$ such that $1<k_1
<k< k_0 (m+2)/(m+1)$ and that $S_{f_0}(t)$ has only $N_0$ equilibrium points
in the closed ball $B_{H^s}(0,km)$. Moreover, we may choose the neighbourhood
$\Vm(f_0)$ of $f_0$ small enough so that, for any $f \in \Vm(f_0)$,
$S_{f}(t)$ has only $N_0$ equilibrium points in the closed balls 
$B_{H^s}(0,km)$
and has no equilibrium points in $\{ v \in H^s(S^1) \, | \, k_1m \leq
\| v\|_{H^s} \leq km  \}$. Let $r_0$ and $R_0$ be chosen as in Lemma
\ref{Bunifvois}. We notice that $r_0$ and $R_0$ can be chosen small enough so
that, for any $1 \leq j \leq N_0$, $B_{H^s}(e_j(f_0),R_0) \subset
B_{H^s}(0,km)$. For later use, we also fix $\rho_0$ and $\rho_1$ such that
$r_0<\rho_0<\rho_1<R_0$.

For any (even) integer $d$, we introduce the set
${\cal E}_{d}$ of equilibria of  Eq.  \eqref{eq} for
$f=f_0$ in $B_{H^s}(0, m)$, the Morse indices of which  are equal to
$d$ and we set
$$
{\cal D}_{d} \, = \,  \cup_{e_j(f_0) \in {\cal E}_{d}}
B_{H^s}(e_j(f_0),r_0)~.
$$
For any integer $\ell$ and (even) integer $d$, we denote
by ${\cal G}_{m}^{\ell,d}$ the set of all $f \in {\cal V}(f_0)$ that have
the following property: every  solution $u$ of
\eqref{eq}  satisfying
\begin{equation}
\label{Gmnd}
\begin{split}
& u(t) \in \overline B_m~, \quad \forall t \in \Rm \cr
& u(t) \in {\cal D}_{d}~, \quad \forall t \in
(-\infty, -\ell] \cup [\ell, +\infty) ~,\cr
\end{split}
\end{equation}
is transverse. We notice that, due to Lemma \ref{Bunifvois},
any non-trivial orbit satisfying the
conditions \eqref{Gmnd} is a connecting orbit, connecting two
equilibria contained in ${\cal D}_{d}$. Due to the choice of ${\cal
V}(f_0)$, the set
$$
{\cal G}_{m} \, = \, \cap_{\ell, d} {\cal G}_{m}^{\ell,d}
$$
satisfies the transversality conditions stated in Proposition
\ref{th2Reduc}. It is therefore sufficient to prove that each ${\cal
G}_{m}^{\ell,d}$ is open and dense in ${\cal V}(f_0)$.

\vspace{3mm}

\noindent {\bf Step 2: Proof of the openess of ${\cal G}_{m}^{\ell,d}$}
\vskip 1mm

\noindent Assume that $(f_{\nu})$ is a sequence in ${\cal V}(f_0) \setminus
{\cal G}_{m}^{\ell,d}$ which converges to $f_{\infty} \in
{\cal V}(f_0)$. We want to show that $f_{\infty}$ does not belong to
${\cal G}_{m}^{\ell,d}$. Since $f_{\nu} \in {\cal V}(f_0) \setminus
{\cal G}_{m}^{\ell,d}$, there exists a solution of
\eqref{eq} for $f=f_{\nu}$, distinct from any equilibrium and
connecting two equilibria in ${\cal
D}_{d}$. Since there is only
a finite number of sets $B_{H^s}(e_j(f_0) r_0)$ in ${\cal D}_{d}$,
passing to a subsequence, we may suppose that there exist indices
$i_1$, $i_2$, with $1 \leq i_1\leq N_0$, $1 \leq i_2 \leq N_0$ such that
   \begin{equation*}
\begin{split}
u_{\nu}(t) \in  B_{H^s}(e^{i_1}(f_0), r_0)~, \quad \forall t \in
(-\infty,-\ell] , \cr
   u_{\nu}(t) \in  B_{H^s}(e^{i_2}(f_0), r_0)~, \quad \forall t \in
[\ell, + \infty) ~.\cr
\end{split}
\end{equation*}
Moreover, if $i_1=i_2 \equiv i$, since the solution
$u_{\nu}(t)$ is not an equilibrium point, there exists a time
$\tau_{\nu}$ such that $u_{\nu}(\tau_{\nu}) \notin B_{H^s}(e^i(f_0),
R_0)$. By Lemma
\ref{5orbitn}, there exists a subsequence $u_{\nu_j}$ that converges in
$C^0_b(\Rm,H^s(S^1))$ to a non trivial connecting orbit $u_{\infty}$
of \eqref{eq} for $f=f_{\infty}$, connecting the equilibria
$e^{i_1}(f_{\infty})$ and $e^{i_2}(f_{\infty})$. Since this
non-trivial orbit connects two equilibrium points with same Morse
index, it cannot be a transverse orbit, which proves the openess of
${\cal G}_{m}^{\ell,d}$.
\vskip 2mm

It remains to show that ${\cal G}_{m}^{\ell,d}$ is dense in ${\cal
V}(f_0)$. This will be done in the next (and remaining) steps of the
proof, by introducing a discrete version of the functional $\Phi$
described in Section 5.1 and applying the Sard-Smale theorem to it.
To this end, we first need to discretize the semi-flow $S_f(t)$.

\vspace{2mm}

\noindent {\bf Step 3: Discretization of the semi-flow $S_f(t)$}
\vskip 1mm

\noindent
For any $u_0 \in H^s(S^1)$ and $f$ close to $f_0$, we consider the image by the
time $\tau$-map $S_f(\tau)u_0 \in H^s(S^1)$ of $u_0$,
\begin{equation}
\label{AG}
G(u_0) \, := \, G_f(u_0) \, = \, e^{A\tau}u_0 + \int_{0}^{\tau}
e^{A(\tau - \sigma)}f(x,S_f(\sigma)u_0,
\partial_x(S_f(\sigma)u_0))\,d\sigma~,
\end{equation}
where $A= \partial_{xx}$.  As we do not assume global existence of
solutions in this section, $G(u_0)$ may not be defined if $S_f(t)u_0$
blows up in a time shorter than $\tau$.  To overcome this difficulty,
for any given $m>0$, we choose a time $\tau_m>0$ such that $G$ is well
defined for any $u_0\in B_{H^s}(0,m)$ and any $f$ in a neighbourhood
of $f_0$.  Then, since $u\mapsto f(.,u,u_x)$ belongs to ${\cal
C}^r(H^s(S^1), L^2(S^1))$, $r \geq1$, the mapping $G$ belongs to
${\cal C}^r(B_{H^s}(0,m),H^s(S^1))$ and, for any $v_0 \in H^s(S^1)$,
\begin{equation}
\label{ADG}
\begin{split}
DG(u_0)v_0 \, = \, e^{A\tau_m}v_0 + \int_{0}^{\tau_m} &e^{A(\tau_m -
\sigma)} \Big(
D_uf(x, S_f(\tilde{s})u_0,\partial_x(S_f(\sigma)u_0))
((DS_f(\sigma)u_0)v_0) \cr
&+D_{u_x}f(x, S_f(\sigma)u_0,\partial_x(S_f(\sigma)u_0))
((DS_f(\sigma)u_0)v_0)_x \Big)
\, d\sigma~,
\end{split}
\end{equation}
that is, $DG(u_0)v_0$ is the image at the time $t=\tau_m$ of the
(classical) solution of the linearized equation,
\begin{equation*}
\begin{split}
       &\partial_t v  \,=\, Av + D_uf(x,S_f(t)u_0,
       \partial_x(S_f(t)u_0))v + D_{u_x}f(x,S_f(t)u_0,
       \partial_x(S_f(t)u_0))v_x ~,
       t>0, \cr
       &v(0)\,=\, v_0~.
    \end{split}
\end{equation*}
In other words, if $\tilde{u}(t)$ is a bounded orbit of
\eqref{eq} with $\sup_{t \in \Rm} \|\tilde{u}(t)\|_{X} \leq m$, then, 
for any $n
\in \Zr$,
\begin{equation}
\label{5DGT}
G(\tilde{u}(n\tau_m)) \, = \, \tilde{u}((n +1)\tau_m)~, \quad
DG(\tilde{u}(n \tau_m))v_0
\, = \, T_{\tilde{u}}((n+1)\tau_m, n\tau_m)v_0~,
\end{equation}
where $T_{\tilde{u}}(t,s)$ is the evolution operator (on $L^2(S^1)$)
defined by the linearized equation along the bounded orbit ${\tilde u}(t)$ (see
Eq. \eqref{Tutilde} in Appendix B.3).

In the next step, we shall introduce the discretized version of the
functional $\Phi$ defined in \eqref{PhiBr} and the ``discretized"
open set corresponding to ${\cal U}_0$. We require several smallness conditions
on the time step $\tau_m$. We assume that
$\tau\equiv\tau_m\equiv\tau_{m,k,k_1}$ and $\Vm(f_0)$ are small enough such
that:\\
(i) if $\| u_0\|_{H^s}\leq m$, then $S_{f_0}(t)u_0 \in B_{H^s}(0,
k_1m)$, for $0 \leq t \leq \tau_m$,\\
(ii) if $\| u_0\|_{H^s}\leq k_1m$, then $S_{f_0}(t)u_0 \in B_{H^s}(0,
km)$, for $0 \leq t \leq \tau_m$,\\
(iii) if $u_0 \in B_{H^s}(e_j(f_0),\rho_0)$, $1 \leq j \leq N_0$ and $f \in
\Vm(f_0)$, then $S_{f}(t) u_0 \in  B_{H^s}(e_j(f_0), \rho_1)$ for $0
\leq t \leq \tau_0$.\\
With these conditions, we control the behavior to the continuous solution
$S_f(t)u_0$ between two time steps $n\tau$ and $(n+1)\tau$. For example, (iii)
ensures that if $u(n\tau)$ belongs to $B_{H^s}(e_j(f_0),\rho_0)$ for
any large enough $n$, then $u(t)$ belongs to $B_{H^s}(e_j(f_0), R_0)$ for  $t$
large enough and thus $u(t)$ belongs to the local stable manifold of 
$e_j(f_0)$.
\vskip 2mm

\noindent {\bf Step 4: A functional characterisation of the transversality}
\vskip 1mm

We are now ready to define the discretized version of the functional
$\Phi$ introduced in Section 5.1. This follows the lines of \cite{BrunRaug}.

We recall that the integer $m$, the function $f_0$ and the
neighbourhood $\Vm(f_0)$ are fixed. Since the Sard-Smale theorem requires that
$\Phi$ is defined on open sets and  the set $\overline B_m$ used in
the definition of ${\cal G}_{m}^{\ell,d}$ is closed, we need to introduce
the following set
$$
B_m^* \, = \, \{ v \in H^s(S^1) \, | \, \|v \|_{H^s} < k_1m \}~.
$$
Let ${\cal E}_d$ be the set of equilibrium points of $S_{f_0}(t)$
in $\overline B_m$ of Morse index $d$. We set
$$
{\cal D}_d^* \, = \, \cup_{e_j(f_0) \in {\cal E}_d} B_{H^s}(e_j(f_0),\rho_0)~.
$$
For any integer $\ell$ and any (even) integer $d$, we finally
introduce the following subspace of $\ell^{\infty}( \Zr,H^s(S^1))$,
\begin{equation}
\label{Xmld}
{\cal X} \equiv {\cal X}_{m, \ell,d} = \{ w (\cdot \tau) \in
\ell^{\infty}(\Zr,H^s(S^1)) \, | \,\forall |n| \geq \ell,~
   w(n\tau) \in  {\cal D}_d^* \text{ and }\forall n \in \Zr,~
  w(n\tau) \in B_m^*\}~.
\end{equation}
We notice that ${\cal X}$ is open in $\ell^{\infty}( \Zr,H^s(S^1))$
and contains the discretizations of all connecting orbits of $S_f(t)$
satisfying \eqref{Gmnd}.
We next define the discretized map $\Phi\equiv\Phi_{m,\ell,d}:
{\cal X}_{m, \ell,d} \times
\Vm(f_0)  \rightarrow l^{\infty}(\Zr,H^s(S^1))$ by
    \begin{equation}
\label{Phimld}
\Phi(w, f) (n) \, \equiv \, \Phi_{m, \ell,d}(w, f) (n) \, = \,
w((n+1)\tau) -G_f(w(n\tau))~,
\quad \forall n \in
\Zr~,
\end{equation}
where $G_f$ has been defined in \eqref{AG}.

Arguing as in \cite[Section 4.2]{BrunRaug}, we obtain the following
characterization of the transversality. Its proof is based
  on the abstract formulation of
transversality given in Appendix \ref{AppendixB}. Without loss of
generality, we may replace $\Vm (f_0)$ by a smaller  neighborhood and
thus assume that $\Vm (f_0)$ is actually a convex neighborhood.

\begin{theorem} \label{BPhivareg}
The above map $\Phi: {\cal X}_{m, \ell,d} \times
\Vm(f_0)  \rightarrow l^{\infty}(\Zr,H^s(S^1))$
is of class ${\cal C}^1$.
A pair $(u,f)$ belongs to
$\Phi^{-1}(0)$ if and only if $u$ is the
discretization of a connecting orbit ${\tilde u}(t)$ (or an 
equilibrium point) of $S_f(t)$
contained in $B_{H^s}(0,km)$ whose discretization belongs to
$B_{m}^{*}$.
Moreover, for any $(u,f) \in \Phi^-1(0)$  the mapping $D_{u} 
\Phi(u,f)$ is a Fredholm
operator of index $0$. \HB
If $0$ is a regular value of the map $u \in
{\cal X}_{m, \ell,d} \mapsto \Phi(u,f)$, then all the
connecting orbits the discretizations of which are contained in 
${\cal X}_{m, \ell,d}$ are
transverse, i.e. $f\in{\cal G}^{\ell,d}_m$.
\end{theorem}

\begin{demo}
The description of the zeros $(u,f)$ of $\Phi$ is
obvious. Indeed,
$G_f(u(n\tau))=S_f(\tau)u(n\tau)$ and thus $u$ is the discretization of a
trajectory ${\tilde u}(t)$ of $S_f(t)$. Moreover, due to the definition of
${\cal X} \equiv {\cal X}_{m, \ell,d}$ and to Lemma
\ref{Bunifvois}, it is a connection between two equilibrium points
$e^-(f)$ and $e^+(f)$ of same index $d$ and
is a constant
sequence if and only if it coincides with an equilibrium
$e^{-}(f)=e^{+}(f)$.

It is straightforward to prove that the
  mapping $\Phi : {\cal X} \times \Vm \mapsto l^{\infty}(\Zr,H^s(S^1))$
is of class $C^1$ (see \cite[Lemma 4.c.4]{BrunPola} or \cite[Lemma 
4.13]{BrunRaug}).
Moreover, the first derivative $D \Phi (u,f)$, for
$(u,f) \in {\cal X} \times \Vm(f_0)$, is given by
\begin{equation}
\label{DwPhi}
\begin{split}
D \Phi (u,f)(Y,h)(n) &\, = \, Y((n+1)\tau)
- D_uG_f(u(n\tau))) Y(n\tau)
   - D_f G_f(u(n\tau))\cdot h \cr
&\, = \, Y((n+1)\tau)- T_{u,f}((n+1)\tau,n\tau) Y(n \tau) - D_f
G_f(u(n \tau))\cdot h \cr
&\, \equiv \, ({\cal L}_{u,f}Y)(n\tau)- D_f
G_f(u(n\tau))\cdot h~,
\end{split}
\end{equation}
where $(Y,h)$ is any element of $l^{\infty}(\Zr,H^s(S^1)) \times \Gm$
and where $T_{u,f}$, $t \geq s$, is the evolution operator defined by
the linearized equation \eqref{Tutilde} and ${\tilde u}$ is
the solution of \eqref{eq}, the discretization of
which is given by $u$.

The expression \eqref{DwPhi} of the derivative $D\Phi$ shows that $u$ is a
regular zero of the mapping $u \in {\cal X}_{m, \ell,d} \mapsto
\Phi(u,f)$ if and only if the mapping ${\cal L}_{u,f}$ is surjective.
Corollary  \ref{ApFredB} implies that
  the map ${\cal L}_{u,f}$ is surjective if
and only if $\tilde{u}$ is transverse.  Corollary  \ref{ApFredB} also
tells that ${\cal L}_{u,f}$ is a Fredholm operator of
index equal to $i(e^{-}(f)) - i(e^{+}(f))=0$.
\HB
As noticed above, if $\tilde{u}=\{\ldots, e,e, \ldots.,\}$ is a constant
sequence, then, by the construction of the various neighbourhoods
made in Lemma \ref{Bunifvois}, $e$ is a hyperbolic equilibrium point of
\eqref{eq}, which implies that $e = e^{-}(f) =e^{+}(f)$.
Again, by Theorem \ref{AFredholm}, the surjectivity of the map
${\cal L}_{u,f}$ is then equivalent to the hyperbolicity of $e^{-}(f)$.
\end{demo}

\vspace{2mm}

\noindent {\bf Step 5:  Surjectivity of $D\Phi$}
\vskip 1mm

\noindent As already explained, we will apply the Sard-Smale theorem to the
functional $\Phi$ introduced in Step 4 and consider the set
$\Phi^{-1}(0)$ in particular. One of the main hypotheses of the
Sard-Smale is the fact that $0$ is a regular value of the map $(w,f)
\in {\cal X} \times \Vm(f_0)
\mapsto \Phi(w,f)$. This property will be shown in the next theorem
as consequence of Corollary \ref{ApFredB} and the one-to-one property
of homoindexed orbits proved in Proposition \ref{injectif}.

\begin{theorem} \label{BDPhisurj} Assume that ${\cal X}$ and $\Phi$
are given as in Step 4. \HB
1) The pair $(\tilde u,f)$ is a regular zero of the map $(w,f)\in
{\cal X} \times \Vm (f_0)
\mapsto \Phi (w,f)$, if and only if, for any nontrivial
bounded solution $\psi(t) \in C^0_b(\Rm, L^2(S^1))$ of
the adjoint equation \eqref{Tutilde*}, there exists $\tilde{g} \in
R\Gm$ such that
$$
\int_{-\infty}^{+\infty} \langle \psi(t), \tilde{g}(\tilde{u}(t))
\rangle_{L^2(S^1)} \,{\rm
d}t   \, \ne \, 0~,
$$
where $\tilde u(t)=S_f(t)\tilde u(0)$ is the (continuous) trajectory
corresponding to the sequence $(\tilde u(n\tau))$. \HB
2) As a consequence of the first statement and of Proposition
\ref{injectif}, $0$ is a regular value of the map
$\Phi$.
\end{theorem}

\begin{demo}
We first prove the second statement of the theorem, which is a direct
consequence of Proposition \ref{injectif}.  If $\Phi(\tilde u,f)=0$ 
and $\tilde u$ is
the discretization of an (hyperbolic) equilibrium point, then as
explained in the proof of Theorem \ref{BPhivareg}, the map ${\cal
L}_{\tilde u,f}$ is
surjective and thus $D \Phi(\tilde u,f)$ is also surjective. Thus, it
remains to consider the case where $(\tilde u,f) \in \Phi^{-1}(0)$ and $\tilde
u$
is not the discretization of an equilibrium point. By the first
statement, the operator $D\Phi (\tilde u,f) \in L({\cal X} \times
\Vm,
l^{\infty}(\Zr,H^s(S^1))$ is surjective if and only if, for any
nontrivial bounded solution $\psi(t) \in C_b^0(\Rm, L^2(S^1)$ of the
adjoint equation \eqref{Tutilde*}, there exists $\tilde{g} \in
R\Gm$ such that
\begin{equation}
\label{5intnnul}
\int_{S^1} \int_{-\infty}^{+\infty} \psi(x,t), \tilde{g}(x,
\tilde{u}(x,t), \tilde{u}_x(x,t)) dx dt \, \ne \, 0~.
\end{equation}
Since $\psi(t)$ is a non trivial solution of the adjoint equation,
  there exist $x_0 \in S^1$ and $t_0$ such that $\psi(x_0,t_0)
\ne 0$. Due to the injectivity property of Proposition \ref{injectif}, for $x_0$
fixed, there exists no other time $t_1$ such that
$\tilde{u}(x_0,t_1)=\tilde{u}(x_0,t_0)$ and
$\tilde{u}_x(x_0,t_1)=\tilde{u}_x(x_0,t_0)$. Moreover, Proposition
\ref{injectif} also implies that $({\tilde u}(x_0,t), {\tilde u}_x(x_0,t))$
stays outside a small neighbourhood of $({\tilde u}(x_0,t_0), {\tilde
u}_x(x_0,t_0))$ for $t$ close to $\pm\infty$. Therefore, one easily constructs a
regular bump function $ {\tilde g}$ which vanishes outside a small
neighbourhood of $(x_0,{\tilde u}(x_0,t_0), {\tilde u}_x(x_0,t_0))$
and is positive in this neighbourhood; so that the function
$(x,s)\mapsto {\tilde g}(x, \tilde{u}(x,s),\tilde{u}_x(x,s))$ is a
regular bump function concentrated around $(x_0,t_0)$. For such a
choice of ${\tilde g}$, the condition \eqref{5intnnul} is thus satisfied.

\vspace{2mm}

We now prove the first statement of the theorem. This proof is
nothing else as the proof of  \cite[Theorem 4.7]{BrunRaug}. For the reader's
convenience, we reproduce it here.

As already explained,  if $({\tilde u}, f)$ belongs to $\Phi^{-1}(0)$,
then $\tilde u$ is a discretization of a trajectory $\tilde{u}(t)$,
$t \in \Rm$, of \eqref{eq}, connecting two equilibria $e^{-}(f)$ and
$e^{+}(f)$.
Without loss of generality, we may assume that $\tilde{u}$ is
a nonconstant sequence. Indeed, if $\tilde{u}=\{\ldots, e,e,
\ldots.,\}$, then $e= e^{-}(f)=e^{+}(f)$ is a hyperbolic
equilibrium point and so
$\tilde{u}$ is a regular zero of $\Phi$. On
the other hand,
the adjoint equation \eqref{Tutilde*} has no nontrivial bounded
solution.

Thus, we assume that $\tilde{u}$ is not a constant sequence. We
recall that, by \eqref{DwPhi}, for any $(Y, \tilde{g}) \in
l^{\infty}(\Zr,H^s(S^1)) \times R\Gm$
$$
(D \Phi(\tilde{u},f)\cdot (Y,\tilde{g}))(n\tau) \, = \,
({\cal L}_{\tilde u, f}Y)(n\tau) -
D_fG_f(\tilde{u}(n\tau))\cdot \tilde{g}~.
$$
We notice that ${\cal L}_{\tilde{u},f}$ corresponds to the operator ${\cal L}$
defined in \eqref{AcalL}. A sequence $H\in\ell^{\infty}(\Zr,H^s(S^1))$
is in the range of $D\Phi$ if and only if one can choose
$\tilde{g} \in  R\Gm$ such that $H+D_fG_f(\tilde u).\tilde g$ is in 
the range of
${\cal L}_{\tilde u}$. According to Corollary \ref{ApFredB}, this is
equivalent to finding $\tilde g$ such that
    \begin{equation}
\label{BDPHIaux2}
    \sum_{n = -\infty}^{+\infty} \langle \psi((n+1)\tau), D_f
        G_f(\tilde{u}(n\tau))\cdot \tilde{g} \rangle_{L^2(S^1)}\,=\,
        -  \sum_{n = -\infty}^{+\infty} \langle \psi((n+1)\tau,
        H(n\tau)\rangle_{L^2(S^1)}~,
\end{equation}
for every nontrivial sequence $\psi(n\tau)=T^{\ast}(n\tau,0) \psi_0$, $\psi_0
\in L^2(S^1)$, which is bounded in $L^2(S^1)$.  We
can choose such a $\tilde{g} \in R\Gm$ if, given a basis $\psi_1,
\psi_2, \ldots ,
\psi_q$ of the (necessarily) finite-dimensional
vector space of bounded sequences $\psi(n)= T^{\ast}(n\tau,0)
\psi_0$, the mapping
\begin{equation}
\label{BDPHIaux3}
\tilde{g} \in \Gm \mapsto
(\, \sum_{-\infty}^{+\infty} \langle \psi_j((n+1)\tau), D_f
G_f(\tilde{u}(n\tau))\cdot
\tilde{g} \rangle_{L^2(S^1)} \, )_{1 \leq j \leq q} \in \Rm^q~
\end{equation}
is surjective.  If the range of the mapping \eqref{BDPHIaux3} is not the
whole vector space $\Rm^q$, there exists a vector $(\alpha_1,\ldots,\alpha_q)$
orthogonal to the range, that is, there
exists a bounded sequence $\psi =\sum \alpha_j\psi_j \ne 0$ such that,
for any $\tilde{g} \in R\Gm$,
$$
     \sum_{n = -\infty}^{+\infty} \langle \psi((n+1)\tau), D_f
    G_f(\tilde{u}(n\tau))\cdot \tilde{g} \rangle_{L^2(S^1)}\,=\, 0~.
$$
Thus, $D\Phi$ is surjective if and only if, for
any bounded sequence $\psi(n\tau)=T^{\ast}(n\tau,0) \psi_0$, there
exists $\tilde{g}
\in R\Gm$ such that,
    \begin{equation}
\label{BDPHIaux4}
\sum_{n = -\infty}^{+\infty} \langle \psi((n+1)\tau,
       D_f G_f(\tilde{U}(n\tau))\cdot \tilde{g} \rangle_{L^2(S^1}) \, \ne \, 0~.
\end{equation}
Since the solution of \eqref{eq} is differentiable
with respect to $f$, we can differentiate \eqref{eq} formally
with respect to  $f$ to deduce that, for any $\tilde{g}$,
    \begin{equation}
\label{BDPHlaux5}
D_f G_f(\tilde{u}(n\tau))\cdot \tilde{g}  \, = \, \int_{n\tau}^{(n+1)\tau}
T_{\tilde{u}}((n+1)\tau,\sigma) \tilde{g}(\tilde{u}(\sigma))\,d\sigma~.
\end{equation}
Using the expression \eqref{BDPHlaux5} in the condition
\eqref{BDPHIaux4} yields
\begin{equation}
\label{BDPHIaux6}
\begin{split}
\sum_{n = -\infty}^{+\infty}
\langle \psi(n+1),
       D_f G_f(\tilde{u}&(n\tau)) \cdot \tilde{g} \rangle_{L^2(S^1)}  \cr
       &=\sum_{n =
-\infty}^{+\infty} \int_{n\tau}^{(n+1)\tau}
\langle \psi((n+1)\tau), T_{\tilde{u}}((n+1)\tau,\sigma)
\tilde{g}(\tilde{u}(\sigma))d\sigma \rangle_{L^2(S^1)} \cr
       &=\sum_{n = -\infty}^{+\infty} \int_{n\tau}^{(n+1)\tau} \langle
        T_{\tilde{u}}^\ast(\sigma,(n+1)\tau) \psi((n+1)\tau),
        \tilde{g}(\tilde{u}(\sigma)) d\sigma \rangle_{L^2(S^1)} \cr
     &=\int_{ -\infty}^{+\infty} \langle
        \psi(\sigma), \tilde{g}(\tilde{u}(\sigma))d\sigma \rangle_{L^2(S^1)} ~,
\end{split}
\end{equation}
where, for any $\sigma \in \Rm$, $\psi(\sigma) =
T_{\tilde{u}}^\ast(\sigma,(n+1)\tau)
\psi((n+1)\tau)= T_{\tilde{u}}^\ast(\sigma,0) \psi_0$. We remark that
$\psi(\cdot)$ belongs to ${\cal C}^0_b(\Rm,L^2(S^1))$ and is
a bounded solution of \eqref{Tutilde}.
Theorem \ref{BDPhisurj} is thus proved.
\end{demo}

\vskip 2mm

\noindent {\bf Step 6:  Application of the Sard-Smale theorem and density of
${\cal G}_{m}^{\ell,d}$}
\vskip 1mm

\noindent 
After all the preliminaries given in the previous steps, we are now
ready to apply the Sard-Smale theorem (in the form recalled in Appendix
\ref{AppendixA}).
This follows the lines of \cite{BrunPola} and in \cite{BrunRaug}.
\vskip 1mm
We recall that the subspace ${\cal X}$ of
$\ell^{\infty}(\Zr,H^s(S^1))$ has been defined in \eqref{Xmld}. We
introduce the open subset ${\cal Y}= \Vm(f_0)$ of ${\cal RO}_m$ and
the Banach space ${\cal Z}= \ell^{\infty}(\Zr,H^s(S^1))$. \HB
We recall that the mapping $\Phi\equiv \Phi_{m,\ell,d}:
{\cal X} \times {\cal Y} \rightarrow {\cal Z}$ has been given in
\eqref{Phimld} as follows:
$$
\Phi(w, f) (n)  \, = \, w((n+1)\tau) -G_f(w(n\tau))~,
\quad \forall n \in \Zr~,
$$
where $G_f$ has been defined in \eqref{AG}. \HB
We  now check that all the hypotheses of Theorem \ref{SardSmale} are
satisfied with $\xi =0$.
\vskip 1mm

By Theorem \ref{BPhivareg}, $\Phi$ is a $C^1$-mapping from ${\cal X}
\times {\cal Y}$ into ${\cal Z}$ and it is also a Fredholm operator
of index $0$ for any $(w,f) \in \Phi^{-1}(0)$. Thus Hypothesis 1 of
Theorem \ref{SardSmale} holds.

By Theorem \ref{BDPhisurj}, $D\Phi(w,f) : T_w{\cal X} \times
T_{f}{\cal Y} \rightarrow T_0{\cal Z}$ is surjective, for
$(w,f) \in \Phi^{-1}(0)$. Thus Hypothesis 2 of
Theorem \ref{SardSmale} also holds.

Taking into account Lemma \ref{5orbitn} and the
remark following this lemma, we can prove Hypothesis 3(b) by 
following the lines
of the proof of the corresponding property in \cite{BrunRaug} (see \cite[Step 3
of the proof of Proposition 4.12]{BrunRaug}).

Thus, the functional $\Phi$ satisfies all the assumptions of the 
Sard-Smale theorem. Due
to the relation between $\Phi$ and the transversality of connecting orbits
(Theorem \ref{BPhivareg}), this shows the genericity of ${\cal
G}_{m}^{\ell,d}$ in $\Vm(f_0)$ and concludes the proof of Theorem \ref{th2}.

\section{The non-wandering set}

To prove Theorem \ref{corointro}, it remains to show that generically
there does not exist non-wandering elements which are not critical
elements and that generically the number of  critical elements is
finite. We emphasize that the dynamics of \eqref{eq} may have
non-trivial non-wandering elements. Indeed, as shown in
\cite{SandFied}, every two-dimensional flow can be realized
in the dynamics of \eqref{eq}. Thus, one can for example create a
sequence of periodic orbits, which piles up on a homoclinic orbit.
This orbit is then a non-wandering orbit which is not critical.
However, non-trivial non-wandering orbits are generically precluded.

\begin{prop}\label{pnoerrant}
Assume that $f$ is a non-linearity such that the dynamics of
\eqref{eq} satisfy the following properties:\\
- there exists a compact global attractor for \eqref{eq}.\\
- All the equilibria and periodic orbits  of \eqref{eq} are hyperbolic.\\
- There is no homoclinic orbit and all the heteroclinic orbits are
transversal.\\
Then, the set of non-wandering elements consists in a finite number
of equilibrium points and periodic orbits.
\end{prop}
As usual in this article, the property stated in Proposition
\ref{pnoerrant} has its equivalent for two-dimensional dynamical
systems.  It mainly relies on Poincar\'e-Bendixson property, proved in
\cite{FMP} for \eqref{eq}.  Proposition \ref{pnoerrant} is the
key point to deduce the genericity of Morse-Smale property from the
genericity of Kupka-Smale property.  The genericity of Morse-Smale
property for dynamical systems of orientable surfaces shown in
\cite{Pei} also relies on a similar property, see \cite{PalisMelo}.

We enhance that, if $f$ is such that \eqref{eq} admits a compact
global attractor and that any equilibrium point and any periodic
orbit are hyperbolic, then
there is at most a finite number of equilibrium points.  However, as
we explained above, there could exist an infinite number of
hyperbolic periodic orbits:
think of a sequence of hyperbolic periodic orbits piling up to a
homoclinic orbit.  One can only ensure that there is a finite number
of hyperbolic periodic orbits with a period less than a given number.

\vskip 2mm

We begin the proof of Proposition \ref{pnoerrant} by several
lemmas.  We assume in the whole section that $f$ has been chosen so that the
hypotheses of Proposition \ref{pnoerrant} are satisfied.

Let $\Cm$ be a hyperbolic equilibrium point or periodic orbit of
\eqref{eq}. We recall that
there exists an open neighbourhood $B$ of $\Cm$ in $H^s(S^1)$ such
that each global solution $u(t)$ of \eqref{eq}, satisfying $u(t)\in
\overline B$ for all $t\leq 0$, belongs to the local unstable manifold
$W^u_{loc}(\Cm)$. We refer for example to \cite{HJR}, \cite{Ruelle}.

\begin{lemma}\label{lnoerrant2}
Let $\Cm$ be a hyperbolic equilibrium point or periodic orbit of
\eqref{eq} and let $B$ be
the neighbourhood of $\Cm$ as described above. Let
$\left(u_n(t)\right)_{n\in\Nm}$ be a sequence of solutions of
\eqref{eq} such that, for each $n\in\Nm$, there exist three times
$\sigma_n<t_n<\tau_n$ such that the following properties hold. For
all $t\in (\sigma_n,\tau_n)$, $u_n(t)\in B$, $u_n(\sigma_n)\in
\partial B$, $u_n(\tau_n)\in \partial B$ and,
$$d\left(u_n(t_n),\Cm\right):=\inf_{c\in \Cm}
\|u_n(t_n)-c\|_{H^s(S^1)} \xrightarrow[~ n \longrightarrow +\infty
~]{} 0 ~ . $$
Then, there exist an extraction $\varphi$ and a globally defined and
bounded solution $u_\infty(t)$ of \eqref{eq} such that
$u_\infty(t)\in W^u_{loc}(\Cm)$, $t \leq 0$, and
\begin{equation}
\label{convTT}
\forall~T>0,~~\sup_{t\in[-T,T]}\left\|u_{\varphi(n)}(\tau_{\varphi(n)}
+ t)- u_\infty(t)\right\|_{H^s(S^1)} \xrightarrow[~n\longrightarrow
+\infty~]{} 0~.
\end{equation}
\end{lemma}
\begin{demo}
First, we claim that $\tau_n-t_n\longrightarrow +\infty$ when
$n\rightarrow +\infty$. Indeed, if this is not true, since $\Cm$ is
compact, we can extract a subsequence $u_{\psi(n)}(t_{\psi(n)})$
converging to some $c\in\Cm$ and such that
$\tau_{\psi(n)}-t_{\psi(n)}$ converges to some $t\geq 0$. Then, by
continuity of the Cauchy problem related to \eqref{eq},
$u_{\psi(n)}(\tau_{\psi(n)})$ converges to a point of $\Cm$,
which contradicts the fact that $u_{\psi(n)}(\tau_{\psi(n)})\in
\partial B$.\\
We set $T=1$. Since \eqref{eq} admits a compact global attractor and
that $(\tau_n)$ converges to $+\infty$, the sequence $u_n(\tau_n-T)$
is precompact in $H^s(S^1)$ and there is an extraction $\varphi_1$
such that $u_{\varphi_1(n)}(\tau_{\varphi_1(n)}-T)$ converges to some
$u_\infty(-T) \in H^s(S^1)$. Let $u_\infty (t)= S(t+T)u_\infty(-T)$,
$t\geq -T$, be the solution of
\eqref{eq} associated to $u_\infty(-T)$. By continuity of the Cauchy
problem related to \eqref{eq},
$u_{\varphi_1(n)}(\tau_{\varphi_1(n)} +t)$ converges to $u_\infty(t)$
uniformly with respect to $t\in [-T,T]$. To achieve the proof of
uniform convergence of $u_n(t)$ to $u_\infty(t)$ on any compact set
of time, it is sufficient to repeat the argument for all $T\in\Nm$
and to use the diagonal extraction $\varphi(n)=\varphi_n\circ ...
\circ\varphi_1(n)$.\\
Finally, let us notice that $u_\infty(t)$ belongs to $W^u_{loc}(\Cm)$.
Indeed, since $u_n(t)\in B$ for all $t\in[t_n,\tau_n)$ and that
$\tau_n-t_n\longrightarrow +\infty$, $u_\infty(t)\in
\overline B$ for all $t\leq 0$. Due to the choice of $B$, this
implies that $u_\infty(t)\in W^u_{loc}(\Cm)$. \end{demo}

Let $M\in\Rm\cup\{+\infty\}$. We use the notation $\lcro 1,M +1 \rcro
=\{k\in\Nm, 1\leq k \leq M +1\}$. We say that a sequence of critical
elements $(\Cm_k)_{k\in \lcro 1,M+1 \rcro}$ is connected if for any
$k\in \lcro 1,M \rcro$, there exists a heteroclinic orbit $u_k(t)$
such that the $\alpha-$limit set of $u_k(t)$ is $\Cm_k$ and its
$\omega-$limit set is $\Cm_{k+1}$. We recall that a chain of
heteroclinic orbits denotes the sequence of heteroclinic orbits
corresponding to a connected sequence of critical elements $(\Cm_k)_{k\in
\lcro 1,p+1 \rcro}$ with $\Cm_{p+1}=\Cm_1$.

\begin{lemma}\label{lnoerrant3}
Assume that $f$ is as in Proposition \ref{pnoerrant}. Then,
there is no connected sequence of critical elements with infinite length.
As a consequence, there is no chain of heteroclinic orbits and, every
$\omega-$limit set and every non-empty $\alpha-$limit set of trajectories
of \eqref{eq} consist exactly of one critical element
\end{lemma}

\begin{demo} Let $M\in \Nm\cup \{+\infty\}$ and let $(\Cm_k)_{k\in
\lcro 1,M+1 \rcro}$ be a connected sequence of closed orbits with
heteroclinic connections $(u_k(t))_{k\in\lcro 1,M\rcro}$. We consider
the Morse indices $i(\Cm_k)$ of the closed orbits. We have several
cases:\\
- if $\Cm_k$ and $\Cm_{k+1}$ are both periodic orbits, then Theorem
\ref{CzaRocha} shows that $i(\Cm_k)>i(\Cm_{k+1})$.\\
- if $\Cm_k$ is an equilibrium point and if $\Cm_{k+1}$ is an
equilibrium point or a periodic orbit, then
dim($W^u(\Cm_k)$)=i($\Cm_k$) and
codim($W^s(\Cm_{k+1})$)=i($\Cm_{k+1}$). Thus, since the intersection
of $W^u(\Cm_k)$ and $W^s(\Cm_{k+1})$ is non-empty and transversal,
one must have $i(\Cm_k)>i(\Cm_{k+1})$.\\
- if $\Cm_k$ is a periodic orbit and $\Cm_{k+1}$ is an equilibrium,
then dim($W^u(\Cm_k)$)=i($\Cm_k$)$+1$ and
codim($W^s(\Cm_{k+1})$)=i($\Cm_{k+1}$). Therefore, $i(\Cm_k)\geq
i(\Cm_{k+1})$.\\
Hence, the Morse index of $\Cm_k$ is non-increasing and decreases
except if $u_k$ goes from a periodic orbit to an equilibrium point.
However, a sequence cannot consist only of connections from a
periodic orbit to an equilibrium and the Morse index must decrease at
least every two steps. Therefore, since the Morse indices are
non-negative, $M$ is bounded by $2i(\Cm_0)$.

The non-existence of connected sequence of critical elements of infinite
length precludes the existence of chains of heteroclinic orbits since
every chain $(\Cm_k)_{k\in\lcro 1,p+1\rcro}$ with $\Cm_{p+1}=\Cm_1$
can be extended to a periodic connected sequence of critical
elements and thus to a connected sequence of infinite length.

Let $u_0 \in H^s(S^1)$ be chosen such that its $\omega-$limit
set $\omega(u_0)$ is not a unique periodic orbit. Then, the
Poincar\'e-Bendixson
property stated in Theorem \ref{Poincare-B} shows that
$\omega(u_0)$ consists of equilibrium points and homoclinic or
heteroclinic orbits connecting them. We know that homoclinic orbits
are precluded. Moreover, since there is no connected sequence of equilibria
of infinite length, there exists an equilibrium point $e$ where no
connected sequence can be extended, that is, such that $W^u(e)\cap
\omega(u_0)=\{e\}$. Let $B_e$ be a small neighbourhood of $e$ in
$H^s(S^1)$ such that any solution $u(t)$ of \eqref{eq} satisfying
$u(t) \in \overline{B}_e$ for any $t \leq 0$, belongs to the unstable
manifold $W^{u}(e)$.  If $\omega(u_0) \ne \{e\}$, then one easily
constructs three sequences of times $(\sigma_n)$, $(t_n)$ and
$(\tau_n)$ going to $+\infty$ and satisfying the hypotheses of Lemma
\ref{lnoerrant2} with $\Cm =\{e\}$ and
$u_n(t)=S(t)u_0$. But then Lemma \ref{lnoerrant2} implies that there exists
a solution $u_{\infty}(t)$ of \eqref{eq} belonging to the unstable
manifold $W^u(e)$ and with $u(\tau_{\varphi(n)})$ converging
$u_{\infty}(0)$.  Therefore $u_{\infty}(0)\in\omega(u_0)\cap \partial
B_e$ and $W^u(e) \cap \omega(u_0) \ne \{e\}$, which leads to a
contradiction.  Therefore, $\omega (u_0) = \{e\}$.
\end{demo}

\noindent{\bf Proof of Proposition \ref{pnoerrant}:} Let $\tilde
u_0\in H^s(S^1)$ be a non-wandering element and let $\tilde u(t)= S(t)
\tilde u_0$.  Using the definition of a non-wandering element, we
easily construct a sequence of trajectories $u_n(t)$ such that
$u_n(0)$
converges to $\tilde u_0$ and, for any sequence $(t_n)$, there exists
a sequence $(t'_n)$ such that $t'_n>t_n$ and $u_n(t'_n)$ converges to
$\tilde u_0$. By Lemma \ref{lnoerrant3}, there exists a
hyperbolic critical element $\Cm_1$ such that $\omega(\tilde
u_0)=\{\Cm_1\}$. Let $B_1$ be a neighbourhood of $\Cm_1$ as in Lemma
\ref{lnoerrant2}. Assume that $\tilde u_0 \not\in \Cm_1$.
Replacing $B_1$ by a smaller neighbourhood if needed,
we may assume that $\tilde u_0 \not \in \overline
B_1$. There is a sequence of times $(t_n)$ and a point $c\in\Cm_1$
such that $\tilde u(t_n)\rightarrow c$. By continuity of the Cauchy
problem, we may assume without loss of generality that
$u_n(t_n)\rightarrow c$. As we can find a sequence of times $t'_n$
such that $t'_n>t_n$
and $u_n(t'_n)\rightarrow \tilde u_0\not\in \overline B_1$, there
exists a sequence of  times $\tau^1_n$ such that $u_n(t)\in B_1$ for $t\in
[t_n,\tau^1_n)$ and $u_n(\tau^1_n)\in \partial B_1$. By Lemma
\ref{lnoerrant2}, we may assume without loss of generality
that $u_n(\tau^1_n+t)$ converges to some $\tilde u_1(t)\in
W^u(\Cm_1)$. Now, we can repeat the arguments: there exist a critical
element $\Cm_2$ such that $\omega(\tilde u_1(t))=\{\Cm_2\}$ and a
sequence of times
$\tau^2_n$ such that, up to an extraction, $u_n(\tau^2_n+t)$
converges to some solution $\tilde u_2(t)\in W^u(\Cm_ 2)$ and so
on... Thus, we are constructing a connected sequence of critical
elements of infinite length, which is precluded by Lemma
\ref{lnoerrant3}. This means that $\tilde u_0$ belongs to
$\Cm_1$, that is, that our non-wandering element either is an equilibrium point
or belongs to a periodic orbit.

To finish the proof of Proposition \ref{pnoerrant}, it suffices to
show that the number of critical elements is finite. First, as
noticed earlier, due to the compactness of the global attractor and
the hyperbolicity of the equilibrium points and the periodic orbits,
the number of equilibrium points and periodic orbits of smallest period
bounded by a given number is finite. Thus we only need to show that
there does not exist an infinite sequence of periodic orbits
$\gamma_n(x,t)$ of smallest period $p_n$, where $p_n$ tends to infinity when
$n$ goes to infinity. If we had such a sequence, we would be able to
repeat the arguments of the first part of the proof and, by using
lemmas \ref{lnoerrant2} and \ref{lnoerrant3}, construct a
connected sequence of critical elements of infinite length, which
leads to a contradiction.
\hfill $\square$


\section*{Appendices}
\setcounter{section}{0}
\renewcommand{\thesection}{{\Alph{section}}}

\section{Fredholm operator and Sard-Smale Theorem}\label{AppendixA}

As we have explained in the introduction, an ingredient of the proof of the
generic non-existence of homoindexed orbits is the Sard-Smale
theorem. We
recall the precise statement of the version of the Sard-Smale
theorem, that we are applying in Section 5.

Let ${\cal X}$, ${\cal Z}$ be two differentiable Banach manifolds and 
$\Phi: {\cal X}
\rightarrow {\cal Z}$ be a $C^1$ map. A point $z$ is a {\sl regular
value} of $\Phi$ if, for any $x \in \Phi^{-1}(z)$ the derivative
$D\Phi(x)$ is surjective and its kernel splits, i.e. has a closed
complement in $T_x {\cal X}$ (sometimes this property is denoted by
$\Phi \, \pitchfork  z$). A point $z \in {\cal Z}$ which is
not regular is
called a {\sl critical value} of $\Phi$.  A subset in a topological
space is {\sl
generic} or {\sl residual} if it is a countable intersection of
open dense sets.
\vskip 1mm

We recall that a continuous linear map $L: X \rightarrow Z$ between two
Banach spaces $X$ and $Z$, is a Fredholm map if its range $R(L)$ is closed
and if both $\dim \, \ker (L)$ and $\codim \, R(L)$ are finite. The {\sl
index} $\ind (L)$ is the integer $\ind (L) = \dim \, \ker (L) -
\codim \, R(L)$.
\vskip 1mm

The version of the Sard-Smale theorem given here has been proved in
\cite{He05} (for weaker versions, we also refer to \cite{Quinn} and
\cite{SautTem}). The next theorem has been widely used in the
genericity proofs in \cite{BrunPola} and \cite{BrunRaug}.

\begin{theorem} {\bf Sard-Smale Theorem}
\label{SardSmale}

Let ${\cal X}, {\cal
Y}, {\cal Z}$ be three smooth
Banach manifolds, $\Phi: {\cal X} \times {\cal Y} \rightarrow
{\cal Z}$ be a mapping of class $C^r$, $r \geq 1$ and $\xi$ an element
of ${\cal Z}$.
Assume that the following hypotheses hold:
\begin{enumerate}
\item  for each $(x,y) \in \Phi^{-1}(\xi)$, $D_x\Phi(x,y)$ is a
Fredholm operator of index strictly less than $r$;

\item for each $(x,y) \in \Phi^{-1}(\xi)$, $D\Phi(x,y): T_x
{\cal X}  \times T_y{\cal Y} \rightarrow T_{\xi}{\cal Z}$ is surjective;

\item  one of the following properties is satisfied: \HB
(a) ${\cal X}$ and ${\cal Y}$ are separable metric spaces; \HB
(b) the map $(x,y)  \in \Phi^{-1}(\xi) \mapsto y \in {\cal
Y}$ is $\sigma$-proper, that is, there is a countable system of
subsets $V_n \subset \Phi^{-1}(\xi)$ such that $\cup_n V_n=
\Phi^{-1}(\xi)$ and for each $n$ the map $(x,y) \in V_n \mapsto y \in
{\cal Y}$ is proper (i.e. any sequence $(x_k,y_k) \in V_n$ such that
$y_k$ is convergent in ${\cal Y}$ has a convergent subsequence in
$V_n$).
\end{enumerate}
Then, the set $\{y \in {\cal Y}\, | \, \xi \hbox{ is a regular value
of } \Phi(.,y) \}$ contains a countable intersection of open dense sets
(and hence dense) in ${\cal Y}$.
\end{theorem}

\section{Exponential dichotomies, shifted exponential
dichotomies and applications to transversality}\label{AppendixB}

In the proof of the generic non-existence of homoindexed connections
between equilibria, we will use a functional characterization of the
transversality property. A main tool in this proof is the notion of
dichotomy.
Also, when studying the asymptotics
of the solutions of linearized equations along connecting orbits
connecting a hyperbolic periodic orbit to another critical element,
we will need to consider iterates of maps and thus, in particular,
discrete shifted dichotomies.
We begin this appendix by recalling the definition and the basic
properties of the exponential dichotomies and shifted exponential
dichotomies. Then, we give
some applications to the scalar parabolic equation on $S^1$.\HB
The results, that we recall here, are all contained in
\cite{He81}, \cite{Li86}, \cite{HaLi86}, \cite{He94}, \cite{ChChHa},
\cite{BrunPola}, \cite{BrunRaug}, \cite{Pal84} and \cite{Pal88}.

\subsection{Generalities}

Let $X$ be a Banach space and
$J$ be an interval in $\Zr$.
    Let $\{T_n \in L(X,X)\, | \,n \in J \}$ be a
family of continuous maps from $X$ into $X$. We define the family of
evolution operators
$$
T(m,m) \, = \, I~,\quad T(n,m) \, = \, T_{n -1} \circ
\ldots \circ T_m~, \forall
\, n \geq m  ~\hbox{ in }J~,
$$
where $I$ is the identity in $X$.

\begin{defi} \label{Adefdich} We say that the family of linear
operators $T_n$, $n \in J$, or the family of evolution operators
    $\{T(n,m) \, | \, n \geq m \hbox{ in }J\}$, admits an {\rm
exponential dichotomy} (or {\rm discrete dichotomy}) on the interval
$J$ with {\rm exponent}
$\beta>0$ (or {\rm constant} $e^{-\beta}$), {\rm bound} $M >0$ and
projections $P(n)$ if there
is a family of continuous projections $P(n)$, $n\in J$, such that
the following properties hold for any $n$ in $J$:
\begin{description}
\item[(i)] $T(n,m)P(m)=P(n)T(n,m)$ for $n \geq m$ in $J$,

\item[(ii)]  the restriction $T(n,m)|_{ R(P(m))}$ is an isomorphism of
$R(P(m))$ onto $R(P(n))$, for $n \geq m$ in $J$, and $T(m,n)$ is
defined as the inverse from $R(P(n))$ onto $R(P(m))$,

\item[(iii)] $\|T(n,m)(I-P(m))\|_{L(X)} \leq Me^{-\beta(n-m)}$ for $n \geq
m$ in $J$,

\item[(iv)] $\|T(n,m)P(m)\|_{L(X)} \leq Me^{-\beta(m-n)}$ for $n \leq
m$.
\end{description}
\end{defi}

If $\dim \, R(P(n))=k$ is finite for some $n \in J$, the equality holds
for all $n \in J$ by (ii) and we say that the dichotomy has finite {\sl rank}
$k$.
\vskip 1mm

\begin{rems} \label{AexPara} \HB
1) We could also have defined (continuous) exponential dichotomies on 
an interval
$J\subset\Rm$ (see \cite{He81} for instance). However,
here we restrict our discussion to discrete dichotomies on
time intervals $J\subset \Zr$ for two reasons. First, as it
was already pointed out by Henry in \cite{He94}, the theory of
dichotomies is much simpler with discrete time and there is little
loss in restricting to this case. Moreover, in our applications to
the asymptotics near periodic orbits, we really need to work with
maps and not with continuous evolution operators.
Finally, we remark
that constructing discrete dichotomies from continuous ones or
conversely is an easy task (see Theorem 1.3 of
    \cite{He94} for example). \HB
2) Let $A$ be a sectorial operator on a Banach space $Y$. For any $J \subset
\Zr$ and $ n \geq m \hbox{ in }J $, we define
$T_n = e^A$ (independent of $n$) and
$T(n,m)= e^{A(n-m)}$ on the Banach space $X=Y^{\alpha}$, $\alpha \in [0,1]$.
Then $\{T(n,m) \, | \, n\geq m \hbox{ in }J\}$ is a family
of evolution operators on $X$.  If the spectrum $\sigma(A)$ satisfies
$\sigma(A) \cap \{ \mu \, | \, \mathrm{Re}\, \mu=0 \} = \emptyset$,
then, for any $t_0 >0$,  we can define the projection $P$ by
\begin{equation}
\label{Adichconst}
P \, = \,I - {\frac 1 {2 i\pi}} \int_{|z| =1} (z I -
e^{At_0})^{-1}\,{\rm d}z  ~.
\end{equation}
And $T(n,m)$ has an exponential dichotomy with projection $P$.
If the spectrum
$\sigma(A)$ satisfies $\sigma(A) \cap \{ \mu \, | \, -\beta \leq
\mathrm{Re}\, \mu \leq \beta\} = \emptyset$ for
some $\beta>0$, then there exists a
positive constant $M$ such that $T(n,m)$ has an exponential dichotomy
with projection $P$, exponent $\beta$ and constant $M$.  If the
essential spectrum of $e^{At}$ is strictly inside the unit circle, the
    dichotomy has finite rank.  This is the case for the linear
    parabolic equation.
\end{rems}

    Let again $\{T_n \in L(X,X)\, | \,n \in J \}$ be a
family of continuous maps from $X$ into $X$. We define the family of
operators on
$X^*$, given by
$$
T^*(m,n) \, = \, (T(n,m))^*~.
$$

\begin{defi} \label{Adefstar} We say that the family of maps
$T_{n}^*$, $n \in J$, or the family of evolution
operators $ \{T^*(m,n) \, | \, n \geq m \hbox{
in }J\}$, admits a {\rm reverse
exponential dichotomy} on the interval $J$ with {\rm exponent}
$\beta>0$, {\rm bound} $M >0$ and projections $P^*(t)$ if there is
a family of continuous projections $P^*(n)$, $n \in J$, such that the
following properties hold, for any $n$ in $J$:
\begin{description}
\item[(i)] $T^*(m,n)P^*(n)=P^*(m)T^*(m,n)$ for $n\geq m$ in $J$,

\item[(ii)] the restriction $T^*(m,n)|_{ R(P^*(n))}$ is an isomorphism of
$R(P^*(n))$ onto $R(P^*(m))$, for $n \geq m$ in $J$, and $T^*(n,m)$ is
defined as the inverse from $R(P^*(m))$ onto $R(P^*(n))$,

\item[(iii)] $\|T^*(m,n)(I-P^*(n))\|_{X} \leq Me^{-\beta(n-m)}$
for $n\geq m$ in $J$,

\item[(iv)] $\|T^*(m,n)P^*(n)\|_{X} \leq Me^{-\beta(m-n)}$ for
$n \leq m$.
\end{description}
\end{defi}

The following natural property is proved in \cite{BrunPola} for
example.

\begin{lemma} \label{Alemstar} If the family of evolution operators $T(n,m)$,
$n,m\in J$ on the Banach space $X$ admits an exponential dichotomy on
the interval $J$, with projections $P(n)$, exponent $\beta$ and
bound $M$, then $T^*(m,n) = (T(n,m))^*$ admits reverse
exponential dichotomy on $J$ with the same exponent and bound and
with the projections $P^*(n)= (P(n))^*$.
\end{lemma}

In our applications to the asymptotics near a periodic orbit, we
cannot directly use the concept of exponential dichotomy since $1$
always belongs to the spectrum of the period map $\Pi(p,0)$
associated to a periodic orbit $\gamma (t)$ of \eqref{eq} of minimal
period $p$. For this reason, we also recall the notion of {\sl
shifted exponential dichotomy}, which is a generalization of the
notion of {\sl exponential dichotomy}.
Calling these dichotomies {\sl shifted}, we follow the terminology of
\cite{HaLi86}; alternatively it is sometimes called {\sl
pseudodichotomy}. If $\lambda_1 < 1 < \lambda_2$, it reduces to the
usual exponential dichotomy. \HB

\begin{defi} \label{Cshifted} We say that the family of linear
operators $T_n$,
$ n\in J$, or the family of linear evolution operators $\{T(n,m) \, | \,
n \geq m \hbox{ in }J\}$, admits a {\rm shifted exponential dichotomy}
on the interval $J$ with {\rm gap} $[\lambda_1,\lambda_2]$, {\rm
bound} $K >0$ and projections $P(n)$, $Q(n)= I -P(n)$ if there is a
family of continuous projections $P(n)$, $n \in J$, such that the
following properties hold for any $n$ in $J$:
\begin{description}
\item[i)]   $T(n,m)P(m)=P(n)T(n,m)$ for $n\geq m$ in $J$,

\item[ii)]   the restriction $T(n,m)|_{ R(P(m))}$ is an isomorphism of
$R(P(m))$ onto $R(P(n))$, for $n \geq m$ in $J$, and
on $R(P(n))$, $T(m,n)$ is
defined as the inverse from $R(P(n))$ onto $R(P(m))$,

\item[iii)]  $\|T(n,m)(I-P(m))\|_{X} \leq K \lambda_1^{n-m}$ for $n\geq
m$ in $J$,

\item[iv)]    $\|T(n,m)P(m)\|_{X} \leq K\lambda_2^{n-m}$ for $n \leq
m$.
\end{description}
We also say that the family of operators $T_n$, $n \in J$, or the family
of evolution operators $\{\Phi(n,m) \, | \, n \leq m \hbox{ in }J\}$,
admits a {\rm reverse shifted exponential dichotomy} on the interval $J$ with
{\rm gap} $[\lambda_1,\lambda_2]$, {\rm bound} $K >0$ and
projections $P(n)$, $Q(n)= I -P(n)$ if the above properties hold
with $n\geq m$ in $J$ (resp. $n\leq m$ in $J$) replaced by
$n\leq m$ in $J$ (resp. $n \geq m$ in $J$).
\end{defi}

\vskip 1mm
Before describing the applications of these abstract notions to our problem,
we recall two properties, which are very useful in
proving that linearized equations along connecting orbits
of \eqref{eq} admit exponential dichotomies.
    The next roughness property is given in \cite[Corollary 1.9]{He94}
and is a consequence of \cite[Theorem 1.5 and Lemma
1.6]{He94} (see also \cite[Theorem 7.6.7]{He81} as well as
Proposition \ref{C2} below). The extension of the result below to
trichotomies is stated in \cite{HaLi86}.

\begin{theorem} \label{Adichpertur} {\bf Perturbation of exponential
dichotomies}

Let $n_0 >0$ (resp.  $n_0 < 0$) be
a
given integer and let $T(n+1,n)$, $n \in \Zr^{+}$, $n \geq
n_0$, (respectively $n \in \Zr^{-}$, $n \leq n_0$), be a
discrete family of evolution operators on a Banach space $X$
admitting a discrete dichotomy on $[n_0, +\infty)$ (respectively
on $(-\infty, n_0]$), with exponent $\beta$, constant $M$ and
projections $P^{T}(n)$.  Let $M_1 > M$, $0 < \beta_1 < \beta$ and
$0 < \varepsilon \leq ({\frac 1 M} - {\frac 1 {M_1}}) \frac{e^{-\beta_1} -
e^{-\beta}}{1 + e^{-(\beta + \beta_1)}}$.
If $S(n+1,n)$, $n \in {\Zr}^{+}$, $n \geq n_0$
(respectively $n \in \Zr^{-}$, $n \leq n_0$), is a discrete
family of evolution operators on $X$ with $\|S(n+1,n) -
T(n+1,n)\|_{L(X,X)} \leq \varepsilon$, for all $n \geq n_0$ in
$\Zr^{+}$ (respectively for all $n \leq n_0$ in $\Zr^{-}$),
then $S(n+1,n)$ admits a discrete dichotomy on $[n_0,
+\infty)$ (respectively on $(-\infty, n_0]$) with exponent
$\beta_1$, constant $M_1$ and projections $P^{S}(n)$. Moreover,
the projections $P^{S}(n)$ satisfy $\sup_{n} \|P^{T}(n) -
P^{S}(n)\|_{L(X,X)} =O( \sup_{n}\|T(n+1,n) - S(n+1,n)\|_{L(X,X)})$
as $\sup_{n} \|T(n+1,n) - S(n+1,n)\|_{L(X,X)}$ tends to $0$.
Furthermore, there exists $\varepsilon_0 >0$ such that, for $0<
\varepsilon \leq \varepsilon_0$, if $T(n+1,n)$ has a dichotomy of
finite rank $m$, then the dichotomy of $S(n+1,n)$ is also of
finite rank $m$.
\end{theorem}

The next result, which is proved in \cite[Theorem 1.14]{He94},
allows to extend dichotomies from smaller to larger ``time intervals".
The continuous version of it is proved in \cite{Li86} and the
extension to trichotomies is stated in \cite{HaLi86}.

\begin{theorem} \label{Adichprolon} {\bf Extension of exponential dichotomies}

   Let $T(n+1,n)$, $n \in \Zr^{-}$, $n
< n_1$, be a discrete family of evolution operators on a Banach
space $X$, and suppose that, for $n < n_0$, with $n_0 <n_1$,
$T(n+1,n)$ admits a discrete dichotomy with finite rank $m$, exponent
$\beta$, constant $M$ and projections $P(n)$, $n \leq n_0$.  Assume
moreover that $T(n_1,n_0)_{\big|R(P(n_0))}$ is injective. Then,
$T(n+1,n)$, for $n < n_1$, admits a discrete dichotomy with the
same rank $m$, same exponent and projections $\tilde{P}(n)$, $n \leq
n_1$,
such that $\|P(n) - \tilde{P}(n)\|_{L(X)}\rightarrow 0$ exponentially as
$n$ goes to $-\infty$.  The constant $M$ has to be replaced by a
larger one.\HB
Let $T(n+1,n)$, $n \in \Zr^{+}$, $n \geq n_0$, be a discrete
family of evolution operators on a Banach space $X$. Suppose that,
for $n \geq n_1$, with $n_0 <n_1$, $T(n+1,n)$ admits a discrete
dichotomy with finite rank $m$, exponent $\beta$, constant $M$ and
projections $P(n)$, $n \geq n_1$.  Assume moreover that the adjoint
operator $T^\ast(n_0,n_1)_{\big|R(P^\ast(n_1))}$ is injective, then
$T(n+1,n)$, for $n \geq n_0$, admits a discrete dichotomy with the
same rank $m$, same exponent $\beta$ and projections $\tilde{P}(n)$,
$n \geq n_0$,
such that $\|P(n) - \tilde{P}(n)\|_{L(X)}\rightarrow 0$ exponentially as
$n$ goes to $+\infty$.  The constant $M$ has to be replaced by a
larger one. \HB
In both cases, the convergence of $\|P(n) - \tilde{P}(n)\|_{L(X)}$ is
of order $O(e^{-2\beta |n|})$.
\end{theorem}

\subsection{Dichotomies and Fredholm property}\label{AppendixC}

In Section 5, in proving the generic non-existence of
homoindexed connecting orbits, we use a functional characterization
of the transversality, which is based on dichotomies and a Fredholm
property. In this appendix, we quickly recall the tools and basic
facts, which lead to this functional characterization.

Let $X$ be a Banach space. We introduce the following spaces
$$
\Zm \, = \, \ell^{\infty}({\bf \Zr}, X)~,
\hbox{ (respectively }
{\cal Z}^{\pm} \, = \, \ell^{\infty}({\bf \Zr}^{\pm}, X))~.
$$
Given  a family $T(n,m)$ of evolution operators on $X$ for $n,m \in {\bf Z}$,
we define the mapping ${\cal L}$ from $X^{{\bf Z}}$ into $X^{{\bf Z}}$
(respectively the mapping ${\cal L}^{\pm}$ from $X^{{\bf Z}^{\pm}}$
into $X^{{\bf Z}^{\pm}}$) by
     \begin{equation}
\label{AcalL}
({\cal L}Y)(n) \, = \, Y(n+1) -T(n+1,n)Y(n)~, \quad \forall n \in {\bf
Z}~,
\end{equation}
(respectively $({\cal L}^{\pm}Y)(n) = Y(n+1) -T(n+1,n)Y(n)$,
$\forall n \in {\bf Z}^{\pm}$). \HB
We say that $Y=\{Y(n)\}_{n \in {\bf Z}}$ belongs to the domain
$D({\cal L})$ if  $\sup_{n \in {\bf Z}} \|Y(n+1)
-T(n+1,n)Y(n)\|_{X} < \infty$ (likewise, we define $D({\cal
L}^{\pm})$). This allows to define
the operator ${\cal L}: D({\cal L}) \subset {\cal Z}
\rightarrow {\cal Z}$ by \eqref{AcalL} (likewise, we may define the
operator ${\cal L}^{\pm}: D({\cal L}^{\pm}) \subset {\cal Z}^{\pm}
\rightarrow {\cal Z}^{\pm}$).

\vskip 2mm
In \cite [Theorem 7.6.5]{He81}, Henry has given the
following characterization of the existence of a discrete
dichotomy for $T(n+1,n)$
     (see also \cite{He94}; for a finite-dimensional version, see
\cite{Pal84}, \cite{Pal88}). The family of evolution operators $T(n+1,n)$ has a
discrete dichotomy if and only if, for every bounded sequence $F \in {\cal
Z}$, there is a unique bounded sequence $Y \in {\cal Z}$ with $({\cal
L}Y)(n) := Y(n+1) - T(n+1,n)Y(n)=F(n)$, for any $n \in {\bf Z}$.
Moreover, the unique bounded solution is given by
     \begin{equation}
     \label{AGreen}
Y(n) \, = \, \sum_{k=-\infty}^{+\infty} \mathcal{G}(n, k+1) F(k)~,
\end{equation}
where $\mathcal{G}(n,m)=T(n,m)(I - P(m))$ for $n \geq m$, $\mathcal{G}(n,m) =
-T(n,m)P(m)$ for $n <m$, is called the Green function.
\par
Henry has also proved in  Theorem 1.13 of \cite{He94} that
any discrete family of evolution operators $T(n+1,n)$
admits a discrete dichotomy on ${\bf Z}$ if and only if the
restrictions to both ${\bf Z}^{+}$ and ${\bf Z}^{-}$ have dichotomies
and also $X=S_0 \oplus U_0$ where
\begin{equation}
\label{AU0S0}
\begin{split}
         U_0 \,&=\, \{x_0 \, | \, \exists \{x_n\}_{n \leq 0} \in {\cal
Z}^{-} \hbox{
         with }x_{n+1}=T(n+1,n)x_n\, ,\, n<0 \} \cr
         S_0 \,&=\, \{x_0 \, | \, \exists \{x_n\}_{n \geq 0} \in {\cal
Z}^{+} \hbox{
         with }x_{n+1}=T(n+1,n)x_n\, ,\, n \geq 0 \}~.
\end{split}
\end{equation}
When the dichotomies in ${\bf Z}^{+}$ and ${\bf Z}^{-}$ have finite
rank, the equality $X=S_0 \oplus U_0$ means that they have the same rank.

The previous existence result of a dichotomy on ${\bf Z}$ is actually
a particular case of the following more general result, which is also
proved in Theorem 1.15 of \cite{He94} (see also \cite{S78}, \cite{Pal84},
\cite{HaLi86} and  \cite[Theorem 4.a.4]{BrunPola} in the case of
ordinary differential, functional differential and parabolic equations).

We recall that $\langle \cdot, \cdot \rangle_{X,X^*}$ denotes the
duality pairing between $X$ and $X^*$.

\begin{theorem} \label{AFredholm} {\bf Fredholm alternative}

Let $T(n+1,n)$ be a discrete family
of evolution
operators on a Banach space $X$, admitting discrete dichotomies of
finite rank on both $\Zr^{+}$ and
$\Zr^{-}$, with projections $P^+(n)$ and $P^-(n)$. Then the
operator ${\cal L}$, defined by \eqref{AcalL}, belongs to $L(D({\cal
L}),{\cal Z})$ and is a Fredholm operator with index $\ind ({\cal
L})$ given by
\begin{equation}
\label{AFredh0}
      \ind ({\cal L}) \, = \, \dim (R(P^{-}(0))) -
\dim(R(P^{+}(0)))~.
\end{equation}
The codimension $\codim R({\cal L})$ of $R({\cal L})$ is given
by $\codim R({\cal L}) = \dim [R(I -P^{-\ast}(0)) \cap
R(P^{+\ast}(0))]$. \HB
A sequence $F \in {\cal Z}$ belongs to $R({\cal L})$ if and only if
\begin{equation}
\label{AFredh1}
         \sum_{n = -\infty}^{+\infty} \langle
         F(n), \Psi(n+1)\rangle_{X,X^*}\,=\,
         0~,
     \end{equation}
for every sequence $\Psi(n)=T^{\ast}(n,0) \Psi_0$, $\Psi_0
\in X^{\ast}$, which is bounded.
\end{theorem}

The proof of Theorem \ref{AFredholm} uses the following two auxiliary
lemmas. First, recall that, for any operator $Q \in L(X)$, one has
\begin{equation}
\label{kerQ*1}
\ker (Q^\ast) \, = \, (R(Q))^{\perp}~,
\end{equation}
where, for any subspace $X_0 \in X$, $X_0^{\perp}= \{\psi \in
X^\ast \, | \, \langle x, \psi  \rangle =0\, , \forall x \in
X_0\}$. \HB
If $Q \in L(X)$ is a projection, we have, in addition,
\begin{equation}
\label{kerQ*2}
R(I-Q^\ast) \, = \, \ker (Q^\ast) \, = \, (R(Q))^{\perp}~.
\end{equation}

\begin{lemma} \label{Tstarpsi0} Let $T(n,m)$ be an evolution operator
admitting discrete dichotomies of finite rank on both
${\bf \Zr}^{+}$ and ${\bf \Zr}^{-}$.
Then, any element $\Psi_0 \in X^\ast$ belongs to $(R(P^{-}(0)))^{\perp} \cap
(R(I-P^{+}(0)))^{\perp}$ if and only if the sequence
$$
\Psi (m) \, = \, T^\ast(m,0) \Psi_0~, \quad m \in \Zr~,
$$
(which is defined for all $m$ due to the property (ii) of the reverse
dichotomy) is bounded (that is belongs to $\ell^{\infty}(\Zr,X^\ast)$). In
this case, $\Psi (m)$ belongs to $R(I -P^{-\ast}(m))$ for $m \leq
0$ and to $R(P^{+\ast}(m))$ for $m \geq 0$.
\end{lemma}

The next lemma emphasizes the formulas given in \eqref{AGreen}.

\begin{lemma} \label{Greenpm} We assume that the hypotheses of
Theorem \ref{AFredholm} hold. Then, \HB
(i) if $F \in {\cal Z}^{-}$, there exists $Y \in {\cal Z}^{-}$ such
that $F= {\cal L}^{-}Y$ if and only if, for any $n \in \Zr^{-}$,
\begin{equation}
\label{AGreen-}
\begin{split}
Y(n)  \, = \,
T(n,0) P^{-}(0)Y(0) - \sum_{k=n}^{-1}&T(n,k+1) P^{-}(k+1) F(k) \cr
+ &\sum_{k=-\infty}^{n-1} T(n,k+1)(I - P^{-}(k+1)) F(k)~;
\end{split}
\end{equation}

(ii) similarly, if $F \in {\cal Z}^{+}$, there exists $Y \in {\cal
Z}^{+}$ such that $F= {\cal L}^{+}Y$ if and only if, for any $n \in {\bf
Z}^{+}$,
\begin{equation}
\label{AGreen+}
\begin{split}
Y(n) = T(n,0)(I-P^{+}(0)) Y(0) + \sum_{k=0}^{n-1} &T(n,k+1)(I -
P^{+}(k+1)) F(k) \cr
- &\sum_{k=n}^{+\infty} T(n,k+1) P^{+}(k+1) F(k)~.
\end{split}
\end{equation}
\end{lemma}

We remark that these ``variation of constants formulas" generalize the
formula \eqref{AGreen}. They have already been given in \cite{Pal88}
under this discrete form in the
finite dimensional context (see  \cite[Formula (13) of Lemma
2.7]{Pal88}) and they are contained in Theorem 1.15 of
\cite{He94}. In the
continuous case for parabolic equations, they are well-known and can be
found in \cite{He81} and in \cite{BrunPola}.

\subsection{Application to the parabolic equation on $S^1$}

In this section, we apply the previous abstract results to the homoclinic and
heteroclinic orbits between equilibrium points of the scalar 
parabolic equation \eqref{eq} and
    we give some equivalent formulations of transversality.

We assume in this section that $\tilde{u}(t) \in {\cal C}^0_b(\Rm,
H^s(S^1))$ is
a bounded trajectory of $S(t)=S_f(t)$ satisfying $\lim_{t
\rightarrow \pm}\tilde{u}(t)=e^{\pm}$, where $e^{\pm}$ are hyperbolic
equilibria of finite Morse index $i(e^{\pm})$. We recall that ${\tilde u}$
belongs to $C^0_b(\Rm, H^2(S^1))
\cap C^{\theta}(\Rm, H^s(S^1))$, where $0 < \theta \leq 1$.
We consider the linearized equation along $\tilde{u}$, that
is, the linear equation for $t \geq s$,
\begin{equation}
\label{Tutilde}
    v_t(t) \,=\, v_{xx}(t)+
      D_uf(x,\tilde{u},\tilde{u}_x) v(t)+ D_{u_x}f(x,\tilde{u},
      \tilde{u}_x)v_x(t)
     \equiv C_{{\tilde u}}(t) v(t)~, ~ t>\sigma~, ~ v(\sigma ) \,=\, v_0~.
\end{equation}
We recall that, for any $v_0
\in L^2(S^1)$, for any $\sigma \in \Rm$, there exists a unique classical
solution $v(t) \in {\cal C}^0 ([\sigma,+ \infty),L^2(S^1)) \cap {\cal
C}^0 ((\sigma,+ \infty), H^s(S^1))$ of
\eqref{Tutilde} such that $v(\sigma)=v_0$.
We set $T(t,\sigma)v_0=T_{\tilde{u}}(t,\sigma)v_0=v(t)$.

We next introduce the adjoint linearized equation to \eqref{Tutilde},
that is, the linear equation for $\sigma \leq t$,
\begin{equation}
\label{Tutilde*}
       \partial_t \psi (\sigma)  \,=\, -C_{{\tilde u}}^{*}(\sigma)
       \psi(\sigma)~, \quad \sigma \leq t~, \quad \psi(t) \,=\, \psi_0~.
\end{equation}
Since $\tilde{u}$ belongs to $C^{\theta}(\Rm, H^s(S^1))$ for any
$\theta \leq 1$ and thus that $(C_{{\tilde u}}(t)- \partial_{xx})^*$ is locally
H\"{o}lder continuous with
exponent $\gamma >s/2$ as a mapping from $\Rm$ into $L(H^s(S^1),
L^2(S^1))$, \eqref{Tutilde*}
has a unique classical solution $\psi(\sigma) = \psi(\sigma,t;
\psi_0)$
on $(-\infty, t)$, for any $\psi_0 \in L^2(S^1)$ (see \cite[Theorem
7.3.1]{He81} for example).
We denote this solution $\psi(s) := \psi(s,t; \psi_0)$.

With $T_{\tilde{u}}(t,s)$, we associate the adjoint evolution operator
on $L^2(S^1)$, given by
\begin{equation}
\label{ATadjoint}
T_{\tilde{u}}^*(\sigma,t) \, = \, (T_{\tilde{u}}(t,\sigma))^*~, \quad t
\geq \sigma~.
\end{equation}
It is well-known (see \cite[Theorem 7.3.1]{He81}) that,
for any $\psi_0 \in L^2(S^1)$,
\begin{equation}
\label{AT*psi}
T_{\tilde{u}}^*(\sigma,t_0) \psi_0\, = \, \psi(\sigma,t_0; \psi_0)~, \quad
\forall \sigma \leq t_0~.
\end{equation}
We also remark that
the adjoint operator $(T_{\tilde{u}}(t,\sigma))^*$ is injective and
that its range is dense in $H^s(S^1)$.

\vskip 2mm
{}From now on, we discretize the evolution operators. We fix a time step
$\tau>0$ and consider the discretizations $S(n\tau)$ and
$T_{\tilde{u}}(n\tau,m\tau)$, with $n,m\in\Zr$. The hyperbolic equilibria
$e^{\pm}$ of $S(t)$, their stable and unstable sets coincide with those of
the discretization $S(n\tau)$; the
discretization of the trajectory $\tilde{u}(t)$ connecting $e^-$ to $e^+$
is a heteroclinic or homoclinic orbit connecting these equilibria for the
discretized semi-flow. Let $\beta^{\pm} >0$ be chosen such that
$$
\sigma(e^{L_{e^{\pm}}}) \cap \{z \, | \, e^{-\beta^{\pm} } \leq |z| \leq
      e^{\beta^{\pm}}\} = \emptyset~,
      $$
      where the linearized operator $L_{e^{\pm}}$ has been defined in
      \eqref{MapLe}.
As explained in the remarks of Section B.1,
$e^{L_{e^{\pm}}\tau}$ admits an exponential dichotomy with
projection $P^{\pm}$
(see \eqref{Adichconst}), exponent $\beta^{\pm}$ and constant $M$ in
$H^{2\alpha}(S^1)$ for any $\alpha \in [0,1)$.
Thus, we will be able to deduce from Theorems \ref{Adichpertur} and
\ref{Adichprolon}
that $T_{\tilde{u}}(n\tau,(n-1)\tau)$ admits
exponential dichotomies on $\Zr^-$ and on $\Zr^+$ of respective index $i(e^-)$
and $i(e^+)$ in $H^{2\alpha}(S^1)$, for any $\alpha \in [0, 1)$. We
will only give a sketch of the proof. For
a more detailed proof in the case of ordinary differential
equations (resp. functional differential equations, resp. parabolic
equations, resp. in the case damped wave equations), we refer the
reader to \cite{Pal84}, \cite{Pal88} (resp. to \cite{Li86},
\cite{BrunPola} and \cite{BrunRaug}).

\begin{theorem} \label{ATpmdicho} For any $\beta_1^{\pm}\in (0,\beta^{\pm})$,
the discrete family of evolution operators $T(n,m)=T_{\tilde{u}}(n\tau,m\tau)$
admits exponential dichotomies on $\Zr^\pm$ in $L^2(S^1)$ (resp.
$H^s(S^1)$) of finite rank equal to the index
$i(e^{\pm})$ of $e^{\pm}$, with exponent $\beta_1^{\pm}$, constant
$M^{\pm}$ and
projections ${\tilde P}_{\tilde{u}}^{\pm}(n)$ (resp.
$P_{\tilde{u}}^{\pm}(n)$), satisfying
\begin{equation}
\label{APnPpm}
\lim_{n \rightarrow \pm \infty} \|{\tilde P}_{\tilde{u}}^{\pm}(n) -
P^{\pm}\|_{L(L^2,L^2)}\,=\,0~, \hbox{ (resp. }
\lim_{n \rightarrow \pm \infty} \|P_{\tilde{u}}^{\pm}(n) -
P^{\pm}\|_{L(H^s,H^s)}\,=\,0 \hbox{ )}~.
\end{equation}
Moreover, $P_{\tilde{u}}^{\pm}$ is the restriction of ${\tilde
P}_{\tilde{u}}^{\pm}$ to $H^s(S^1)$, that is,
\begin{equation}
\label{PnL2}
{\tilde P}_{\tilde{u}}^{\pm}(n)_{|_{H^s(S^1)}} =
P_{\tilde{u}}^{\pm}(n)~.
\end{equation}
\end{theorem}

\begin{demo} In the proof of Corollary \ref{CappD0} below, (see
\eqref{Dcoaux3}), we show
that, for $n \in \Zr^{\pm}$, with $| n |$ large enough, we have the
exponential asymptotic convergence
   \begin{equation}
\label{Basexp}
\|T_{\tilde{u}}((n+1)\tau,n\tau) - e^{L_{e^{\pm}}}\|_{L(H^s,H^s)}
\leq Ce^{- C |n|}~,
\end{equation}
   where $C$ is a positive constant. Thus we may apply Theorem
\ref{Adichpertur}, which implies that there exists $n_0  \in \Zr^{+}$
such that
$T_{\tilde{u}}((n+1)\tau,n \tau)$ admits an exponential dichotomy on
   $\Zr^{+}$ (resp. $\Zr^{-}$) in $H^s(S^1)$, for $n \geq n_0$ (resp.
$n \leq -n_0$)
   of finite rank $i(e^{+})$ (resp. $i(e^{-})$) and projections
   $P^{+}_{T}$ (resp. $P^{-}_{T}$) satisfying the properties of Theorem
\ref{Adichpertur}.
   Applying Theorem \ref{Adichprolon} with $X=H^s(S^1)$,
   we next extend these dichotomies in $H^s(S^1)$ to $\Zr^{\pm}$ and thus
   prove that
   $T_{\tilde{u}}((n+1)\tau,n\tau)$
   admits exponential dichotomies on $\Zr^{\pm}$ in
   $H^s(S^1)$ of finite rank equal to
   $i(e^{\pm})$, with exponent $\beta_1^{\pm}$, constant $M^{\pm}$ and
   projections $P_{\tilde{u}}^{\pm}(n)$, satisfying the property
\eqref{APnPpm}.
   Using next  \cite[Lemma 7.6.2 and Exercise 5 of Chapter 7]{He81}, we
   may extend these projections
$P_{\tilde{u}}^{\pm}(n)$ in $H^s(S^1)$ to projections
$P_{\tilde{u}}^{\pm}$ in $L^2(S^1)$ satisfying the properties
\eqref{APnPpm} and \eqref{PnL2}.\HB
An alternative proof consists in showing first that, for any $n \in \Zr^{\pm}$,
with $| n |$ large enough, we have the
exponential asymptotic convergence
   \begin{equation}
\label{BasexpB}
\|T_{\tilde{u}}((n+1)\tau,n\tau) - e^{L_{e^{\pm}}}\|_{L(L^2,H^s)}
\leq Ce^{- C |n|}~,
\end{equation}
where $C$ is a positive constant. Then we may apply Theorem
\ref{Adichpertur} in the space $L^2(S^1)$, which implies that there
exists $n_0  \in \Zr^{+}$
such that $T_{\tilde{u}}((n+1)\tau,n \tau)$ admits an exponential dichotomy on
   $\Zr^{+}$ (resp. $\Zr^{-}$) in $L^2(S^1)$, for $n \geq n_0$ (resp.
$n \leq -n_0$)
   of finite rank $i(e^{+})$ (resp. $i(e^{-})$) and projections
   ${\tilde P}^{+}_{T}$ (resp. ${\tilde P}^{-}_{T}$) satisfying the
properties of Theorem
\ref{Adichpertur}.
   Applying Theorem \ref{Adichprolon} with $X=L^2(S^1)$,
   we next extend these dichotomies in $L^2(S^1)$ to $\Zr^{\pm}$ and thus
   prove that  $T_{\tilde{u}}((n+1)\tau,n\tau)$
   admits exponential dichotomies on $\Zr^{\pm}$ in
   $L^2(S^1)$ of finite rank equal to
   $i(e^{\pm})$, with exponent $\beta_1^{\pm}$, constant $M^{\pm}$ and
   projections ${\tilde P}_{\tilde{u}}^{\pm}(n)$, satisfying the property
\eqref{APnPpm}. We remark that, by the property (ii) of the
definition \ref{Adefdich} and by the property \eqref{BasexpB}, the
image of ${\tilde P}_{\tilde{u}}^{\pm}(n)$ belongs to $H^s(S^1)$,
which implies, together with  \cite[Lemma 7.6.2 of Chapter 7]{He81},
that the restrictions $P_{\tilde{u}}^{\pm}$ of the
projections ${\tilde P}_{\tilde{u}}^{\pm}(n)$ to $H^s(S^1)$ define an
exponential dichotomy of $T_{\tilde{u}}((n+1)\tau,n\tau)$
on $\Zr^{\pm}$ in $H^s(S^1)$ of finite rank equal to
   $i(e^{\pm})$, with exponent $\beta_1^{\pm}$, constant $M^{\pm}$.
\end{demo}
\vskip 2mm

We notice that Theorem \ref{ATpmdicho} and Lemma \ref{Alemstar} imply
that $T^*(n, n+1)=T_{\tilde{u}}((n+1)\tau,n\tau)^*$ admits a reverse
exponential dichotomy on $\Zr^\pm$ in $L^2(S^1)$ (resp.
$H^{-s}(S^1)$) with rank $i(e^{\pm})$, exponent $\beta_1^{\pm}$ and
projections  $({\tilde P}_{\tilde{u}}^{\pm}(n))^*$ (resp.
$(P_{\tilde{u}}^{\pm}(n))^*$).
\vskip 2mm

Lemma 4.2 (on page 376) and Appendix C of \cite{ChChHa} yield the
  important characterization of the range of $P_{\tilde{u}}^{\pm}(n)$
  given in the next proposition.

   \begin{prop} \label{Atangent} We have the following equalities in
   $H^s(S^1)$,
\begin{equation}
\label{tangentP}
\begin{split}
R(P_{\tilde{u}}^-(n))  \,&=\, T_{\tilde{u}(n)} W^u(e^-)~, \quad \forall n \in
\Zr^{-}\cr
   R(I-P_{\tilde{u}}^+(n)) \, & = \, T_{\tilde{u}(n)} W^s(e^+)~, \quad
   \forall n \in \Zr^{+}~.
\end{split}
   \end{equation}
   \end{prop}

Let ${\tilde u}(t)$ belongs to
$W^u(e^{-}) \cap W^s(e^{+})$.
We say that the bounded orbit ${\tilde u}$ is  {\sl transverse}
at ${\tilde u}(0)$  if
$$
W^u(e^{-}) \pitchfork_{{\tilde u}(0)} W^s_{loc}(e^{+})~,
$$
which means that $T_{{\tilde u}(0)} W^u(e^{-})$ contains a closed
complement of $T_{{\tilde u}(0)}  W^s(e^{+})$ in $H^s(S^1)$ (notice
that, as $W^u(e^{-}) \cap W^s(e^{+})$ are immersed manifolds in
$H^s(S^1)$, this notion is well-defined, see \cite[Page 23]{Lang}).
It is easily seen that, since the linearized operator $T(t,\sigma)$
is injective and has dense range in $H^s(S^1)$, the above condition
implies that, for any $t \in \Rm$,$$
W^u(e^{-}) \pitchfork_{{\tilde u}(t)} W^s_{loc}(e^{+})~,
$$
which allows to simply say that the orbit ${\tilde u}$ is a {\sl
transverse} connecting orbit.

{}From the previous proposition and the equalities \eqref{kerQ*1} and
\eqref{kerQ*2} as well as the property that the range of
${\tilde P}_{\tilde{u}}^{\pm}$ is contained in $H^s(S^1)$, we at once deduce
the following equivalences.

\begin{prop} \label{Atrans1}
(i) The trajectory $\tilde{u}(t)$ is transverse in $H^s(S^1)$ if and only if
\begin{equation}
\label{CaractT1}
R(P_{\tilde{u}}^-(0)) + R(I-P_{\tilde{u}}^+(0)) \, = \, H^s(S^1)~,
\end{equation}
or equivalently,  since $R(P_{\tilde{u}}^-(0))$ is finite-dimensional
and, thus, this sum is closed,
\begin{equation}
\label{CaractT2}
[R(P_{\tilde{u}}^-(0))]^{\perp} \cap
[R(I- P_{\tilde{u}}^+(0))]^{\perp} \, = \, \{0 \}~,
\end{equation}
or also
\begin{equation}
\label{CaractT3}
R(I - P_{\tilde{u}}^-(0)^{\ast}) \cap
R(P_{\tilde{u}}^+(0)^{\ast}) \, = \, \{0 \}~.
\end{equation}
(ii) Moreover, the trajectory $\tilde{u}(t)$ is transverse in
$H^s(S^1)$ if and only if
\begin{equation}
\label{Caract4}
R({\tilde P}_{\tilde{u}}^-(0)) + R(I- {\tilde P}_{\tilde{u}}^+(0)) \,
= \, L^2(S^1)~,
\end{equation}
or equivalently,
\begin{equation}
\label{CaractT5}
[R({\tilde P}_{\tilde{u}}^-(0))]^{\perp} \cap
[R(I-{\tilde P}_{\tilde{u}}^+(0))]^{\perp} \, = \, \{0 \}~,
\end{equation}
or also
\begin{equation}
\label{CaractT6}
R(I - {\tilde P}_{\tilde{u}}^-(0)^{\ast}) \cap
R({\tilde P}_{\tilde{u}}^+(0)^{\ast}) \, = \, \{0 \}~.
\end{equation}
\end{prop}

Applying Lemma \ref{Tstarpsi0} and Proposition \ref{Atrans1}, we
obtain the next characterization of transversality.
\begin{coro} \label{Atrans2}
The trajectory $\tilde{u}(t)$ is transverse
if and only if there does not exist any element $\psi_0 \in (H^s(S^1))^\ast$,
$\psi_0 \ne 0$, such that the sequence $(T_{\tilde{u}}^\ast(n,0) \psi_0)_{n \in
\Zr}$ is bounded in $(H^s(S^1))^\ast$ or equivalently, if and only if
there does not
exist any element $\psi_1 \in L^2(S^1)$,
$\psi_1 \ne 0$, such that the sequence $(T_{\tilde{u}}^*(n,0) \psi_1)_{n \in
\Zr}$ is bounded in $L^2(S^1))$
\end{coro}

Finally we show how to apply Theorem \ref{AFredholm} to obtain a
corollary, which plays a crucial role in the proof of the genericity
of the non-existence of homoindexed  orbits between equilibrium points.
We introduce the operator ${\cal L}_{\tilde u} \equiv \cal L$ defined
in  \eqref{AcalL} with
$T(n,m)=T_{\tilde u}(n\tau,m\tau)$, $X=H^s(S^1)$ and ${\cal
Z}=\ell^\infty(\Zr,X)$. Likewise, we introduce the extension ${\tilde {\cal
L}}_{\tilde u}$ of ${\cal L}_{\tilde u}$ to ${\tilde {\cal
Z}} =\ell^\infty(\Zr,L^2(S^1))$. We notice that, due to the smoothing
properties of  $T_{\tilde u}(n\tau,m\tau)$, a sequence $F \in {\cal
Z}$ belongs to
$R({\cal L}_{\tilde u})$ if and only if $F  \in {\cal Z}$ belongs to
$R({\tilde {\cal L}}_{\tilde u})$. This remark, Theorem
\ref{AFredholm} and Corollary \ref{Atrans2} imply the following
results.

\begin{coro}\label{ApFredB}{\bf Functional characterization of
transversality} \HB
The operators ${\cal L}_{\tilde u}:D({\cal L}_{\tilde u})\rightarrow {\cal
Z}$ and ${\tilde {\cal L}}_{\tilde u}: D({\tilde {\cal L}}_{\tilde
u})\rightarrow
{\tilde {\cal Z}}$ defined above are Fredholm operators of index
$i(e^-)-i(e^+)$. In particular, the codimension of $R({\cal
L}_{\tilde u})$ in $H^s(S^1)$ and of $R({\tilde {\cal L}}_{\tilde
u})$ in $L^2(S^1)$ is equal to
$\hbox{codim }[R({\tilde P}_{\tilde{u}}^-(0)) + R(I - {\tilde
P}_{\tilde{u}}^+(0))]$.
\HB
Moreover, a sequence $F\in {\cal Z}$ belongs to $R({\cal L}_{\tilde u})$ if and
only if
$$
\sum_{n=-\infty}^{+\infty} \langle F(n),\Psi(n+1)\rangle_{L^2(S^1)}
\,=\, 0~,
$$
for every bounded sequence $\Psi(n)=T^*_{\tilde u}(n\tau,0)\Psi_0$, $\Psi_0\in
L^2(S^1)$. Finally, the connecting orbit $\tilde u(t)$ is transverse
if and only if  ${\cal L}_{\tilde u}$ is surjective or also if and
only if ${\tilde {\cal L}}_{\tilde u}$ is surjective.
\end{coro}


\section{Asymptotics of solutions of perturbations of linear
autonomous equations}\label{AppendixD}

In this section, $X$ denotes a Banach space and $J$ an interval of
$\Zr$.

In several proofs of the previous sections, we need to know the
asymptotics of the bounded solutions (as $t \rightarrow \pm \infty$) of
the linearized equations along orbits, connecting hyperbolic
equilibrium points or periodic orbits of \eqref{eq}. These
asymptotics will be described in the second and third sections of
this appendix.

\subsection{Abstract results}

Here, we are
going to describe these asymptotics for a general linearized equation
or for iterates of a general linearized mapping.
Thus, in the
first place, we are interested in the asymptotic behaviour of the
bounded sequences $u(n)$, $n \in \Zr^{+}$ (resp. $n \in \Zr^{-}$),
defined by
\begin{equation}
\label{CTSigma}
u(n+1) = Tu(n) + \Sigma(n) u(n) \equiv L(n) u(n) ~,
\end{equation}
where $T \in L(X)$ and $\Sigma(n) \in L(X)$, for any $n \in \Zr^{+}$
(resp. $n \in \Zr^{-}$). \HB

All the statements given in this appendix are already known and are
mainly results or generalizations of results  of \cite{He85},
\cite{He94}, \cite{ChChHa} and \cite{BrunRaug}. The main theorems
\ref{C6} and \ref{C8}
are a refinement of Theorem B.5 of \cite{ChChHa} and have been
proved in Appendix B of \cite{BrunRaug}. Here, we closely follow
Appendix B of \cite{BrunRaug}. We want to point out that all the
statements contained in this appendix are more or less common
knowledge, at least in finite dimensions. For additional references, see
\cite{Co}, \cite{Pal84}, \cite{Pal88}, \cite{Pal00},
\cite{HaLi86} and \cite{S78} for example.

For $n \geq m$ in $J$, we define the evolution
operator  $\Phi(n,m)= L(n-1) \circ \ldots \circ L(m)$.

For any linear operator $T:
X \rightarrow X$, we
denote by ${\cal R}(T)$ the set of all nonnegative numbers $\rho$
for which
$$ \sigma (T) \cap \{z \in \Com \, | \, |z| = \rho \}
\ne \emptyset~.
$$
For $\rho \notin {\cal R}(T)$, we denote by
$P_{\rho}$, $Q_{\rho}$ the spectral projections associated with
the partition of the spectrum $\sigma(T)$ into its subsets $\sigma(T)
\cap \{|z| >\rho\}$ and $\sigma(T) \cap \{|z| <\rho\}$
respectively.  \HB
The following proposition gives a sufficient condition for $L(t)$
defined by \eqref{CTSigma} to admit a shifted dichotomy.

\begin{prop} \label{C2} Let $L(n)= T +\Sigma(n)$ where
$T$, $\Sigma(n)$, $n \in \Zr^{+}$, belong to $L(X)$ and
$\Sigma(n)$ satisfies the asymptotic condition
\begin{equation}
\label{CSigrho}
\|\Sigma(n)\|_{L(X)} \, = \, O(r^n)~,
\hbox{ for some }r <1, \hbox{ when }n\rightarrow +\infty~.
\end{equation}
Assume that $0 < \rho_1^*<\rho_1 \leq \rho_2 <
\rho_2^\ast$ are such that $\sigma(T) \cap \{\rho_1^\ast
\leq |z| \leq \rho_2^\ast\} = \emptyset$, which implies that
$Q_{\rho_2^\ast}-Q_{\rho_1^\ast}=0$.  Suppose also that $T$
admits shifted exponential dichotomy with gap $[\rho_1^\ast,
\rho_2^\ast]$ and with projections $P_{\rho_1^\ast}$,
$Q_{\rho_1^\ast}$. Then
the family $L(\cdot)$ admits shifted dichotomy
on $J=\Zr^{+}$ with gap
$[\rho_1,\rho_2]$ and projections $P(n)$, $Q(n)$.  Moreover,
when $n \in \Zr^{+}$ is sufficiently large,
\begin{equation}
\label{CQtQ1}
\|Q(n) -
Q_{\rho_1}\|_{L(X)} \, = \, O(r^n)~, ~ (\hbox{ resp. }
\|P(n) - P_{\rho_2}\|_{L(X)} \, = \, O(r^n))~,
\end{equation}
and $\tilde Q(n) := Q(n)_{\big|R(Q_{\rho_1})}: R(Q_{\rho_1})
\rightarrow R(Q(n))$ (resp.  $\tilde Q_{1}(n)=Q_{\rho_1\big|R(Q(n)}: R(Q(n))
\rightarrow
R(Q_{\rho_1})$) is an isomorphism satisfying
\begin{equation}
\label{CQtQ2}
\begin{split}
\max
(\| \tilde Q(n)- I\|_{L(X)}, \| \tilde Q_1(n)- I\|_{L(X)}) \,
&\leq \, \|Q(n) - Q_{\rho_1}\|_{L(X)}=O(r^{n})~, \cr
\max(
\|{\tilde Q}^{-1}(n) - I\|_{L(X)},  \|{\tilde Q}_1^{-1}(n) -
I\|_{L(X)}) \, &\leq \, \frac {\|Q(n) - Q_{\rho_1}\|_{L(X)}} {1
-\|Q(n) - Q_{\rho_1}\|_{L(X)}}=O(r^n)~.
\end{split}
\end{equation}
The
same statement holds on $J=\Zr^{-}$ if the condition ``$n
\in \Zr^{+}$" is replaced by ``$n \in \Zr^{-}$" and
$r <1$ by $r >1$.
\end{prop}

We next recall three results about the asymptotic behaviour of $u(n)$ for $n$
large enough, where $u(n)$ is given by \eqref{CTSigma}. The first theorem
has been proved by D. Henry (see Theorem 2 in \cite{He85}) ; here we
state it under the form given by Chen, Chen and Hale in
\cite[Theorem B.2]{ChChHa}.

\begin{theorem} \label{CHenryQP}
   Let $T\in L(X,X)$  and let $0<
\rho_1 \leq \rho_2$ be such that $[\rho_1,\rho_2] \cap
{\cal R}(T) = \emptyset$.
Let $u(n) \ne 0$, $n \geq n_0$ in $\Zr^{+}$, be a sequence in
$X$ such that
\begin{equation}
\label{CHenryQP1}
\lim_{n \rightarrow +\infty}{ \frac{\|u(n+1)
-Tu(n)\|_{X}}{\|u(n)\|_{X}}} \, =\, 0~.
\end{equation}
Then, either
$$
(i) \lim_{n \rightarrow +\infty} \frac{\|Pu(n)\|_{X}} {\|Qu(n)\|_{X}} \, = \, +
\infty~, \hbox{ and }\quad
\liminf_{n \rightarrow +\infty} \|u(n)\|_{X}^{1/n} \geq
\rho_2~;
$$
or
$$
(ii) \lim_{n \rightarrow +\infty}
\frac{\|Pu(n)\|_{X}}{\|Qu(n)\|_{X}} \, = \, 0~, \hbox{ and }\quad
\limsup_{n \rightarrow +\infty} \|u(n)\|_{X}^{1/n} \leq
\rho_1~;
$$
where $P=P_{\rho_1}=P_{\rho_2}$ and
$Q=Q_{\rho_1}=Q_{\rho_2}$. The same property holds if, in
the above statements, $n \geq n_0$ in $\Zr^{+}$ and $n
\rightarrow +\infty$ are replaced by $n \leq -n_0$ in
$\Zr^{-}$ and $n \rightarrow -\infty$ respectively.
\end{theorem}

Theorem \ref{CHenryQP} gives a close relation between the spectrum of $T$ and
the growth rate of $u(n)$. As direct corollary, Chen, Chen and Hale
(see \cite[Corollary B.3]{ChChHa}) have proved the following property.

\begin{theorem} \label{Cvitesseconv} Let $T$ be a continuous linear operator
from $X$ into $X$ such that ${\cal R}(T)$ is nowhere dense in $[0,
+\infty)$ and let $u(n) \ne 0$, $n \geq n_0$ in
$\Zr^{+}$, be a sequence
in $X$ satisfying the property \eqref{CHenryQP1}. Then, there exists
$\rho \in {\cal R}(T)$ such that
\begin{equation}
\label{Cvitessconv}
\lim_{n \in \Zr^{+} \rightarrow + \infty} \|u(n)\|_{X}^{1/n}
\, = \, \rho~.
\end{equation}
The same property holds when $\Zr^{+}$ is replaced by
$\Zr^{-}$.
\end{theorem}

\begin{rem} For sequences $u(n)$, $n \in \Zr^{+}$, given by the
recursion formula
\eqref{CTSigma} with $\Sigma(n)$ satisfying the hypothesis
\eqref{CSigrho}, the condition \eqref{CHenryQP1} obviously holds. We
thus deduce from Theorem \ref{Cvitesseconv} that there exists $\rho \in
{\cal R}(T)$ such that $\lim_{n \in \Zr^{+} \rightarrow + \infty}
\|u(n)\|_{X}^{1/n} = \rho$.
\end{rem}

In this paper, we cannot directly apply this corollary since
${\cal R}(T)$ has an accumulation point at $0$. However, we know that
the sequences $u(n)$ that we will consider do not converge faster to
zero than an exponential.

For any $\psi \in X$, for any integers $m$, $n$, with $n \geq m$,
following the notations of \cite{ChChHa}, we set,
$$
u(n,m;\psi) =\Phi (n,m)\psi ~,
$$
where $\Phi (n,m) = L(n-1) \circ \ldots \circ L(m)$ and $L(n)$ is
defined by \eqref{CTSigma}. We also introduce the quantity
$$
r_{\infty} (m, \psi) = \limsup_{n \rightarrow +\infty} \|u(n,m;
\psi)\|_{X}^{1/n}~.
$$
Let $T$ be a continuous linear operator from $X$ into $X$ such that
${\cal R}(T)$ is a bounded sequence converging to $0$ and that
$\rho_{j+1} \leq
\rho_j$, for any integer $j$. We then introduce the spaces
$$
E_j^+(m)=\{\psi \in X \, | \, r_{\infty} (m, \psi) \leq
\rho_j\}~.
$$
Then,
$$
\ldots E_2^+(m) \subset E_1^+(m) \subset E_0^+(m) ~.
$$
If we assume that, for any $\psi \in X$, $\psi \ne 0$ and any integer
$m$, $r_{\infty} (m, \psi) >0$ (that is $\Phi (n,m)\psi$ does not
decay faster to $0$ than an exponential), then
$$
\bigcap_{j=0}^{j= \infty} E_j^+(m) = \{ 0\}~.
$$
Taking into account the above considerations and following the proof
of Corollary B.3 of \cite{ChChHa}, we obtain the following theorem.

\begin{theorem} \label{C5}

Let $T$ be a continuous linear operator
from $X$ into $X$ such that ${\cal R}(T)$ is a bounded sequence
converging to $0$ and that $\rho_{j+1} \leq
\rho_j$, for any integer $j$. We also
      assume that the property \eqref{CSigrho} holds.
Let  $u(n) \equiv u(n,m;\psi)$,
      be a sequence in $X$ such that $r_{\infty} (m, \psi) >0$.
Then, there exists
$\rho \in {\cal R}(T)$ such that
\begin{equation}
\label{CvitessC5}
\lim_{n \in \Zr^{+} \rightarrow + \infty} \|u(n)\|_{X}^{1/n}
\, = \, \rho~.
\end{equation}
The same property holds when $\Zr^{+}$ is replaced by
$\Zr^{-}$. \HB
In particular, if for any $\psi \in X$, $\psi \ne 0$ and any integer
$m$, $r_{\infty} (m, \psi) >0$, then
$$
X = E_0^+(m) = \{\psi \in X \, | \, r_{\infty} (m, \psi) \leq
\rho_1\}~,
$$
and
$$
E_j^+(m) - E_{j+1}^+(m) =\{\psi \in X \, | \, r_{\infty} (m, \psi)
= \rho_{j+1}\}~.
$$
\end{theorem}

The next theorem
is nothing else as Theorem B.6 of \cite{BrunRaug}
and is actually a refinement of Theorem B.5 of
\cite{ChChHa} about the asymptotics of sequences
$u(n)$ given by the recurrence formula
\eqref{CTSigma} when $n$ goes to $\pm\infty$.

\begin{theorem} \label{C6}{\bf Convergence to a solution of the asymptotic
equation}

Let $T\in L(X,X)$
and suppose that there exist positive numbers
$\delta_1$, $\tilde{\delta}_1$, $\delta$, $\tilde{\delta}$, with
$0 < \delta_1 < \delta$ and $0 < \tilde{\delta}_1 <
\tilde{\delta}$, and $\rho \in {\cal R}$ such that
\begin{equation}
\label{Clambdadelta}
\emptyset \, \ne \, \sigma(T) \cap \{z \in \Com \, | \, \rho -\delta
\leq |z| \leq \rho +\tilde{\delta}\} \subset \{z \in \Com\, |
\, \rho
-\delta_1 < |z|  < \rho +\tilde{\delta}_1\}~.
\end{equation}
Suppose also that $T$ admits shifted dichotomy with gap $[\rho
+\tilde{\delta}_1,\rho + \delta^\ast]$, for some $\delta^\ast
>\tilde{\delta}$ (resp. with gap $[\rho -\delta^\ast,
\rho - \delta_1]$, for some $\delta^*> \delta$).
Let $u(n)$, $n \in \Zr^{+}$ (resp. $n\in\Zr^-$) be a sequence given by the
recurrence formula
\eqref{CTSigma},  with $\Sigma(n)$ satisfying the hypothesis
\eqref{CSigrho}, where $(\rho +\tilde{\delta})r <\rho
-\delta$ (resp. where $(\rho - \delta)r
<\rho +\tilde{\delta}$), such that
\begin{equation}
\label{Climdelta1u}
\rho -\delta_1 \leq \lim_{n \rightarrow +\infty} \|u(n)\|_{X}^{1/n}
\leq  \rho + \tilde{\delta}_1~~~~\text{(resp. }\rho -\delta_1 \leq \lim_{n
\rightarrow -\infty} \|u(n)\|_{X}^{1/n}
\leq  \rho +  \tilde{\delta}_1\text{)}~.
\end{equation}
Denote by $T_{\rho}$ the operator
$$
T_{\rho} \, = \, [Q_{\rho +\tilde{\delta}} - Q_{\rho -\delta}] T
[Q_{\rho + \tilde{\delta}} -Q_{\rho -\delta}]~.
$$
Then, there exists a non-vanishing sequence $u_{+\infty}(n)$,
$n\in\Zr^+$, (resp.
$u_{-\infty}(n)$, $n\in\Zr^-$,) in $R(Q_{\rho
+\tilde{\delta}} - Q_{\rho -\delta})$ and satisfying
\begin{equation}
\label{Cuinfini1}
u_{+\infty}(n +1) \, = \, T_{\rho}u_{+\infty}(n)~~~~\text{(resp.}
u_{-\infty}(n +1) \, = \, T_{\rho}u_{-\infty}(n)\text{)}~,
\end{equation}
and
\begin{align}
\label{Cuinfini2}
\|u(n) - u_{+\infty}(n)\|_{X} \, &= \, O((\rho - \delta)^{n})~, \quad
\hbox{ as }n \longrightarrow + \infty\\
\text{(resp. }~\|u(n) - u_{-\infty}(n)\|_{X} \, &= \, O((\rho +
\tilde\delta)^{n})~, \quad \hbox{ as }n \longrightarrow - \infty \text{ ).}
\label{Cuinfini2-}
\end{align}
\end{theorem}

The previous theorem allows to specify the asymptotics of the
bounded sequences $u(n)$ given by \eqref{CTSigma}.
As consequences of Theorem \ref{C6}, one obtains
corresponding results for solutions of evolutionary partial differential
equations. More precisely, let $Y$ be a Banach space and $A$ be the
infinitesimal
generator of an analytic semigroup on $Y$. Let $\alpha \in [0,1)$ be a real
number. We introduce the fractional space $X=Y^{\alpha}$ and consider
      the equation
\begin{equation}
\label{CAGU}
         \partial_t U(t)  \,=\, (A + G(t))U(t)~, \quad t>0~,
         \quad U(0) \,=\, U_0~,
         \end{equation}
where $U(t)$, $t \geq0$, and $U_0$ belong to $X$,
and $G: t \in \Rm \mapsto G(t) \in L(Y^{\alpha}, Y)$
is such that
\begin{equation}
\label{CGr}
\|G(t)\|_{L(Y^{\alpha}, Y)} \, = \, O(e^{-rt}) ~, \quad \hbox{ as }t \in \Rm
\rightarrow +\infty~, \hbox{ where }r >0~.
\end{equation}
If $\alpha=0$, it suffices to assume that $A$ is the generator of a
$C^0$- semigroup.

The proof of the next theorem, which is a consequence of Theorem
\ref{C6}, follows the lines of the proof of Theorem B.8 of
\cite{BrunRaug}.

\begin{theorem} \label{C8}
Suppose that there exist positive
constants $d_1$,
$\tilde{d}_1 $, $d$, $\tilde{d}$ with $0< d_1 < d$, $0 < \tilde{d}_1
< \tilde{d}$, and $\mu \in \Rm$  such that
      \begin{equation}
      \label{Cmud}
\emptyset \, \ne \, \sigma(A) \cap \{z \in \Com\, | \, \mu -d
\leq \Re z \leq \mu + \tilde{d} \} \subset \{z \in \Com\, |
\, \mu -d_1 < \Re z  < \mu + \tilde{d}_1\}~.
\end{equation}
Suppose also that $e^A$ admits shifted dichotomy
      with gap $[e^{\mu + \tilde{d}_1}, e^{\mu +d^\ast}]$,
for some $d^*>\tilde{d}$ (resp. with gap $[e^{\mu -d^\ast},
e^{\mu -d_1}]$, for some $d^*>d$). Let $U(t)$, $t \in \Rm^{+}$ (resp. $t\in\Rm
^-$), be a
solution of \eqref{CAGU} with $G(t)$ satisfying the hypothesis
\eqref{CGr} and with $e^{(\mu + \tilde{d}-r)} <e^{\mu -d}$ (resp. with $e^{(\mu
- d -r)} > e^{\mu +\tilde{d}}$), such
that
\begin{equation}
\label{CBlimd1U}
\mu - d_1 \leq \lim_{t \rightarrow +\infty} \ln
(\|U(t)\|_{X}^{1/t}) \leq  \mu + \tilde{d}_1 ~.
\end{equation}
Denote
       $$
       A_{\mu} \, = \, [Q_{\mu +\tilde{d}} - Q_{\mu -d}] A
[Q_{\mu + \tilde{d}} -Q_{\mu -d}]~,
$$
where now $Q_{\mu
+\tilde{d}}$ and $Q_{\mu -d}$ denote the spectral projections
associated with the parts of the spectrum $\sigma(A) \cap \{\Re z >
\mu +\tilde{d}\}$ and $\sigma(A) \cap \{\Re z > \mu - d\}$.
\HB
Then, there exists a non-vanishing solution $U_{+\infty}(t)$, $t\in\Rm^+$,
(resp. $U_{-\infty}(t)$, $t\in\Rm^-$,) in $R(Q_{\mu +\tilde{d}} - Q_{\mu
-d})$ and satisfying
\begin{equation}
\label{CUmuinfini1}
        \partial_t U_{+\infty}(t)  \,=\, A_{\mu}U_{+\infty}(t)~~~~~\text
{ (resp. }~\partial_t U_{-\infty}(t)  \,=\, A_{\mu}U_{-\infty}(t)~\text{ )}~,
\end{equation}
and
\begin{align}
\label{CUmuinfini2}
\|U(t) - U_{+\infty}(t)\|_{X} \, &= \, O(e^{(\mu -d)t})~, \quad
\hbox{ as }t \longrightarrow + \infty~.\\
\text{(resp. }~\|U(t) - U_{-\infty}(t)\|_{X} \, &= \, O(e^{(\mu
+\tilde d)t})~, \quad
\hbox{ as }t \longrightarrow - \infty~\text{ )}~.
\label{CUmuinfini2-}
\end{align}
\end{theorem}

\subsection{Applications to the parabolic equation near an
equilibrium point}

The results of Appendix B and those of the first part of this appendix,
together with those of Section \ref{subsecLin},
will be applied here in order to determine the asymptotics of the solutions
$u(t)$ of \eqref{eq}, which belong to the local unstable or stable
manifolds of hyperbolic equilibria or periodic orbits,
as well as the asymptotics of the solutions of the corresponding
linearized equations.
\vskip 2mm

Let $e$ be a hyperbolic equilibrium point of \eqref{eq}.  In
accordance with Section \ref{subsecLin}, we denote by $L_e$ the
corresponding linearized operator and by $\lambda_i$, $i \geq 1$ its
eigenvalues, counted with their multiplicity.

\begin{coro}\label{CappD0}
Let $u(t)$ be a trajectory of \eqref{eq} belonging to the unstable
manifold $W^u(e)$ of $e$ (resp.  the local stable manifold $W^s_{loc}(e)$) and
let $v(t) =u(t) -e$. \HB
Then, there exist $(a,b)\in\Rm^2-\{(0,0)\}$ and an eigenvalue
$\lambda_i$ of $L_e$ such that $\Re(\lambda_i)>0$ (resp.
$\Re(\lambda_i)<0$) and
$$
\lim_{t \to -\infty} \ln \|v(t)\|_{H^s(S^1)} ={\Re(\lambda_i)} \quad
\hbox{( resp. } \lim_{t \to +\infty} \ln \|v(t)\|_{H^s(S^1)} ={\Re(\lambda_i)}
\hbox{)}~.
$$
More precisely, the asymptotic behavior of $v$ in $H^s(S^1)$ is as
follows:\\
(i) if $\lambda_i$ is a simple real eigenvalue with eigenfunction $\varphi_i$,
then there exists $a\in\Rm-\{0\}$ such that $v(t)=a e^{\lambda_i t}\varphi_i
+ o(e^{\lambda_i t})$.\\
(ii) If $\lambda_i=\lambda_{i+1}$ is a double real eigenvalue with two
independent eigenfunctions $\varphi_i$ and $\varphi_{i+1}$, then there exist
$(a,b)\in\Rm^2-\{(0,0)\}$ such that $v(t)=a e^{\lambda_i
t}\varphi_i + be^{\lambda_i t}\varphi_{i+1} + o(e^{\lambda_i t})$. \\
(iii) If $\lambda_i=\lambda_{i+1}$ is an algebraically double real eigenvalue
with eigenfunction $\varphi_i$ and with generalized
eigenfunction $\varphi_{i+1}$, then there exist
$(a,b)\in\Rm^2-\{(0,0)\}$ such that $v(t)=(a+bt)e^{\lambda_i
t}\varphi_i + be^{\lambda_i t}\varphi_{i+1} + o(e^{\lambda_i t})$.\\
(iv) If $\lambda_{i +1}=\overline{\lambda_{i}}$ is a (simple) nonreal 
eigenvalue
with eigenfunction $\varphi_{i +1}=\overline{\varphi_{i}}$, then there exist
$(a,b)\in\Rm^2-\{(0,0)\}$ such
that $v(t)=e^{\Re(\lambda_i)t} \big[(a\cos(\Im(\lambda_i)t) - b \sin(
\Im(\lambda_i) t))
\Re(\varphi_i)  - (a\sin( \Im(\lambda_i) t) + b \cos(\Im(\lambda_i)t))
\Im( \varphi_i )) \big] + o(e^{\Re(\lambda_i) t})$.

Let $j\in\Nm \setminus \{0\}$ be such that $\lambda_i$ belongs to the pair of
eigenvalues $(\lambda_{2j-1},\lambda_{2j})$, or let $j=0$ if
$\lambda_i=\lambda_0$.  As a consequence of the asymptotic behaviour,
there exists $t_0\in\Rm$ such that, for all $t\leq t_0$ (resp.  $t\geq
t_0$), $v(t)$ has exactly $2j$ zeros which are simple.
\end{coro}

\begin{demo}
Since the proofs are very similar when $t$ tends to $\pm \infty$, we
shall only prove the corollary when $u(t)$ belongs to the local
stable manifold $W^s_{loc}(e)$. To prove the corollary, we shall apply
Theorem \ref{C8}, so we have to check that the hypotheses of Theorem
\ref{C8} are satisfied. \HB
Since $u(t)$ belongs to
$W^s_{loc}(e)$, there exist two positive constants $c_1$ and $\kappa$
such that
    \begin{equation}
\label{Dconv1}
\| u(t) -e \|_{H^s} \leq c_1 e^{-\kappa t}~, \quad \hbox{ as } t
\rightarrow +\infty ~.
\end{equation}
The function $v(t) \equiv u(t) -e$ is a classical solution of the
equation
    \begin{equation}
\label{Deqv}
v_t = v_{xx} + D_{u_x}f(x,e,e_x)v_x + D_uf(x,e,e_x)v + a(x,t)v_x +
b(x,t)v~,
\end{equation}
where
    \begin{equation*}
\begin{split}
a(x,t) =& \int_{0}^{1} \big( f'_{u_x}(x, e+\theta (u-e),e_x+\theta
(u_x-e_x)) -f'_{u_x}(x, e,e_x) \big) d\theta \cr
b(x,t)= & \int_{0}^{1} \big( f'_{u}(x, e+\theta (u-e),e_x+\theta
(u_x-e_x)) -f'_{u}(x, e,e_x) d\theta \big)
\end{split}
\end{equation*}
One at once checks that $\| a(x,t)\|_{C^0} + \| b(x,t)\|_{C^0}  \leq
c_2 \| u(t) -e
\|_{H^s}$, which implies that
    \begin{equation}
\label{Dconv2}
\| a(x,t)v_x + b(x,t)v \|_{L^2} \leq c_3  e^{-\kappa t} \|v \|_{H^s}
\end{equation}
Thus, $v$ is the solution of an equation of the form \eqref{CAGU},
with $G(t) v=a (x,t)v_x + b(x,t)v$ satisfying the condition
\eqref{CGr}. \HB
We next remark that
$$
v(n+1)= Tv(n) + \Sigma (n) v(n)~,
$$
where $T= T(1)$, $\Sigma (t)v(t) = \int_0^1 T(1-\sigma)
G(t +\sigma)v(t + \sigma) d\sigma$ and $T(t)$ is the linear semigroup
associated with the linear equation
$v_t = v_{xx} + f'_{u_x}(x,e,e_x)v_x + f'_u(x,e,e_x)v \equiv L_e v$. 
We next verify
that $\Sigma (n)$ satisfies the condition \eqref{CSigrho}, for some
$r <1$, when $n$ goes to infinity. First, $\Sigma (n)$ is a
continuous linear operator from $H^s(S^1)$ into $H^s(S^1)$. Moreover,
since $T(t)$ is an analytic linear semigroup, we obtain the following
inequality, for
$0 <\tau \leq 1$,
    \begin{equation}
\label{Dcoaux1}
\| v(t+\tau)\|_{H^s} \leq C e^{\alpha \tau}\| v(t)\|_{H^s} + C
\int_{0}^{\tau} e^{\alpha(\tau
-\sigma)} (\tau-\sigma)^{-s/2} \| v(t+\sigma)\|_{H^1}~,
\end{equation}
where $ \alpha >0$. Using a generalized Gronwall inequality (see \cite[Lemma
7.1.1]{He81}), we deduce from \eqref{Dcoaux1} that, for $0 <\tau \leq 1$,
    \begin{equation}
\label{Dcoaux2}
\| v(t+\tau)\|_{H^s} \leq C e^{(\alpha +K)\tau} \| v(t)\|_{H^s}~,
\end{equation}
where $K$ is a positive constant.
The definition of $\Sigma(n)$ and the properties \eqref{Dconv2} and
\eqref{Dcoaux2} imply
that, for $n >0$ large enough,
    \begin{equation}
\label{Dcoaux3}
\|\Sigma (n) v(n)\|_{H^s} \leq C \int_0^1
e^{\alpha(1-\sigma)}(1-\sigma)^{-s/2} c_3e^{-\kappa n} \| v(n
+\sigma)\|_{H^1}
\leq C^*  e^{(\alpha +K)}e^{-\kappa n} \| v(n)\|_{H^s}~,
\end{equation}
and hence $\Sigma(n)$ satisfies the property \eqref{CSigrho}. \HB
Since $\Sigma(n)$ satisfies the property \eqref{Dcoaux3} and that $T$
admits a shifted exponential dichotomy on $\Zr^+$, the family $L(\cdot)= T
+\Sigma(\cdot)$ admits a shifted dichotomy on $\Zr^+$ by Proposition
\ref{C2}. As, by a result of Agmon (see \cite{Agmon}), every non-zero
solution of a linear parabolic equation does not go to zero faster
than an exponential when $t$ tends to infinity, the hypotheses of
Theorem \ref{C5} are satisfied, that is, for any integer $m$ and any
$\psi \in H^s(S^1)$, $r_{\infty}(m, \psi) >0$. Hence, by Theorem
\ref{C5}, there exists $\rho_i \in {\cal R}(T)$, $\rho_i <1$, (and thus $\mu_i$
belonging to the spectrum of $T$), such that
$$
\lim_{n \to \infty} \| v(n)\|_{H^s}^{1/n}= \rho_i = | \mu_i| ~,
$$
or, in other terms, there exists an eigenvalue $\lambda_i$ of the
linearized operator $L_e$ such that
$$
\lim_{t \to \infty} \ln (\| v(t)\|_{H^s}^{1/t})=| \Re (\lambda_i)| ~.
$$
We can now apply Theorem \ref{C8}. Let $E_i$ be the generalized
(real) eigenspace associated with the eigenvalue $\lambda_i$ and
$L_{e,i}$ be the restriction of the operator $L_e$ to this eigenspace
$E_i$. By Theorem \ref{C8}, there exists a nonvanishing solution
$\psi_{\infty, i}(t) \in E_i$ of the equation
\begin{equation}
\label{Dcoaux4}
\partial_t \psi_{\infty, i}(t) = L_{e,i} \psi_{\infty, i}(t)~,
\end{equation}
such that
\begin{equation}
\label{Dcoaux5}
\| v(t) - \psi_{\infty, i}(t) \|_{H^s}=o(e^{\Re (\lambda_i) t})~.
\end{equation}
Now the corollary is an elementary consequence of the properties
\eqref{Dcoaux4} and \eqref{Dcoaux5} and of Proposition
\ref{spectreLe}. According to this proposition,
    $\lambda_i$ is either a simple real eigenvalue (and the dimension of
    $E_i$ is one), or a double real eigenvalue or a simple
non-real eigenvalue (in which cases, the dimension of $E_i$ is equal
to two).
\end{demo}
\vskip 2mm

In the same way, we prove the following corollary.

\begin{coro}\label{CappD1}
Let $u(t)$ be a trajectory of \eqref{eq} belonging to the unstable
manifold $W^u(e)$ of $e$ (resp.  the local stable manifold
$W^s_{loc}(e)$).  Let
$v_0\in T_{u(0)}W^u(e)$ (resp.  $v_0\in T_{u(0)}W^s_{loc}(e)$) and let
$v(t)$ be the solution for $t\leq 0$ (resp. $t\geq 0$) of the linearized
equation
\begin{equation}\label{eq-AppD}
v_t=v_{xx}+D_uf(x,u,u_x)v+ D_{u_x}f(x,u,u_x)v_x~, \quad v(0)=v_0~.
\end{equation}
Then, there exist $(a,b)\in\Rm^2-\{(0,0)\}$ and an eigenvalue
$\lambda_i$ of $L_e$ such that $\Re(\lambda_i)>0$ (resp.
$\Re(\lambda_i)<0$) and
$$
\lim_{t \to -\infty} \ln \|v(t)\|_{H^s(S^1)} ={\Re(\lambda_i)} \quad
\hbox{( resp. } \lim_{t \to +\infty} \ln \|v(t)\|_{H^s(S^1)} ={\Re(\lambda_i)}
\hbox{)}~.
$$
Moreover, all the possible asymptotic behaviors are the same as those
described in Corollary \ref{CappD0}. \HB
In addition, let $j\in\Nm \setminus \{0\}$ be such that $\lambda_i$ 
belongs to the pair of
eigenvalues $(\lambda_{2j-1},\lambda_{2j})$, or let $j=0$ if
$\lambda_i=\lambda_0$.  Then,
there exists $t_0\in\Rm$ such that, for all $t\leq t_0$ (resp.  $t\geq
t_0$), $v(t)$ has exactly $2j$ zeros which are simple.
\end{coro}

\begin{demo} If $u(t)$ belongs to the local stable manifold
$W^s_{loc}(e)$, $u(t)$ satisfies the property \eqref{Dconv1} and
Eq.  \eqref{eq-AppD} can be written in the form \eqref{Deqv},
where the functions $a(x,t)$ and $b(x,t)$ satisfy the properties
\eqref{Dconv2}. We remark that, by \cite[Theorem C2]{ChChHa}, we
already know
that $\limsup_{n \to \infty} \| v(n)\|_{H^s}^{1/n} <1$.
We thus obtain the asymptotic behavior of $v(t)$ by
following the lines of the proof of Corollary \ref{CappD0}.
\end{demo}

\subsection{Application to the parabolic equation near a periodic
orbit}\label{App-per}

Before proving analogous corollaries in the case of periodic orbits,
we briefly recall the known properties of the local stable and
unstable manifolds of hyperbolic periodic orbits.

Let $\gamma(x,t)$ be a hyperbolic periodic solution of \eqref{eq} of
minimal period $p$ and let $\Gamma=\{\gamma(t),~t\in\Rm\}$.
As in the introduction and in section 2, we introduce the linearized
equation \eqref{linearized} along the periodic solution $\gamma(t)$
and introduce the associated evolution operator $\Pi(t,0): H^s(S^1)
\longrightarrow H^s(S^1)$, defined by $\Pi(t,0)\varphi_0= \varphi(t)$ where
$\varphi(t)$ is the solution of the linearized equation
\eqref{linearized}. We recall that the operator
$\Pi(p,0)=D_u(S_{f}(p,0)\gamma(0))$ is called the period
map and we denote $(\mu_i)$ its eigenvalues (the spectral properties
of $\Pi(p,0)$ have been given in Proposition \ref{spectre-Pi}). Since
$\gamma(x,t)$ is a hyperbolic periodic solution, the intersection of
the spectrum of $\Pi(p,0)$ with the unit circle of $\Com$ reduces to
the eigenvalue $1$, which a simple (isolated) eigenvalue. We remark
that, if $\gamma(a)$, $a \in [0,p)$, is another point of the periodic
orbit, the spectrum of $D_u(S_{f}(p,0)\gamma(a))$ coincides with the
one of $\Pi(p,0)$ whereas the corresponding eigenfunctions depend on
the point $\gamma(a)$.

We  denote  $P_u(a)$ (resp. $P_c(a)$, resp. $P_s(a)$) the projection
in $H^s(S^1)$ onto the space
generated by the (generalized) eigenfunctions of $D_u(S_{f}(p,0)\gamma(a))$
corresponding to the eigenvalues with
modulus strictly larger than $1$ (resp. equal to $1$, resp. with
modulus strictly smaller than $1$).

Since a hyperbolic periodic orbit is a particular case of a normally
hyperbolic $C^1$ manifold, we may apply, for example, the existence
results of  \cite{HPS70}, \cite{HPS77} or
\cite[Theorem 14.2 and Remark 14.3]{Ruelle} (see also
\cite{HaRa09}) and thus, we may state the following theorem.

\begin{theorem}\label{th-Wuper}
Let $\Gamma=\{\gamma(t),~t\in\Rm\}$ be a hyperbolic periodic orbit of
Eq.  \eqref{eq}. \HB
1) There exists a small
neighbourhood $U_{\Gamma}$ of $\Gamma$ in $H^s(S^1)$ such that the
local stable and unstable sets
\begin{equation*}
\begin{split}
&W^s_{loc}(\Gamma) \equiv W^s(\Gamma,U_{\Gamma}) = \{u_0 \in H^s(S^1)
\, | \, S_f(t)u_0
\in U_{\Gamma}\, , \, \forall t \geq 0\} \cr
&W^u_{loc}(\Gamma) \equiv W^u(\Gamma,U_{\Gamma}) = \{u_0 \in H^s(S^1)
\, | \, S_f(t)u_0
\in U_{\Gamma}\, , \, \forall t \leq 0\} \cr
\end{split}
\end{equation*}
are (embedded) $C^1$-submanifolds of $H^s(S^1)$ of codimension $i(\Gamma)$ and
dimension $i(\Gamma) +1$ respectively.  \HB
2) Moreover,
$W^{s}_{loc}(\Gamma)$ and $W^u_{loc}(\Gamma)$ are fibrated by the local
strong stable (resp. unstable) manifolds at each point $\gamma(a)
\in \Gamma$, that is,
$$
W^{s}_{loc}(\Gamma) = \cup_{a \in [0,p)} W^{ss}_{loc}(\gamma(a))~, \quad
W^{u}_{loc}(\Gamma) = \cup_{a \in [0,p)} W^{su}_{loc}(\gamma(a))~,
$$
where there exist positive constants ${\tilde r}$, $\kappa$ and $\kappa^*$ such
that
\begin{equation}
\label{Dperaux1}
\begin{split}
W^{ss}_{loc}(\gamma(a))= &\{u_0 \in H^s(S^1) \, | \, \|S_f(t) u_0 -
\gamma(a+t)\|_{H^s} < {\tilde r}~, ~\forall t \geq 0 \, ,\cr
& ~\lim_{t \to \infty} e^{\kappa t} \|S_f(t) u_0 - \gamma(a+t) \|_{H^s}
=0\} ~, \cr
W^{su}_{loc}(\gamma(a))=& \{u_0 \in H^s(S^1) \, | \, \|S_f(t) u_0 -
\gamma(a +t)\|_{H^s} < {\tilde r} ~, ~ \forall t \leq 0 \, , \cr
&~\lim_{t \to  -\infty} e^{\kappa^* t} \|S_f(t) u_0 - \gamma(a+t) \|_{H^s}
=0\}~.
\end{split}
\end{equation}
For any $a \in [0, p)$, $W^{ss}_{loc}(\gamma(a))$ (resp.
$W^{su}_{loc}(\gamma(a))$) is a $C^1$-submanifold of $H^s(S^1)$
tangent at $\gamma(a)$ to $P_s(a) H^s(S^1)$ (resp. $P_u(a) H^s(S^1)$).
\end{theorem}

In the introduction, we have also defined the global stable and
unstable sets as follows
\begin{equation*}
\begin{split}
&W^s(\Gamma) = \{u_0 \in H^s(S^1) \, | \, S_f(t)u_0
\mathrel {\mathop\rightarrow_{t \rightarrow
   +\infty}}  \Gamma\} ~, \cr
& W^u(\Gamma) = \{u_0 \in H^s(S^1) \, | \, S_f(t)u_0
\hbox{ is well-defined for }t\leq 0 \hbox{ and }
S_f(t)u_0 \mathrel {\mathop\rightarrow_{t \rightarrow
   -\infty}} \Gamma\} ~.
\end{split}
\end{equation*}
We recall that $W^s(\Gamma)$ and $W^u(\Gamma)$ are injectively
immersed $C^1$-manifolds of codimension $i(\Gamma)$ and dimension
$i(\Gamma) +1$ respectively. Moreover,
$$
W^u(\Gamma) = \cup_{t \geq 0} S_f(t) W^u_{loc}(\Gamma) ~,
$$
is a union of embedded $C^1$-submanifolds of $H^s(S^1)$ of dimension
$i(\Gamma)+1$ (see \cite{CR} or \cite{HJR} for example).

We are now ready to prove the following corollary.

\begin{coro}\label{CappD02} Let $\Gamma=\{\gamma(t),~t\in\Rm\}$ be a
hyperbolic periodic orbit of Eq.  \eqref{eq}.  \HB Let $u(t)$
be a trajectory of \eqref{eq} belonging to the strong unstable
manifold $W^{su}(\gamma(a)) \setminus \gamma(a)$ (resp.  the local
strong stable manifold $W^{ss}_{loc}(\gamma(a)) \setminus \gamma(a)$)
and let $v(t)=u(t) -\gamma(t+a)$.  \HB
Then, there exists an eigenvalue $\mu_i$ of $\Pi(p,0)$ such that
$|\mu_i|>1$ (resp.  $|\mu_i|<1$) and
$$
\lim_{n \to -\infty}\|v(np)\|_{H^s}^{1/n} = |\mu_i |~, \hbox{ (resp.
} \lim_{n \to \infty}
\|v(np)\|_{H^s}^{1/n} = |\mu_i | \hbox{ ) }
~.
$$
More precisely, the asymptotic behavior of $v(np)$ in $H^s(S^1)$ is
given by one of the following possibilities: \HB
(i) if $\mu_i$ is a simple real eigenvalue with corresponding real
eigenfunction $\varphi_i(a) \in H^s(S^1)$, then there exists $b \in \Rm
- \{0\}$ such that
    $v(np)= b\mu_i^n \varphi_i(a) + o(|\mu_i|^n)$. \HB
(ii) If $\mu_i=\mu_{i+1}$ is a double real eigenvalue with two independent
eigenfunctions $\varphi_i(a)$ and $\varphi_{i+1}(a)$, then there exist
$(b,c)\in\Rm^2-\{(0,0)\}$ such that $v(np) = b\mu_i ^n \varphi_i(a) +
c\mu_i^n \varphi_{i+1}(a)+ o(|\mu_i|^n)$. \HB
(iii) If $\mu_i=\mu_{i+1}$ is an algebraically double real eigenvalue
with eigenfunction $\varphi_i(a)$ and generalized eigenfunction
$\varphi_{i+1}(a)$, then there exist $(b,c)\in\Rm^2-\{(0,0)\}$ such that
$v(t)=(b+cn) \mu_i ^n \varphi_i(a) + c \mu_i ^n \varphi_{i+1}(a)
+  o(|\mu_i|^n)$. \HB
(iv) If $\mu_i=|\mu_i| e^{i\theta}$ is a (simple) complex eigenvalue
with eigenfunction $\varphi_i(a)= \overline{\varphi_{i+1}}(a)$, then
there exist $(b,c)\in\Rm^2-\{(0,0)\}$ such that
$v(np)= |\mu_i|^n  \big[(b\cos(n\theta) - c \sin(
n\theta)) \Re(\varphi_i(a))  - (b\sin( n\theta) + c\cos(n\theta))
\Im( \varphi_i(a) ) \big] +o(|\mu_i|^n)$.

Let $j\in\Nm \setminus \{0\}$ be such that $\mu_i$ belongs to the pair of
eigenvalues $(\mu_{2j-1},\mu_{2j})$, or let $j=0$ if $\mu_i=\mu_0$.  As
a consequence of the asymptotic behaviour, there exists $t_0\in \Rm$
such that, for all $t \leq t_0$ (resp.  $t\geq t_0$), $v(t)$ has
exactly $2j$ zeros which all are
simple.
\end{coro}

\begin{demo} Since $u(t)$ belongs to the local strong stable manifold (or
strong unstable ) manifold of a point $\gamma(a)$ of the non trivial
periodic solution $\gamma(t)$ of \eqref{eq}, the
proof of this corollary is very similar to the
one of Corollary \ref{CappD0}. Thus, we will not repeat the whole
proof, but only give the details of the beginning of the proof, in
order to point out the differences with the proof of Corollary
\ref{CappD0} and also to emphasize the properties of the strong
stable or unstable
manifolds of $\gamma(a)$, that we are using here. \HB
    Since the proofs are very similar when $n$ tends to $\pm
\infty$, we only consider the case where $u(t)$ belongs to the
local strong stable manifold $W^{ss}_{loc}(\gamma(a)) \setminus
\gamma(a)$.  To prove the corollary, we shall apply Theorem \ref{C6},
so we have to check that the hypotheses of Theorem \ref{C6} are
satisfied.  Also, without loss of generality, we may assume that
$a=0$.\HB
Since $u(t)$ belongs to $W^{ss}_{loc}(\gamma(0)) \setminus
\gamma(0)$, there exist two positive constants $c_1$ and $\kappa$
such that
   \begin{equation}
\label{Dconvp1}
\| u(t) -\gamma(t) \|_{H^s} \leq c_1 e^{-\kappa t}~, \quad \hbox{ as } t
\rightarrow +\infty ~.
\end{equation}
The function $v(t) \equiv u(t) -\gamma(t)$ is a classical solution of the
equation
   \begin{equation}
\label{Deqvp}
v_t = v_{xx} + D_{u_x}f(x,\gamma (x,t),\gamma_x(x,t))v_x +
D_uf(x,\gamma(x,t),\gamma_x(x,t))v + a(x,t)v_x +
b(x,t)v~,
\end{equation}
where
   \begin{equation*}
\begin{split}
a(x,t) =& \int_{0}^{1} \big( f'_{u_x}(x, \gamma+\theta
(u-\gamma),\gamma_x+\theta
(u_x-\gamma_x)) -f'_{u_x}(x, \gamma,\gamma_x) \big) d\theta \cr
b(x,t)= & \int_{0}^{1} \big( f'_{u}(x, \gamma+\theta
(u-\gamma),e\gamma_x+\theta
(u_x-\gamma_x)) -f'_{u}(x, \gamma,\gamma_x) d\theta \big)
\end{split}
\end{equation*}
One at once checks that $\| a(x,t)\|_{C^0} + \| b(x,t)\|_{C^0}  \leq
c_2 \| u(t) -\gamma(t) \|_{H^s}$, which implies that  $\|
a(x,t)\|_{L^2} + \| b(x,t)\|_{L^2}$ satisfies the inequality \eqref{Dconv2}.
We next remark that
$$
v((n+1)p)= Tv(pn) + \Sigma (n) v(pn)~,
$$
where $T= \Pi(p,0)$, and
$$
\Sigma (n)v(pn) = \int_0^p \Pi(p,\sigma)
(a(x, np+\sigma)v_x(x, np+\sigma)+ b(x,np+\sigma)v(x,np+\sigma))
   d\sigma~.
   $$
   Arguing as in the proof of Corollary \ref{CappD0}, one checks that
   $\Sigma(n)$ satisfies the condition \eqref{CSigrho}. As the periodic
   orbit $\Gamma$ is hyperbolic, $T= \Pi(p,0)$ admits a shifted
   exponential dichotomy on $\Zr^+$ (see \cite{HaLi86}
for example). Thus, by Proposition \ref{C2},
   the family $L(\cdot)= T
+\Sigma(\cdot)$ admits a shifted dichotomy on $\Zr^+$.
As, by \cite{Agmon}, every non-zero
solution of a linear parabolic equation does not go to zero faster
than an exponential when $t$ tends to infinity, the hypotheses of
Theorem \ref{C5} hold. Since the exponential
decay property \eqref{Dconvp1} holds, it follows from
Theorem \ref{C5}, that there exists $\rho_i \in {\cal R}(T)$, $\rho_i
<1$, (and thus $\mu_i$
belonging to the spectrum of $\Pi(p,0)$), such that
$$
\lim_{n \to \infty} \| v(n)\|_{H^s}^{1/n}=|\rho_i |= | \mu_i| ~.
$$
Corollary \ref{CappD02} is now an easy consequence of Theorem
\ref{C6} and of Proposition \ref{spectre-Pi} (for the details, see
the proof of Corollary \ref{CappD0}).
\end{demo}

\vskip 2mm

\begin{coro}\label{CappD2} Let $\Gamma=\{\gamma(t),~t\in\Rm\}$ be a
hyperbolic periodic orbit of Eq.  \eqref{eq} and $u(t)$ be a
trajectory of \eqref{eq} belonging to the strong unstable manifold
$W^{su}(\gamma(a)) \setminus \gamma(a)$ (resp.  the local strong stable
manifold $W^{s}_{loc}(\gamma(a)) \setminus \gamma(a)$).  Let $v_0 \in
T_{u(0)}W^{su}(\gamma(a))$ (resp.  $v_0 \in
T_{u(0)}W^{ss}_{loc}(\gamma(a))$) and $v(t)$ be the solution for $t\leq 0$
(resp. $t\geq 0$) of the
linearized equation \eqref{eq-AppD}.  Then, there exists an eigenvalue
$\mu_i$ of $\Pi(p,0)$ such that $|\mu_i|>1$ (resp.  $|\mu_i|<1$) and
$$
\lim_{n \to -\infty}\|v(np)\|_{H^s}^{1/n} = |\mu_i |~, \hbox{ (resp.
} \lim_{n \to \infty}
\|v(np)\|_{H^s}^{1/n} = |\mu_i | \hbox{ ) }
~.
$$
Moreover, all the possible asymptotic behaviors are the same as those
described in Corollary \ref{CappD02}. \HB
In addition, let $j\in\Nm \setminus \{0\}$ be such that $\mu_i$ 
belongs to the pair of
eigenvalues $(\mu_{2j-1},\mu_{2j})$, or let $j=0$ if $\mu_i=\mu_0$.
Then, there exists $t_0\in \Rm$
such that, for all $t \leq t_0$ (resp.  $t\geq t_0$), $v(t)$ has
exactly $2j$ zeros which all are
simple.
\end{coro}

\begin{demo} If $u(t)$ belongs to the local strong stable manifold
   $W^{ss}_{loc}(\gamma(a)) \setminus \gamma(a)$, $u(t)$ satisfies the
   property \eqref{Dconvp1} and
   Eq.  \eqref{eq-AppD} can be written in the form \eqref{Deqvp},
where the functions $a(x,t)$ and $b(x,t)$ satisfy the properties
   \eqref{Dconv2}. We emphasize that, by \cite[Theorem C6]{ChChHa}, we
already know
   that $\limsup_{n \to \infty} \| v(np)\|_{H^s}^{1/n} <1$.
We thus obtain the asymptotic behavior of $v(t)$ by
following the lines of the proof of Corollary \ref{CappD02}.
\end{demo}


\end{document}